\documentclass[11pt]{article}
\usepackage{amssymb}
\usepackage{amssymb}
\usepackage{amssymb}
\usepackage{amssymb}

\usepackage{}
\usepackage{latexsym,amsmath}
\usepackage{amsfonts}
\usepackage{amsmath}
\usepackage{amssymb}
\usepackage{latexsym}
\usepackage{indentfirst}
\usepackage{fancyhdr}
\usepackage{stmaryrd}
\usepackage{mathrsfs}
\usepackage{cite}
\usepackage{amsmath,amsfonts,amssymb}
\usepackage{stmaryrd}
\usepackage[all]{xy}

\date{}

\numberwithin{equation}{section}

\newtheorem{theorem}{Theorem}[section]
\newtheorem{lemma}[theorem]{Lemma}

\newcommand{\al}{\alpha}
\newcommand{\ga}{\gamma}

\newcommand{\be}{\beta}
\newcommand{\G}{\Gamma}
\newcommand{\vs}{\varsigma}

\newcommand{\lmd}{\lambda}

\newcommand{\ves}{\varepsilon}
\newcommand{\kn}{\mbox{ker}\:}
\newcommand{\ad}{\mbox{ad}\:}
\newcommand{\for}{\mbox{for}}
\newcommand{\mbb}{\mathbb}
\newcommand{\es}{\epsilon}

\makeatletter
\newcommand{\rmnum}[1]{\romannumeral#1}
\newcommand{\Rmnum}[1]{\expandafter\@slowromancap\romannumeral#1@}
\makeatother

\begin{document}
\newcommand{\parti}{\partial}

\baselineskip=16pt

\title{Multiplicity-Free Representations of \\ Divergence-Free Lie Algebras }
\author{Ling \uppercase{Chen}
\thanks{ E-mail: chenling@amss.ac.cn (L. Chen)}\\
{\small{ Institute of Applied Mathematics, Academy of Mathematics and Systems Science,}}\\
{\small{Chinese Academy of Sciences,
    Beijing
    100190, P. R. China}}}\maketitle

\begin{abstract}

Divergence-free Lie algebras (also known as the special Lie algebras
of Cartan type) are Lie algebras of volume-preserving transformation
groups. They are simple in generic case.  Dokovic and Zhao found a
certain graded generalization of them. In this paper, we classify
all the irreducible and indecomposable multiplicity-free modules of
the simple generalized divergence-free Lie algebras.

\end{abstract}

\section{Introduction}

Lie algebras of Cartan type are closely related to geometry and
dynamics. They also play important roles in the structure theory of
simple Lie algebras.  Kac \cite{Ka} gave an abstract definition of
generalized Cartan type Lie algebras by derivations. After that,
various graded generalizations of Lie algebras of Cartan type were
constructed and studied by Kawamoto \cite{Kaw}, Osborn \cite{O},
Dokovic and Zhao \cite{DZ1,DZ2,DZ3}, and Zhao \cite{Zk1}. The
fundamental ingredients for these graded generalizations are the
pairs of a semi-group algebra and a set of grading operators and
down-grading operators. Motivated from his works on quadratic
conformal algebras \cite{X3} and Hamiltonian operators \cite{X4}, Xu
\cite{X2} found certain nongraded generalizations of the Lie
algebras of Cartan type and determined their simplicity. Xu's
nongraded generalizations are based on a semi-group algebra and a
set of locally-finite derivations. In generic case, the algebras
found by Xu do not have a toral Cartan subalgebra, which makes key
difference from the other generalizations.

 The representations of simple Lie algebras of Cartan
type have attracted the attention of many researchers as well. Shen
\cite{Sg1, Sg2, Sg3} introduced mixed product of graded modules over
graded Lie algebras of Cartan type and obtained certain irreducible
modules over a field with characteristic $p$. Larsson [8] obtained
the same representation of Witt algebras from different motivation
(Later was called ``Larsson functor" by some authors). Rao \cite{R1,
R2} constructed some irreducible weight modules over the derivation
Lie algebra of the algebra of Laurent polynomials virtually based on
Shen's mixed product (Larsson functor).

Lin and Tan [9] constructed some irreducible weight modules over the
Lie algebra of quantum torus via Shen's mixed product. Y. Zhao
\cite{Zy1} determined the module structure of Shen's mixed product
over Xu's nongraded Lie algebras of Witt type (cf. \cite{X2}).
Moreover, she \cite{Zy2} constructed a family of irreducible modules
for Xu's two-derivation nongraded  Lie algebras of Block type (cf.
\cite{X1}) via an embedding of the Lie algebras into the
corresponding generalized Weyl algebra and Shen's mixed product.
Furthermore, she \cite{Zy3} obtained a composition series for a
family of modules with parameters over Xu's nongraded Hamiltonian
Lie algebras (cf. \cite{X2}). These modules are constructed from
finite-dimensional multiplicity-free irreducible modules of
symplectic Lie algebras by Shen's mixed product. In \cite{ZL}, Y.
Zhao and Liang generalized the results of \cite{Zy2} to Xu's
four-derivation nongraded  Lie algebras of Block type (cf.
\cite{X1}).

Howe \cite{Ho} classified all the finite-dimensional
multiplicity-free representations over finite-dimensional simple Lie
algebras. On infinite-dimensional simple Lie
algebra side, Kaplansky \cite{Kap,KS} and Santharoubane \cite{KS} gave classifications of multiplicity-free
representations over classical Virasoro algebras. After that, Su \cite{S1,S2} generalized Kaplansky and Santharoubane's result
to multiplicity-free modules over high rank Virosoro algebras and super-Virasoro algebras. Based on their classifications of multiplicity-free
representations over generalized Virasoro algebras in \cite{SZk},
Zhao \cite{Zk2} classified the multiplicity-free representations
over graded generalized Witt algebras. Moreover, Su and Zhou
\cite{SZh} generalized Zhao's result in \cite{Zk2} to generalized
weight modules over the nongraded generalized Witt algebras
introduced by Xu \cite{X2}.

Divergence-free Lie algebras, which are also called the {\it special
Lie algebras of Cartan type},  are Lie algebras of volume-preserving
transformation groups. They are simple in generic case. In this
paper,  we classify all the irreducible and indecomposable
multiplicity-free graded modules over the simple generalized
divergence-free Lie algebras introduced by Dokovic and Zhao
\cite{DZ3}. Since the algebras we are concerned with do not contain
a centerless Virasoro algebra, we can not use the existing results
and have to find a new way to do our classification. This is a
reason why this paper is so long. Below are some details.

Let $\G$ be an additive group. A Lie algebra ${\cal G}$ over a field
$\mathbb{F}$ is called a $\G$-{\it graded Lie algebra} if ${\cal G}=
\bigoplus_{\al\in\G} {\cal G}_\al$, where each ${\cal G}_\al$ is a
subspace of ${\cal G}$, and $[{\cal G}_\al, {\cal G}_\be]\subseteq
{\cal G}_{\al+\be}$. A module $M$ over a $\G$-graded Lie algebra
${\cal G}$ is called a {\it graded module} if $M= \bigoplus _{\al
\in \G} M_\al$ is a $\G$-graded space and ${\cal G}_\al M_\be
\subset M_{\al+\be}$. Each component $M_\al$ is called a {\it
homogenous subspace} of $M$. The aim of this paper is to give a
complete classification of the irreducible and indecomposable
 graded modules of the simple generalized
divergence-free Lie algebras introduced by Dokovic and Zhao
\cite{DZ3}, whose dimensions of homogenous subspaces are less than
or equal to one. In a subsequent paper, we will classify all the
irreducible and indecomposable multiplicity-free generalized graded
modules over Xu's nongraded divergence-free Lie algebras (cf.
\cite{X2}, \cite{SX}) as Su and Zhou did for Xu's nongraded Witt
algebras in \cite{SZh}.

Throughout this paper, we assume that the base field $\mathbb{F}$ is an
algebraically closed field with characteristic 0. This paper is organized as follows. In Section 2, we review the
definitions and some basic facts about the generalized Witt algebras given by Kawamoto \cite{Kaw}
and the generalized divergence-free Lie algebras introduced by
Dokovic and Zhao \cite{DZ1,DZ2}. Moreover, we construct three
classes of graded modules over Dokovic and Zhao's generalized
divergence-free Lie algebras with 1-dimensional homogenous
subspaces, and give the necessary and sufficient conditions for two
such modules to be isomorphic. Furthermore, we give the main theorem
of classifying all the irreducible and indecomposable graded modules
of these divergence-free Lie algebras, whose dimensions of
homogenous subspaces are less than or equal to one. In Sections 3 to
5, we prove the main theorem under different conditions,
respectively.

\section{The Lie Algebras and Modules}

We first review the definition of generalized Witt algebras introduced by Kawamoto \cite{Kaw}. The definition here is different in form to that in \cite{Kaw}, but they are equivalent.

Let $D$ be a finite-dimensional nonzero vector space and let $\G$ be
an additive subgroup of its dual space $D^\ast$ such that
$\bigcap_{\al\in\G}\kn\al=\{0\}$. Denote by $\mathbb{F}[\G]$ the
vector space with a basis $\{x^\al\mid\al\in\G\}$ and the
multiplication determined by $x^\al x^\be=x^{\al+\be}$ for
$\al,\be\in\G$. Write 1 instead of $x^{0}$ for convenience. The
tensor product ${\cal W}=\mathbb{F}[\G]\otimes_{\mathbb{F}}D $ is a
free left $\mathbb{F}[\G]$-module. Moreover, we denote an arbitrary
element of $D$ by $\parti$, and for simplicity,  we write
$x^\al\parti$ instead of $x^\al\otimes \parti$.
 Under the bracket
\begin{equation}
[x^\al\parti_{1},
x^\be\parti_{2}]=x^{\al+\be}(\be(\parti_{1})\parti_{2}-\al(\parti_{2})\parti_{1})
\textrm{ for any }
\parti_{1}, \parti_{2} \in D \textrm{ and } \al,\be \in \G,
\end{equation}
${\cal W}$ becomes a Lie algebra. It is usually denoted by ${\cal
W}(\G, D)$, and is called a {\it generalized Witt algebra}. Kawamoto \cite{Kaw} showed that ${\cal W}(\G, D)$ is a
simple Lie algebra.

The Lie algebra ${\cal W}={\cal W}(\G, D)$ has a natural
$\G$-gradation ${\cal W}=\bigoplus _{\al\in\G} {\cal W}_\al$, where
${\cal W}_\al=x^\al D$ for $\al \in \G$. In particular, we have
${\cal W}_{0}=D$. It follows that
\begin{equation}
[\parti, x^\al\parti_{1}]=\al(\parti) x^\al\parti_{1}, \quad
\forall\;
\parti, \parti_{1} \in D, \, \al\in \G,
\end{equation}
namely, $\ad\parti$ acts on ${\cal W}_\al$ as the scalar
$\al(\parti)$. Hence $D$ is a toral Cartan subalgebra of ${\cal W}$.

We now review the generalized divergence-free Lie algebras
introduced by Dokovic and Zhao \cite{DZ3}. Define the divergence as
the $\mathbb{F}$-linear map div from ${\cal W}$ to $\mathbb{F}[\G]$
determined by
\begin{equation}
\textrm{div}(x^\al\parti)=\al(\parti) x^\al, \quad \forall \; \al\in
\G, \, \parti \in D.
\end{equation}
Note that $\mathbb{F}[\G]$ becomes a ${\cal W}$-module with the
action:
\begin{equation}x^\al\parti.\:x^\be=\be(\parti)x^{\al+\be}\qquad\mbox{for}\;\parti\in
D, \,\al,\be\in\G.
\end{equation}
It is easy to verify that the
divergence has the following two properties:
\begin{equation}\label{equ:2.1}
\textrm{div}[u,v]=u. \textrm{div}(v)-v. \textrm{div}(u),
\end{equation}
\begin{equation}
\textrm{div}(fw)=f\textrm{div}(w)+ w. f
\end{equation}
for $u,v,w \in {\cal W}$ and $f \in \mathbb{F}[\G]$. So the subspace
\begin{equation}
\tilde{\cal S}= \textrm{ker (div)}
\end{equation}
forms a subalgebra of ${\cal W}$. If $\dim D=1$, then $\tilde{\cal
S}=D$. We assume from now on that $\dim D \geq 2$. Then we have that
\begin{equation}
\tilde{\cal S}=\bigoplus_{\al\in \G} \tilde{\cal S}_\al, \textrm{
where } \tilde{\cal S}_\al=\tilde{\cal S}\cap {\cal W}_\al.
\end{equation}
It is clear that $ \tilde{S}_\al=x^\al(\kn\al)$. Consequently
\begin{equation}
\textrm{codim}_{{\cal W}_\al}(\tilde{\cal S}_\al)=1-\delta_{\al,0}.
\end{equation}
We set ${\cal S}=[\tilde{\cal S},\tilde{\cal S}]$. Then
\begin{equation}
{\cal S}=\bigoplus_{\al\in \G \backslash \{0\}} \tilde{\cal S}_\al
\end{equation}
and ${\cal S}$ is a simple Lie algebra (c.f. \cite{DZ3}). The Lie
algebra ${\cal S}$ is called a {\it simple generalized
divergence-free Lie algebra} (also called a {\it generalized Cartan
type special Lie algebra}). The goal of this paper is to classify
all the irreducible and indecomposable multiplicity-free graded
modules of the Lie algebra ${\cal S}$. We sometimes denote ${\cal
S}$ by ${\cal S}(\G, D)$.

For any $\mu,  \eta \in D^{\ast}$, $\zeta  \in \G$, we define the
$\G$-graded ${\cal S}$-module $\mathscr{M}_{\mu}=\bigoplus_{\al\in \G} \mathbb{F} v_\al$ with the action
\begin{equation}\label{equ:3.1}
x^\al\parti . v_\be=(\be+\mu)(\parti)v_{\al+\be}  \textrm{ for any }
\al \in \G\backslash\{0\},  \be\in \G,  \parti \in \kn\al;
\end{equation}
the $\G$-graded ${\cal S}$-module $\mathscr{A}_{\zeta ,\eta }=\bigoplus_{\al\in \G} \mathbb{F} v_\al$ with the action
\begin{eqnarray}
& &x^\al\parti. v_\be=  (\be+\zeta)(\parti)v_{\al+\be}  \textrm{ for
} \al\in \G \backslash \{0\},  \be\in
\G\backslash\{-\zeta\},  \parti \in \kn\al, \nonumber \\
& & x^\al\parti . v_{-\zeta}  = \eta(\parti)v_{\al-\zeta}  \textrm{
for } \al\in \G\backslash \{0\},   \parti \in
\kn\al;\label{equ:3.2}
\end{eqnarray}
the $\G$-graded ${\cal S}$-module $\mathscr{B}_{\zeta ,\eta }=\bigoplus_{\al\in\G} \mathbb{F} v_\al$ with the action
\begin{eqnarray}
& &x^\al\parti . v_\be = (\be+\zeta)(\parti)v_{\al+\be} \textrm{ for
} \al\in \G \backslash \{0\},  \be\in
\G\backslash\{ -\al-\zeta\},   \parti \in \kn\al, \nonumber \\
& & x^\al\parti . v_{-\al-\zeta}  = \eta(\parti)v_{-\zeta}  \textrm{
for } \al\in \G \backslash \{0\},  \parti \in
\kn\al.\label{equ:3.3}
\end{eqnarray}

It is easy to verify that these three classes of modules are
restrictions to ${\cal S}$ of Zhao's five classes of modules over
${\cal W}$ in \cite{Zk2}. It is straightforward to prove:

\begin{lemma}\label{le:3.1}
The irreducibility of the three classes of modules
$\mathscr{M}_{\mu}$, $\mathscr{A}_{\zeta ,\eta }$ and $\mathscr{B}_{\zeta ,\eta }$ is as follows:

\rmnum{1}) $\mathscr{M}_{\mu}$ is irreducible if  $\mu
\not\in \G$. When $\mu \in \G$, $\mathscr{M}_{\mu}=
(\bigoplus_{\al\in \G \backslash \{-\mu \}} \mathbb{F} v_\al) \bigoplus
\mathbb{F} v_{-\mu}$ is a direct sum of two irreducible submodules.

\rmnum{2}) If $\eta\not=0$, the modules $\mathscr{A}_{\zeta ,\eta }$
and $\mathscr{B}_{\zeta ,\eta }$ are indecomposable but reducible.
In this case, $\mathscr{A}_{\zeta ,\eta }$ has irreducible submodule
$\bigoplus_{\al\in\G \backslash \{-\zeta\}} \mathbb{F} v_\al$, while
$\mathscr{B}_{\zeta ,\eta }$ has one-dimensional trivial submodule
$\mathbb{F} v_{-\zeta}$ and the quotient $\mathscr{B}_{\zeta ,\eta
}/\mathbb{F} v_{-\zeta}$ is irreducible. If $\eta=0$, then
$\mathscr{A}_{\zeta,0}\simeq \mathscr{B}_{\zeta,0}\simeq
\mathscr{M}_{0}$.
\end{lemma}

 We use $\mathscr{M}_{\mu}'$,
$\mathscr{A}_{\zeta ,\eta}'$ and $\mathscr{B}_{\zeta ,\eta}'$ to
denote the nontrivial irreducible submodules or the  nontrivial
irreducible quotients of $\mathscr{M}_{\mu}$, $\mathscr{A}_{\zeta
,\eta}$ and $\mathscr{B}_{\zeta ,\eta}$, respectively. In analogy with
Zhao's result, we have:
\begin{theorem}\label{th:3.2}
Among the ${\cal S}$-modules $\mathscr{M}_{\mu}$,
$\mathscr{A}_{\zeta ,\eta}$, $\mathscr{B}_{\zeta
,\eta}$ for $\mu,  \eta \in D^{\ast}, \ \zeta  \in \G$,
and their nontrivial irreducible submodules or  nontrivial
irreducible quotients, we have only the following module
isomorphisms:
\begin{eqnarray*}
&(\rmnum{1})  &  \mathscr{M}_{\mu} \simeq \mathscr{M}_{\mu'} \textrm{ iff } \mu - \mu' \in \G, \\
&(\rmnum{2})  &  \mathscr{M}_{\mu}' \simeq \mathscr{M}_{\mu'}' \textrm{ iff } \mu - \mu' \in \G,\\
&(\rmnum{3})  &  \mathscr{A}_{\zeta ,\eta}\simeq \mathscr{A}_{0 ,\eta}\simeq
\mathscr{A}_{0,a\eta}  \textrm{ for any } a \in \mathbb{F}\backslash \{0\}, \\
&(\rmnum{4}) &  \mathscr{B}_{\zeta ,\eta}\simeq \mathscr{B}_{0 ,\eta}\simeq
\mathscr{B}_{0,a\eta} \textrm{ for any } a \in \mathbb{F}\backslash \{0\}, \\
&(\rmnum{5})  &  \mathscr{A}_{0,0}\simeq \mathscr{B}_{0,0}\simeq \mathscr{M}_{0},\\
&(\rmnum{6})  &  \mathscr{A}_{\zeta ,\eta}'\simeq \mathscr{B}_{\zeta ,\eta}'\simeq
\mathscr{M}_{0}'.\\
\end{eqnarray*}
\end{theorem}

Then we give our main result of this paper:

\begin{theorem}\label{th:4.1} Suppose $\dim
D\geq 3$. Assume that $V = \bigoplus_{\theta \in \G} V_\theta$ is a
$\G$-graded ${\cal S}(\G, D)$-module with $\dim V_\theta\leq 1$ for
$\theta \in \G$. If $V$ is irreducible or indecomposable, then $V$ is isomorphic to one of the following modules
for appropriate $\mu\in D^{\ast}, \, \eta \in D^{\ast}\backslash\{0\}$ and $\zeta \in A $:

$(\rmnum{1})  \textrm{the trivial module } \mathbb{F}v_{0};\
(\rmnum{2})    \mathscr{M}_{\mu}';\  (\rmnum{3})   \mathscr{A}_{\zeta
,\eta};\ (\rmnum{4}) \mathscr{B}_{\zeta ,\eta}.$

\end{theorem}

We will prove it case by case progressively in the following sections.

\section{The case $\dim D =3$ and $\G \simeq  \mathbb{Z}^{3}$}

In this section, we will prove Theorem \ref{th:4.1} under the condition that $\dim D =3$ and $\G \simeq  \mathbb{Z}^{3}$. Throughout this section, we shall always assume $\dim D =3$ and $\G \simeq  \mathbb{Z}^{3}$.

Let $M = \bigoplus_{\theta \in \G} M_{\theta}$ be a $\G$-graded
${\cal S}$-module with $\dim M_{\theta }= 1$ for each $\theta\in\G$.
Write $ M_{\theta}=\mbb{F}w_\theta$ for $\theta \in \G$. In order to prove Theorem \ref{th:4.1}, we first specify all the possible action of ${\cal S}(\G, D)$ on $M$. We start our analysis with two useful properties.\vspace{0.2cm}

Fix any nonzero $\sigma\in\G$. Take any two nonzero vectors
$\parti,\parti' \in \kn\sigma$. It can be noticed that $x^{\sigma}\parti$ and $x^{-\sigma}\parti'$ preserve $\bigoplus_{i\in\mbb{Z}}M_{\nu + i \sigma}$ for any $\nu \in \G$. Fix $\nu $ and we write
\begin{equation}
x^{-\sigma}\parti'.x^{\sigma}\parti.w_{\nu+ i
\sigma} =c_{i}w_{\nu + i \sigma} \textrm{ with }c_{i}\in\mbb{F} \textrm{ for }i
\in\mathbb{Z}.
 \end{equation}
We demonstrate how $x^{\sigma}\parti$ and $x^{-\sigma}\parti'$ act on $\bigoplus_{i\in\mbb{Z}}M_{\nu + i \sigma}$ in the following lemma.

\begin{lemma}\label{le:5.2}
The constant $c_i=c$ is independent of $i$. Moreover, if $c\not=0$,
for  $a\in\{\sqrt{c},-\sqrt{c}\}$, there exist $\{0\not= v_{\nu + i
\sigma}\in M_{\nu + i \sigma}\mid i\in\mbb{Z}\}$ such that
$x^{\sigma}\parti. v_{\nu+ (i-1)\sigma}=a v_{\nu+i\sigma}$ and
$x^{-\sigma}\parti'.v_{\nu+i\sigma}=a v_{\nu+(i-1)\sigma}$ for
$i\in\mbb{Z}$.\end{lemma}

\noindent{\bf Proof}. Note $x^{\sigma}{\parti}$ and
$x^{-\sigma}{\parti}'$ commute. For any $i\in\mathbb{Z}$, we have
\begin{equation}
c_{i}x^{\sigma}{\parti}.w_{\nu + i
\sigma}=x^{\sigma}{\parti}.(x^{-\sigma}{\parti}'.x^{\sigma}{\parti}.w_{\nu
+i\sigma})=x^{-\sigma}{\parti}'.x^{\sigma}{\parti}.(x^{\sigma}{\parti}.w_{\nu
+ i\sigma})=c_{i+1}x^{\sigma}{\parti}.w_{\nu + i\sigma}.
\end{equation}
Thus  $c_{i}=c_{i+1}$ if $x^{\sigma}{\parti}.w_{\nu + i \sigma}
\not=0$. When $x^{\sigma}{\parti}.w_{\nu + i \sigma} =0$, we obtain
\begin{equation}
x^{-\sigma}\parti'.x^{\sigma}\parti.w_{\nu+ i \sigma} =0
\textrm{ and } x^{-\sigma}\parti'.x^{\sigma}\parti.w_{\nu+ (i+1)
\sigma}=x^{\sigma}\parti.(x^{-\sigma}\parti'.w_{\nu+ (i+1) \sigma})
=0,
\end{equation}
which indicates $c_i=0=c_{i+1}$. To sum up, we
get  $c_i=c_{i+1}$ for any $i\in\mathbb{Z}$, i.e.,
 $c_i=c$ is independent of $i$.

Suppose $c\not=0$. The fact $c_i=c$ implies that $
x^{\sigma}{\parti}.w_{\nu+i\sigma}\not=0$ for all $i \in\mathbb{Z}$.
Take $v_{\nu}=w_{\nu}$. Fixing $a\in\{\sqrt{c},-\sqrt{c}\}$, we
define $v_{\nu+i\sigma}$'s by
$x^{\sigma}{\parti}.v_{\nu+(i-1)\sigma}=a v_{\nu+i\sigma}$ for
$i\in\mathbb{Z}$. Then
\begin{equation}x^{-\sigma}{\parti}'.v_{\nu+i\sigma}=
\frac{1}{a}x^{-\sigma}{\parti}'.(x^{\sigma}{\parti}.v_{\nu+(i-1)\sigma})
=av_{\nu+(i-1)\sigma}.\end{equation}
So the lemma follows.
$\qquad\Box$\vspace{0.2cm}

Let $\sigma, \rho \in \G$ be any $\mathbb{Z}$-linearly independent
elements. Observe that $\dim
(\kn{\sigma} \cap \kn{\rho})=1$. Pick any nonzero vector
$\tilde{\parti}_{1} \in \kn{\sigma} \cap \kn{\rho}$. Take
$\tilde{\parti}_{2}\in\kn\sigma\setminus\mbb{F}\tilde{\parti}_{1}$
and
$\tilde{\parti}_{3}\in\kn\rho\setminus\mbb{F}\tilde{\parti}_{1}$.
Then $\{\tilde{\parti}_{1}, \tilde{\parti}_{2}\}$ forms a basis of
$\kn{\sigma}$ and $\{\tilde{\parti}_{1}, \tilde{\parti}_{3}\}$ forms
a basis of $\kn{\rho}$. Obviously $\rho(\tilde{\parti}_{2})\not=0$ and $\sigma(\tilde{\parti}_{3})\not=0$. Notice that $x^{\pm\sigma}\tilde{\parti}_{1}$, $x^{\pm\sigma}\tilde{\parti}_{2}$, $x^{\pm\rho}\tilde{\parti}_{1}$ and $x^{\pm\rho}\tilde{\parti}_{3}$ preserve $ \bigoplus_{i,k\in\mbb{Z}} M_{\nu + i \sigma + k \rho}$ for any $\nu\in\G$. We derive their action on $ \bigoplus_{i,k\in\mbb{Z}} M_{\nu + i \sigma + k \rho}$ under certain conditions in the following lemma.

\begin{lemma}\label{le:5.3}
 If
$x^{-\sigma}\tilde{\parti}_{1}.x^{\sigma}\tilde{\parti}_{1}.w_{\nu}
\not=0$ and
$x^{-\rho}\tilde{\parti}_{1}.x^{\rho}\tilde{\parti}_{1}.w_{\nu}
\not=0$ for some $\nu\in \G$, then there exist
$\{0\not=v_{\nu+i\sigma+ k \rho}\in M_{\nu + i \sigma + k \rho}\mid i,k\in\mbb{Z}\}$ such that
\begin{equation}
x^{ \pm \sigma}\tilde{\parti}_{1}.v_{\nu+i\sigma+
k \rho}=a_{1} v_{\nu+(i \pm
1)\sigma+ k \rho},\quad x^{\pm \rho}\tilde{\parti}_{1}.v_{\nu+i\sigma+ k
\rho}=b_{1} v_{\nu+i\sigma+(k \pm 1) \rho}
\end{equation}
for $i,k\in\mbb{Z}$, where $a_{1}$ and $b_{1}$ are any nonzero
constants satisfying
$x^{-\sigma}\tilde{\parti}_{1}.x^{\sigma}\tilde{\parti}_{1}.w_{\nu}
= a_{1}^{2}w_{\nu}$ and $
x^{-\rho}\tilde{\parti}_{1}.x^{\rho}\tilde{\parti}_{1}.w_{\nu}=
b_{1}^{2}w_{\nu}$, respectively. For such $v_{\nu+i\sigma+ k
\rho}$'s, writing
$x^{\sigma}\tilde{\parti}_{2}.v_{\nu}=a_{2}v_{\nu+\sigma}$, $x^{
\rho}\tilde{\parti}_{3}.v_{\nu}=a_{3} v_{\nu+ \rho}$ and $
x^{\sigma+\rho}\tilde{\parti}_{1}.v_{\nu}=d v_{\nu+\sigma+\rho}$, we
have the following  relations:
\begin{equation}
x^{\sigma+\rho}\tilde{\parti}_{1}.v_{\nu+i\sigma+ k \rho}=d
v_{\nu+(i+1)\sigma+(k+1) \rho} \quad\textrm{ for }i,k\in\mbb{Z};
\end{equation}
\begin{equation}
x^{ \pm \sigma}\tilde{\parti}_{2}.v_{\nu+i\sigma+
k \rho}=(a_{2}+k \rho(\tilde{\parti}_{2})\frac{d}{b_{1}})
v_{\nu+(i \pm
1)\sigma+ k \rho} \quad\textrm{ for }i,k\in\mbb{Z};
\end{equation}
\begin{equation}
x^{\pm \rho}\tilde{\parti}_{3}.v_{\nu+i\sigma+ k
\rho}=(a_{3}+ i \sigma(\tilde{\parti}_{3})\frac{d}{a_{1}})
v_{\nu+i\sigma+(k \pm 1) \rho} \quad\textrm{ for }i,k\in\mbb{Z}.
\end{equation}
Moreover, we have $a_{1}^{2}=b_{1}^{2}=d^{2}$.
\end{lemma}

\noindent{\bf{Proof.}} Since
$x^{-\sigma}\tilde{\parti}_{1}.x^{\sigma}\tilde{\parti}_{1}.w_{\nu}
\not=0$ and
$x^{-\rho}\tilde{\parti}_{1}.x^{\rho}\tilde{\parti}_{1}.w_{\nu}
\not=0$, Lemma \ref{le:5.2} enables us to choose $\{0\not= v_{\nu+i\sigma}\in
M_{\nu + i \sigma} \mid i \in \mathbb{Z}\}$ and $\{0\not= v_{\nu+ k
\rho}\in M_{\nu + k \rho} \mid k \in \mathbb{Z}\}$  such that
\begin{equation}\label{equ:1.4}
x^{\pm \sigma}\tilde{\parti}_{1}.v_{\nu+i\sigma}=a_{1} v_{\nu+(i \pm
1)\sigma}
\textrm{ and }
x^{\pm \rho}\tilde{\parti}_{1}.v_{\nu+ k \rho}=b_{1} v_{\nu+(k \pm
1) \rho}
\end{equation}
for $i,k\in\mbb{Z}$, where $a_{1}$ and $b_{1}$ are any nonzero
constants satisfying
$x^{-\sigma}\tilde{\parti}_{1}.x^{\sigma}\tilde{\parti}_{1}.w_{\nu}
= a_{1}^{2}w_{\nu}$ and $
x^{-\rho}\tilde{\parti}_{1}.x^{\rho}\tilde{\parti}_{1}.w_{\nu}=
b_{1}^{2}w_{\nu}$, respectively.

In order to determine the other $v$'s, we shall first prove $x^{
\rho}\tilde{\parti}_{1}.w_{\nu + i \sigma+k \rho}\not=0$ for all
$i,k\in\mathbb{Z}$. Consider the case $i,k\geq 0$. We give the proof by induction on $i$. If $i=0$, (\ref{equ:1.4}) shows
\begin{equation}
x^{ \rho}\tilde{\parti}_{1}.w_{\nu +k\rho}\not=0 \ \textrm{ for
all }k\geq 0.
\end{equation}
For some $i\geq 0$, we assume
\begin{equation}\label{e:1}
x^{ \rho}\tilde{\parti}_{1}.w_{\nu + i\sigma+k\rho}\not=0 \textrm{ for
all } k\geq 0.
\end{equation}
Since $x^{- \sigma}\tilde{\parti}_{1}$ and $x^{
\rho}\tilde{\parti}_{1}$ commute, (\ref{equ:1.4}) and (\ref{e:1}) indicate
\begin{equation}\label{equ:1.8}
x^{- \sigma}\tilde{\parti}_{1}.x^{ \rho}\tilde{\parti}_{1}.v_{\nu +(i+1)
\sigma} = x^{ \rho}\tilde{\parti}_{1}.(x^{-
\sigma}\tilde{\parti}_{1}.v_{\nu + (i+1)\sigma})= a_1 x^{ \rho}\tilde{\parti}_{1}.v_{\nu + i\sigma} \not=0,
\end{equation}
which implies $x^{ \rho}\tilde{\parti}_{1}.v_{\nu +
(i+1)\sigma}\not=0$ and $x^{- \sigma}\tilde{\parti}_{1}.w_{\nu +
(i+1)\sigma +\rho}\not=0$. Suppose $x^{
\rho}\tilde{\parti}_{1}.w_{\nu + (i+1)\sigma+k\rho}\not=0$ and $x^{-
\sigma}\tilde{\parti}_{1}.w_{\nu + (i+1)\sigma +(k+1)\rho}\not=0$
for some $k\geq 0$. Then this together with (\ref{e:1}) imply
\begin{equation}
x^{- \sigma}\tilde{\parti}_{1}.x^{ \rho}\tilde{\parti}_{1}.w_{\nu + (i+1)\sigma +(k+1)\rho} = x^{ \rho}\tilde{\parti}_{1}.x^{-
\sigma}\tilde{\parti}_{1}.w_{\nu + (i+1)\sigma +(k+1)\rho}\not=0,
\end{equation}
which implies $x^{ \rho}\tilde{\parti}_{1}.w_{\nu +
(i+1)\sigma+(k+1)\rho}\not=0$ and $x^{-
\sigma}\tilde{\parti}_{1}.w_{\nu + (i+1)\sigma
+(k+2)\rho}\not=0$. By induction on $k$, we have
\begin{equation}\label{equ:7.1}
x^{ \rho}\tilde{\parti}_{1}.w_{\nu + (i+1)\sigma+k\rho}\not=0
\textrm{ and } x^{- \sigma}\tilde{\parti}_{1}.w_{\nu + (i+1)\sigma +(k+1)\rho}\not=0 \ \textrm{ for all } k\geq 0.
\end{equation}
So induction on $i$ gives
\begin{equation}
x^{ \rho}\tilde{\parti}_{1}.w_{\nu + i\sigma+k\rho}\not=0 \quad \textrm{ for } i,k\geq 0.
\end{equation}
It can be proved similarly for the other cases: $i,k\leq 0$; $i\geq
0,k\leq 0$; $i\leq 0,k\geq 0$. We omit the details. So we have
\begin{equation}\label{equ:7.*}x^{ \rho}\tilde{\parti}_{1}.w_{\nu +
i \sigma+k \rho}\not=0\qquad \for\;\;
i,k\in\mathbb{Z}.\end{equation}

Now we are ready to define the $v$'s. Take $v_{\nu+i\sigma}$ and $v_{\nu+ k
\rho}$ for $i,k \in \mathbb{Z}$ as they were in (\ref{equ:1.4}). Define $v_{\nu + i\sigma+k \rho} $'s
with $i, k \in \mathbb{Z}\backslash\{0\}$ by
\begin{equation}\label{equ:4.15}
x^{\rho}\tilde{\parti}_{1}.v_{\nu+i'\sigma+ (k'-1) \rho}=b_{1}
v_{\nu+i'\sigma+k' \rho}\quad \textrm{ for }i', k' \in \mathbb{Z} \textrm{ with } i'\not=0.
\end{equation}
Then it follows from (\ref{equ:1.4}) and (\ref{equ:4.15}) that
\begin{equation}\label{eq:2}
x^{\rho}\tilde{\parti}_{1}.v_{\nu+i\sigma+ (k-1) \rho}=b_{1}
v_{\nu+i\sigma+k \rho}\quad \textrm{ for }i, k \in \mathbb{Z}.
\end{equation}
For $k > 0$, (\ref{eq:2}) implies
\begin{equation}\label{equ:8.23}
x^{\pm \sigma}\tilde{\parti}_{1}.v_{\nu+i\sigma+ k \rho}=x^{\pm
\sigma}\tilde{\parti}_{1}.(\frac{x^{\rho}\tilde{\parti}_{1}}{b_1})^{k}.v_{\nu+i\sigma}=(\frac{x^{\rho}\tilde{\parti}_{1}}{b_1})^{k}.x^{\pm
\sigma}\tilde{\parti}_{1}.v_{\nu+i\sigma}=a_{1} v_{\nu+(i\pm
1)\sigma+ k \rho}.
\end{equation}
Moreover,
\begin{equation}\label{equ:8.24}
(\frac{x^{\rho}\tilde{\parti}_{1}}{b_1})^{k}.(x^{\pm
\sigma}\tilde{\parti}_{1}.v_{\nu+i\sigma- k \rho})= x^{\pm
\sigma}\tilde{\parti}_{1}.(\frac{x^{\rho}\tilde{\parti}_{1}}{b_1})^{k}.v_{\nu+i\sigma-
k \rho} =x^{\pm
\sigma}\tilde{\parti}_{1}.v_{\nu+i\sigma}=a_{1}v_{\nu+(i\pm
1)\sigma}
\end{equation}
for $k>0$, which together with (\ref {eq:2}) indicates $x^{\pm
\sigma}\tilde{\parti}_{1}.v_{\nu+i\sigma-k \rho}=a_{1} v_{\nu+(i\pm
1)\sigma-k \rho}$ . So we have
\begin{equation}\label{equ:1.6}
x^{\pm \sigma}\tilde{\parti}_{1}.v_{\nu+i\sigma+ k \rho}=a_{1}
v_{\nu+(i\pm 1)\sigma+ k \rho} \quad \textrm{for } i,k \in \mathbb{Z}.
\end{equation}
Similarly,  it can be deduced from (\ref{equ:1.4}) and
(\ref{equ:1.6}) that
\begin{equation}\label{equ:1.19}
x^{-\rho}\tilde{\parti}_{1}.v_{\nu +i\sigma+ k \rho}=b_{1}
v_{\nu+i\sigma+(k-1) \rho} \quad \textrm{for } i,k \in \mathbb{Z}.
\end{equation}

Finally we want to show that $x^{\sigma+\rho}\tilde{\parti}_{1}$,
$x^{\pm\sigma}\tilde{\parti}_{2}$ and
$x^{\pm\rho}\tilde{\parti}_{3}$ act on the $v$'s in the desired
expressions.

Note $[x^{\sigma+\rho}\tilde{\parti}_{1},x^{\pm
\sigma}\tilde{\parti}_{1}]=0$ and
$[x^{\sigma+\rho}\tilde{\parti}_{1},x^{\pm
\rho}\tilde{\parti}_{1}]=0$. From the assumption
$x^{\sigma+\rho}\tilde{\parti}_{1}.v_{\nu}=dv_{\nu+\sigma+ \rho}$,
it can be derived that
\begin{equation}\label{equ:1.11}
x^{\sigma+\rho}\tilde{\parti}_{1}.v_{\nu+i\sigma+ k \rho}=d
v_{\nu+(i+1)\sigma+(k+1) \rho}\quad \textrm{for } i,k \in \mathbb{Z}.
\end{equation}
Moreover, since
\begin{equation}
[x^{-\sigma}\tilde{\parti}_{2},x^{\sigma+\rho}\tilde{\parti}_{1}].v_{\nu}=\rho(\tilde{\parti}_{2})x^{\rho}\tilde{\parti}_{1}.v_{\nu}=b_1\rho(\tilde{\parti}_{2})v_{\nu+\rho}\not=0,
\end{equation}
we have $d\not=0$.

Write
\begin{equation}x^{\sigma}\tilde{\parti}_{2}.v_{\nu+i\sigma+ k \rho}=f^{+}(i,k)
v_{\nu+(i+1)\sigma+k
\rho}\qquad\for\;\;i,k\in\mathbb{Z}.\end{equation} Using
(\ref{equ:1.6}) in
$[x^{\sigma}\tilde{\parti}_{2},x^{\sigma}\tilde{\parti}_{1}].v_{\nu+i\sigma+
k \rho}=0$, we get $f^{+}(i+1,k)=f^{+}(i,k)$. Moreover, applying
(\ref{equ:4.15}) and (\ref{equ:1.11}) to
\begin{equation}
[x^{\sigma}\tilde{\parti}_{2},x^{\rho}\tilde{\parti}_{1}].v_{\nu+i\sigma+ k \rho}=\rho(\tilde{\parti}_{2})x^{\sigma+\rho}\tilde{\parti}_{1}.v_{\nu+i\sigma+ k \rho},
\end{equation}
we derive
$f^{+}(i,k+1)-f^{+}(i,k)=\rho(\tilde{\parti}_{2})\frac{d}{b_{1}}$.
Therefore, by induction, we get
\begin{equation}\label{equ:1.13}
f^{+}(i,k)=f^{+}(0,0)+k\rho(
\tilde{\parti}_{2})\frac{d}{b_{1}}\qquad \textrm{for }
i,k\in\mathbb{Z}.
\end{equation}

Write \begin{equation}x^{-\sigma}\tilde{\parti}_{2}.v_{\nu+i\sigma+
k \rho}=f^{-}(i,k) v_{\nu+(i-1)\sigma+k \rho}\qquad \for
\;\;i,k\in\mathbb{Z}.\end{equation} Then using (\ref{equ:1.6}) in
$[x^{-\sigma}\tilde{\parti}_{2},x^{\sigma}\tilde{\parti}_{1}].v_{\nu+i\sigma+
k \rho}=0$, we get $f^{-}(i+1,k)=f^{-}(i,k)$. Moreover, substituting
(\ref{equ:4.15}) and (\ref{equ:1.11}) into
\begin{equation}
[x^{-\sigma}\tilde{\parti}_{2},x^{\sigma+\rho}\tilde{\parti}_{1}].v_{\nu+i\sigma+ k \rho}=\rho(\tilde{\parti}_{2})x^{\rho}\tilde{\parti}_{1}.v_{\nu+i\sigma+ k \rho},
\end{equation}
we see that
$f^{-}(i+1,k+1)-f^{-}(i,k)=\rho(\tilde{\parti}_{2})\frac{b_{1}}{d}$.
Therefore, by induction, we obtain
\begin{equation}\label{equ:1.14}
f^{-}(i,k)=f^{-}(0,0)+k\rho(
\tilde{\parti}_{2})\frac{b_{1}}{d}\qquad \textrm{for }
i,k\in\mathbb{Z}.
\end{equation}

Set \begin{equation} x^{\rho}\tilde{\parti}_{3}.v_{\nu+i\sigma+ k
\rho}=g^{+}(i,k) v_{\nu+i\sigma+(k+1) \rho}\textrm{ and }
x^{-\rho}\tilde{\parti}_{3}.v_{\nu+i\sigma+ k \rho}=g^{-}(i,k)
v_{\nu+i\sigma+(k-1) \rho}\end{equation}for $i,k\in\mathbb{Z}.$
Observe that
$[x^{\rho}\tilde{\parti}_{3},x^{\rho}\tilde{\parti}_{1}]=0$ and
$[x^{\rho}\tilde{\parti}_{3},x^{\sigma}\tilde{\parti}_{1}]
=\sigma(\tilde{\parti}_{3})x^{\sigma+\rho}\tilde{\parti}_{1}$. As
(\ref{equ:1.13}), we deduce
\begin{equation}\label{equ:1.17}
g^{+}(i,k)=g^{+}(0,0)+ i
\sigma(\tilde{\parti}_{3})\frac{d}{a_{1}}\quad \textrm{for }
i,k\in\mathbb{Z}
\end{equation} by induction.
Likewise,
$[x^{-\rho}\tilde{\parti}_{3},x^{\rho}\tilde{\parti}_{1}]=0$,
$[x^{-\rho}\tilde{\parti}_{3},x^{\sigma+\rho}\tilde{\parti}_{1}]=\sigma(\tilde{\parti}_{3})
x^{\sigma}\tilde{\parti}_{1}$ and induction give rise to
\begin{equation}\label{equ:1.18}
g^{-}(i,k)=g^{-}(0,0)+ i
\sigma(\tilde{\parti}_{3})\frac{a_{1}}{d}\qquad \textrm{for }
i,k\in\mathbb{Z}.
\end{equation}

We get $b_{1}^{2}=d^{2}$ by applying (\ref{equ:4.15}),
(\ref{equ:1.6}), (\ref{equ:1.14}) and (\ref{equ:1.18}) to
\begin{equation}
[[x^{-\sigma}\tilde{\parti}_{2},x^{\rho}\tilde{\parti}_{1}],x^{-\rho}\tilde{\parti}_{3}].v_{\nu}=\rho(\tilde{\parti}_{2})\sigma(\tilde{\parti}_{3})x^{-\sigma}\tilde{\parti}_{1}.v_{\nu}.
\end{equation}
Then substituting (\ref{equ:1.13}),
(\ref{equ:1.14}), (\ref{equ:1.17}) and $b_{1}^{2}=d^{2}$ into
\begin{equation}
[x^{-\sigma}\tilde{\parti}_{2},[x^{\sigma}\tilde{\parti}_{2},x^{\rho}\tilde{\parti}_{3}]].v_{\nu}
=\rho^{2}(\tilde{\parti}_{2})x^{\rho}\tilde{\parti}_{3}.v_{\nu},
\end{equation}
we get $f^{-}(0,0)=f^{+}(0,0)$. Moreover, (\ref{equ:1.13}), (\ref{equ:1.14}) and $b_{1}^{2}=d^{2}$ show
\begin{equation}\label{equ:1.22}
f^{-}(i,k)=f^{+}(i,k)=a_{2}+k
\rho(\tilde{\parti}_{2})\frac{d}{b_{1}} \quad \textrm{ for } i, k
\in \mathbb{Z}.
\end{equation}

We obtain $a_{1}^{2}=d^{2}$ by inserting (\ref{equ:1.6}),
(\ref{equ:1.19}), (\ref{equ:1.13}) and (\ref{equ:1.17}) into
\begin{equation}
[x^{\rho}\tilde{\parti}_{3},[x^{-\rho}\tilde{\parti}_{1},x^{\sigma}\tilde{\parti}_{2}]].v_{\nu}=\rho(\tilde{\parti}_{2})\sigma(\tilde{\parti}_{3})x^{\sigma}\tilde{\parti}_{1}.v_{\nu}.
\end{equation}
Applying (\ref{equ:1.17}), (\ref{equ:1.18}), (\ref{equ:1.22}) and $a_{1}^{2}=d^{2}$ to
\begin{equation}
[x^{-\rho}\tilde{\parti}_{3},[x^{\rho}\tilde{\parti}_{3},x^{\sigma}\tilde{\parti}_{2}]].v_{\nu}=\sigma^{2}(\tilde{\parti}_{3})x^{\sigma}\tilde{\parti}_{2}.v_{\nu}, \end{equation}
we deduce $g^{-}(0,0)=g^{+}(0,0)$. Moreover, (\ref{equ:1.17}), (\ref{equ:1.18}) and $a_{1}^{2}=d^{2}$ indicate
\begin{equation}\label{equ:1.23}
g^{-}(i,k)=g^{+}(i,k)=a_{3}+ i
\sigma(\tilde{\parti}_{3})\frac{d}{a_{1}} \quad \textrm{ for } i, k
\in \mathbb{Z}.
\end{equation}
Thus the lemma follows. $\qquad\Box$ \vspace{0.2cm}

Next we begin to analyze the possible action of ${\cal S}$ on $M$.

Fix a $\mathbb{Z}$-basis $\{\ves_1, \ves_2, \ves_3 \}$ of $\G$.  Then $
\kn{\ves_1}\cap \kn{\ves_2} \cap \kn{\ves_3}=\{ 0 \}$. Picking any nonzero
vectors $
\parti_{1} \in \kn{\ves_2} \cap \kn{\ves_3}$, $ \parti_{2} \in \kn{\ves_1} \cap
\kn{\ves_3}$ and $ \parti_{3} \in \kn{\ves_1} \cap \kn{\ves_2}$, we get a
$\mathbb{F}$-basis of $D$. Obviously $\ves_1(\parti_{1})\not=0$,
$\ves_2(\parti_{2})\not=0$ and $\ves_3(\parti_{3})\not=0$. It can be verified that
 \begin{equation}\label{eq:31}
{\cal S}(\G, D)\;\mbox{is generated by}\; {\cal X}=\{
x^{\pm\ves_1}\parti_{2},  x^{\pm\ves_1}\parti_{3},
x^{\pm\ves_2}\parti_{1}, x^{\pm\ves_2}\parti_{3},
x^{\pm\ves_3}\parti_{1}, x^{\pm\ves_3}\parti_{2}\}.\end{equation}
Thus, we only need to derive how ${\cal X}$ act on
$M$. We analyze the action of ${\cal X}$ case by case in the following lemmas.

\begin{lemma}\label{le:5.4}
 If there exists some $\nu \in \G$ such that
\begin{eqnarray}
& & x^{-\ves_2}\parti_{1}.x^{\ves_2}\parti_{1}.w_{\nu} \not=0 \textrm{ and }
x^{-\ves_3}\parti_{1}.x^{\ves_3}\parti_{1}.w_{\nu} \not=0,\nonumber\\
& \textrm{or,} & x^{-\ves_1}\parti_{2}.x^{\ves_1}\parti_{2}.w_{\nu} \not=0 \textrm{ and }
x^{-\ves_3}\parti_{2}.x^{\ves_3}\parti_{2}.w_{\nu} \not=0,\nonumber\\
& \textrm{or,} & x^{-\ves_1}\parti_{3}.x^{\ves_1}\parti_{3}.w_{\nu} \not=0 \textrm { and }
x^{-\ves_2}\parti_{3}.x^{\ves_2}\parti_{3}.w_{\nu} \not=0,\nonumber
\end{eqnarray}
then $M$ is isomorphic to one of the following modules
for appropriate $\mu\in D^{\ast}$ and $\eta \in D^{\ast}\backslash\{0\}$:
$(\rmnum{1})    \mathscr{M}_{\mu};\  (\rmnum{2})   \mathscr{A}_{0,\eta};\ (\rmnum{3}) \mathscr{B}_{0,\eta}.$
\end{lemma}

\noindent{\bf{Proof.}} Assume
$x^{-\ves_2}\parti_{1}.x^{\ves_2}\parti_{1}.w_{\nu} \not=0$ and
$x^{-\ves_3}\parti_{1}.x^{\ves_3}\parti_{1}.w_{\nu} \not=0$ for some $\nu
\in \G$; the other cases can be proved similarly. By a translation
of the indices, we may assume $\nu=0$. Our two main goals in the proof are to properly choose
\begin{equation}
\{0\not=v_{i\ves_1+j\ves_2+k\ves_3}\in M_{i\ves_1+j\ves_2+k\ves_3}\mid i,j,k\in\mbb{Z}\}
 \end{equation}
 and to derive how ${\cal X}$ act on them (c.f. (\ref{eq:31})). The main idea is that we first specify the action of ${\cal X}$ on the three subspaces $\bigoplus_{j,k\in\mbb{Z}}M_{j\ves_2+k\ves_3}$, $\bigoplus_{i,k\in\mbb{Z}}M_{i\ves_1+k\ves_3}$ and
$\bigoplus_{i,j\in\mbb{Z}}M_{i\ves_1+j\ves_2}$, and then extend it
to the action of ${\cal X}$ on $M$. We process our proof in several
steps.\vspace{0.2cm}

{\it Step 1}. Elements $\{v_{j\ves_2+k\ves_3}\mid j,k \in
\mathbb{Z}\}$ and action of
$x^{\pm\ves_2}\parti_{1},\,x^{\pm\ves_2}\parti_{3},\,x^{\pm\ves_3}\parti_{1},\,x^{\pm\ves_3}\parti_{2}$.\vspace{0.2cm}

Since $x^{-\ves_2}\parti_{1}.x^{\ves_2}\parti_{1}.w_{0} \not=0$ and
$x^{-\ves_3}\parti_{1}.x^{\ves_3}\parti_{1}.w_{0} \not=0$, Lemma \ref{le:5.3} enables us to
choose $\{0\not=v_{j\ves_2+k\ves_3} \in M_{ j\ves_2+k\ves_3}\mid j,k \in
\mathbb{Z}\}$ such that
\begin{equation}
x^{\pm \ves_2}\parti_{1}.v_{j\ves_2+k\ves_3}=a_{1} v_{(j \pm 1)\ves_2+k\ves_3} \textrm{ and } x^{ \pm \ves_3}\parti_{1}.v_{j\ves_2+k\ves_3}=b_{1} v_{j \ves_2+(k \pm 1)\ves_3},
\end{equation}
where $a_{1}$ and $b_{1}$ are nonzero scalars satisfying
$x^{-\ves_2}\parti_{1}.x^{\ves_2}\parti_{1}.v_{0} =a_{1}^{2} v_{0}$ and $
x^{-\ves_3}\parti_{1}.x^{\ves_3}\parti_{1}.v_{0}=b_{1}^{2} v_{0}$,
respectively. Write $ x^{\ves_2+\ves_3}\parti_{1}.v_{0}=d_{1}
v_{\ves_2+\ves_3}$. Lemma \ref{le:5.3} says $a_{1}^{2}=b_{1}^{2}=d_{1}^{2}$.
Changing the sign of $v_{j\ves_2+k\ves_3}$ for $ j,k \in \mathbb{Z}$ if
necessary, we can take  $a_{1}=b_{1}=d_{1}$; equivalently, we get
$\{0\not=v_{j\ves_2+k\ves_3} \in M_{ j\ves_2+k\ves_3}\mid j,k \in \mathbb{Z}\}$
such that
\begin{eqnarray}
&&x^{\pm \ves_2}\parti_{1}.v_{j\ves_2+k\ves_3}=a_{1} v_{(j \pm 1)\ves_2+k\ves_3},\label{equ:3.32}\\
&&x^{ \pm \ves_3}\parti_{1}.v_{j\ves_2+k\ves_3}=a_{1} v_{j \ves_2+(k \pm 1)\ves_3},\label{equ:3.31}\\
&&x^{\ves_2+\ves_3}\parti_{1}.v_{j\ves_2+k\ves_3}=a_{1}
v_{(j+1) \ves_2+(k+1)\ves_3},\label{equ:3.33}\\
&& x^{ \pm \ves_3}\parti_{2}.v_{j\ves_2+k\ves_3}=(a_{2}+j\ves_2(
\parti_{2})) v_{j \ves_2+(k \pm
1)\ves_3},\label{equ:3.34}\\
&&  x^{\pm \ves_2}\parti_{3}.v_{j\ves_2+k\ves_3}=(a_{3}+ k\ves_3(
\parti_{3})) v_{(j \pm 1)\ves_2+k\ves_3},\label{equ:3.35}
\end{eqnarray}
where $a_{2}$ and $a_{3}$ are constants
determined by $x^{\ves_3}\parti_{2}.v_{0}=a_{2}v_{\ves_3}$ and $x^{
\ves_2}\parti_{3}.v_{0}=a_{3} v_{\ves_2}$.\vspace{0.2cm}

{\it Step 2}.  Elements $\{v_{i\ves_1+k\ves_3},
v_{i\ves_1+j\ves_2}\mid i,j,k \in \mathbb{Z}\}$ and action of
the set ${\cal X}$ (c.f. (\ref{eq:31})).\vspace{0.2cm}

According to (\ref{equ:3.32})
and (\ref{equ:3.31}), $x^{\pm \ves_2}\parti_{1}.v_{j\ves_2+k\ves_3}\neq 0$
and $x^{ \pm \ves_3}\parti_{1}.v_{j\ves_2+k\ves_3}\neq 0$ for any
$j,k\in\mbb{Z}$. By another proper translation of the indices, we
may assume that
\begin{equation}\label{eq:1}
a_{2}+j\ves_2(\parti_{2})\not=0 \textrm{ for }|j|\leq 1 \textrm{ and } a_{3}+ k\ves_3(
\parti_{3})\not=0 \textrm{ for } |k|\leq 1.
\end{equation}
Then we consider applying Lemma \ref{le:5.3} to the action of ${\cal X}$ on the subspaces
$\bigoplus_{i,k\in\mbb{Z}}M_{i\ves_1+k\ves_3}$ and
$\bigoplus_{i,j\in\mbb{Z}}M_{i\ves_1+j\ves_2}$. First it can be deduced from (\ref{eq:1}) that
\begin{equation}\label{e:2}
x^{-\ves_3}\parti_{2}.x^{\ves_3}\parti_{2}.v_{0}=a_{2}^{2}v_{0}\not=0 \textrm{ and } x^{
-\ves_2}\parti_{3}.x^{\ves_2}\parti_{3}.v_{0}=a_{3}^{2}v_{0}\not=0.
\end{equation}

Next we need the following claim.

\vspace{0.2cm}

 \emph{ Claim 1}. {\it $x^{
-\ves_1}\parti_{2}.x^{\ves_1}\parti_{2}.v_{0}\not=0$ and $x^{
-\ves_1}\parti_{3}.x^{\ves_1}\parti_{3}.v_{0}\not=0$.}\vspace{0.2cm}

We only give the proof of $x^{
-\ves_1}\parti_{2}.x^{\ves_1}\parti_{2}.v_{0}\not=0$. Inequality
$x^{ -\ves_1}\parti_{3}.x^{\ves_1}\parti_{3}.v_{0}\not=0$ can be similarly
proved.

Firstly, we want to prove $x^{ \ves_1}\parti_{2}.v_{0}\not=0$.
Assume  $x^{ \ves_1}\parti_{2}.v_{0}=0$ and then we will see this
leads to a contradiction. Since $a_2\not=0$, we have
\begin{equation}
x^{\ves_1}\parti_{2}.v_{\ves_3}=\frac{1}{a_2}x^{\ves_1}\parti_{2}.x^{\ves_3}\parti_{2}.v_{0}=\frac{1}{a_2}x^{\ves_3}\parti_{2}.x^{\ves_1}\parti_{2}.v_{0}=0,
\end{equation}
which further implies
\begin{equation}\label{equ:1.40}
x^{\ves_1+\ves_3}\parti_{2}.v_{0}=\frac{1}{\ves_1(\parti_{1})}[x^{\ves_3}\parti_{1},x^{\ves_1}\parti_{2}].v_{0}=-\frac{a_1}{\ves_1(\parti_{1})}x^{\ves_1}\parti_{2}.v_{\ves_3}=0.
\end{equation}
If $x^{ -\ves_1}\parti_{3}.v_{0}=0$, it can be derived
\begin{equation}
0\not=a_2 v_{\ves_3}=x^{\ves_3}\parti_{2}.v_{0}=\frac{1}{\ves_3(\parti_{3})}[x^{-\ves_1}\parti_{3},x^{\ves_1+\ves_3}\parti_{2}].v_{0}=0,
\end{equation}
which is absurd. Consider $x^{ -\ves_1}\parti_{3}.v_{0}\not=0$.
Since
\begin{equation}
x^{\ves_1}\parti_{2}.(x^{-\ves_1}\parti_{3}.v_{0})=x^{-\ves_1}\parti_{3}.x^{\ves_1}\parti_{2}.v_{0}=0,
\end{equation}
we obtain $x^{\ves_1}\parti_{2}.w_{-\ves_1}=0$, which further implies
\begin{equation}\label{equ:1.41}
x^{-\ves_3}\parti_{2}.(x^{\ves_1}\parti_{2}.w_{-\ves_1+\ves_3})=x^{\ves_1}\parti_{2}.(x^{-\ves_3}\parti_{2}.w_{-\ves_1+\ves_3})=0.
\end{equation}
 As $x^{-\ves_3}\parti_{2}.v_{\ves_3}=a_2v_{0}\not=0$,  (\ref{equ:1.41})
 leads to
$x^{\ves_1}\parti_{2}.w_{-\ves_1+\ves_3}=0$. Then
\begin{equation}\label{equ:1.44}
x^{\ves_1+\ves_3}\parti_{2}.w_{-\ves_1}=\frac{1}{\ves_1(\parti_{1})}[x^{\ves_3}\parti_{1},x^{\ves_1}\parti_{2}].w_{-\ves_1}=-\frac{1}{\ves_1(\parti_{1})}x^{\ves_1}\parti_{2}.(x^{\ves_3}\parti_{1}.w_{-\ves_1})=0.
\end{equation}
From (\ref{equ:1.40}) and (\ref{equ:1.44}), it can be deduced
\begin{equation}
0\not=a_2 v_{\ves_3}=x^{\ves_3}\parti_{2}.v_{0}=\frac{1}{\ves_3(\parti_{3})}[x^{-\ves_1}\parti_{3},x^{\ves_1+\ves_3}\parti_{2}].v_{0}=0,
\end{equation}
which leads to a contradiction. Hence we have
\begin{equation}\label{equ:1.46}
x^{\ves_1}\parti_{2}.v_{0}\not=0.
\end{equation}

Secondly, we want to prove $x^{ -\ves_1}\parti_{2}.w_{\ves_1}\not=0$. Suppose $x^{ -\ves_1}\parti_{2}.w_{\ves_1}=0$. Then we will get a contradiction. Note
\begin{equation}\label{equ:1.47}
x^{-\ves_3}\parti_{2}.(x^{ -\ves_1}\parti_{2}.w_{\ves_1+\ves_3})=x^{ -\ves_1}\parti_{2}.(x^{-\ves_3}\parti_{2}.w_{\ves_1+\ves_3})=0
\end{equation}
and $x^{-\ves_3}\parti_{2}.v_{\ves_3}=a_{2}v_{0}\not=0$. We therefore get $x^{
-\ves_1}\parti_{2}.w_{\ves_1+\ves_3}=0$. If $x^{ -\ves_1}\parti_{2}.v_{0}=0$, we derive
\begin{equation}
0\not=a_2 v_{\ves_3}=x^{\ves_3}\parti_{2}.v_{0}=\frac{1}{\ves_3(\parti_{3})\ves_1(\parti_{1})}[x^{\ves_1+\ves_3}(\ves_1(\parti_{1})\parti_{3}-\ves_3(\parti_{3})\parti_{1}),x^{-\ves_1}\parti_{2}].v_{0}=0,
\end{equation}
which is a contradiction. On the other hand, we assume $
x^{-\ves_1}\parti_{2}.v_{0}\not=0$. Since $x^{
-\ves_1}\parti_{2}.w_{\ves_1}=0$ leads to
\begin{equation}
x^{\ves_1}\parti_{3}.(x^{-\ves_1}\parti_{2}.v_{0})=x^{-\ves_1}\parti_{2}.(x^{\ves_1}\parti_{3}.v_{0})=0,
\end{equation}
we get $x^{\ves_1}\parti_{3}.w_{-\ves_1}=0$. Since
\begin{equation}
x^{-\ves_3}\parti_{2}.x^{-\ves_1}\parti_{2}.v_{\ves_3}=x^{-\ves_1}\parti_{2}.x^{-\ves_3}\parti_{2}.v_{\ves_3}=a_{2}x^{-\ves_1}\parti_{2}.v_{0}\not=0,
\end{equation}
we obtain $x^{-\ves_1}\parti_{2}.v_{\ves_3}\not=0$. By $x^{
-\ves_1}\parti_{2}.w_{\ves_1+\ves_3}=0$ mentioned above, we get
\begin{equation}
x^{\ves_1}\parti_{3}.(x^{-\ves_1}\parti_{2}.v_{\ves_3})=x^{-\ves_1}\parti_{2}.(x^{\ves_1}\parti_{3}.v_{\ves_3})=0,
\end{equation}
which implies $x^{\ves_1}\parti_{3}.w_{-\ves_1+\ves_3}=0$. Then we have
\begin{eqnarray}
0\not=\ves_3(\parti_{3})\ves_1(\parti_{1})x^{\ves_3}\parti_{2}.v_{0}&=&[[x^{\ves_3}\parti_{1},x^{\ves_1}\parti_{3}],x^{-\ves_1}\parti_{2}].v_{0}\nonumber\\
&=&x^{\ves_3}\parti_{1}.x^{\ves_1}\parti_{3}.(x^{-\ves_1}\parti_{2}.v_{0})-x^{\ves_1}\parti_{3}.(x^{\ves_3}\parti_{1}.x^{-\ves_1}\parti_{2}.v_{0})\nonumber\\
&&-x^{-\ves_1}\parti_{2}.([x^{\ves_3}\parti_{1},x^{\ves_1}\parti_{3}].v_{0})=0,
\end{eqnarray}
which is a contradiction. So we must have
\begin{equation}\label{equ:1.53}
x^{ -\ves_1}\parti_{2}.w_{\ves_1}\not=0.
\end{equation}

In summary, (\ref{equ:1.46}) and (\ref{equ:1.53}) give $x^{
-\ves_1}\parti_{2}.x^{\ves_1}\parti_{2}.v_{0}\not=0$. Symmetrically, considering the action of $x^{\pm\ves_2}\parti_{1}$, $x^{\pm\ves_2}\parti_{3}$, $x^{\pm\ves_1}\parti_{2}$ and $x^{\pm\ves_1}\parti_{3}$ on the subspace $\bigoplus_{i,j\in\mbb{Z}}M_{i\ves_1+j\ves_2}$, we can
similarly prove $x^{
-\ves_1}\parti_{3}.x^{\ves_1}\parti_{3}.v_{0}\not=0$. Hence Claim 1 holds.  \vspace{0.3cm}

Now we can apply Lemma \ref{le:5.3} to the action of ${\cal X}$ on the subspaces
$\bigoplus_{i,k\in\mbb{Z}}M_{i\ves_1+k\ves_3}$ and
$\bigoplus_{i,j\in\mbb{Z}}M_{i\ves_1+j\ves_2}$.

Observe that Claim 1 and (\ref{e:2}) give $x^{ -\ves_1}\parti_{2}.x^{\ves_1}\parti_{2}.v_{0}\not=0$ and $x^{
-\ves_3}\parti_{2}.x^{\ves_3}\parti_{2}.v_{0}\not=0$. Take $\{v_{k\ves_3}\mid k\in
\mathbb{Z}\}$ as they were in Step 1. Then Lemma \ref{le:5.3} enables us
to choose $\{0\not=v_{i\ves_1+k\ves_3} \in M_{i\ves_1+k\ves_3 } \mid i, k \in
\mathbb{Z},\,i\not=0\}$ such that
\begin{equation}\label{equ:3.51}
 x^{ \pm \ves_1}\parti_{2}.v_{i\ves_1+k\ves_3}=a_2
v_{(i\pm 1)\ves_1+k \ves_3}, \quad x^{ \pm \ves_3}\parti_{2}.v_{i\ves_1+k\ves_3}=a_2 v_{i\ves_1+(k \pm
1)\ves_3}
\end{equation}
for $i, k \in \mathbb{Z}$. Write $x^{
\ves_1+\ves_3}\parti_{2}.v_{0}=d_{2}v_{\ves_1+\ves_3}$ and
$x^{\ves_1}\parti_{3}.v_{0}=bv_{\ves_1}$. Lemma \ref{le:5.3} gives that
\begin{equation}\label{equ:3.52}
x^{\ves_1+\ves_3}\parti_{2}.v_{i\ves_1+k\ves_3}=d_{2}v_{(i+1)\ves_1+(k+1)\ves_3},
\end{equation}
\begin{equation}\label{equ:3.53}
 x^{ \pm \ves_3}\parti_{1}.v_{i\ves_1+k\ves_3}=(a_{1}+i\ves_1(\parti_{1})\frac{d_{2}}{a_2}) v_{ i \ves_1+(k \pm
1)\ves_3}
\end{equation}
and
\begin{equation}\label{equ:3.54}
  x^{\pm \ves_1}\parti_{3}.v_{i\ves_1+k\ves_3}=(b+ k\ves_3(
\parti_{3})\frac{d_{2}}{a_2}) v_{(i \pm 1) \ves_1+k\ves_3}
\end{equation}
for $i, k \in \mathbb{Z}$. Moreover, we have $a_{2}^{2}=d^{2}_{2}$.

Note that Claim 1 and (\ref{e:2}) show $x^{ -\ves_1}\parti_{3}.x^{\ves_1}\parti_{3}.v_{0}\not=0$
 and $x^{ -\ves_2}\parti_{3}.x^{\ves_2}\parti_{3}.v_{0}\not=0$. Take $\{v_{j\ves_2} \mid j \in
\mathbb{Z}\}$ and $\{v_{i\ves_1} \mid i \in
\mathbb{Z}\backslash\{0\}\}$ as they were in Step 1 and
(\ref{equ:3.51}), respectively. Then Lemma \ref{le:5.3} says that
there exist $\{0\not=v_{i\ves_1+j\ves_2} \in M_{i\ves_1+j\ves_2 }
\mid i,j \in \mathbb{Z}\backslash\{0\}\}$ such that
\begin{equation}\label{equ:3.55}
 x^{ \pm \ves_1}\parti_{3}.v_{i\ves_1+j\ves_2}=b v_{(i\pm
1)\ves_1+j\ves_2},\quad  x^{ \pm
\ves_2}\parti_{3}.v_{i\ves_1+j\ves_2}=a_{3} v_{i\ves_1+(j\pm 1)\ves_2}
\end{equation}
for $i,j \in \mathbb{Z}$. Writing $x^{
\ves_1+\ves_2}\parti_{3}.v_{0}=d_{3}v_{\ves_1+\ves_2}$, we have
\begin{equation}\label{equ:3.56}
x^{\ves_1+\ves_2}\parti_{3}.v_{i\ves_1+j\ves_2}=d_{3}v_{(i+1)\ves_1+(j+1)\ves_2},
\end{equation}
\begin{equation}\label{equ:3.57}
 x^{ \pm \ves_2}\parti_{1}.v_{i\ves_1+j\ves_2}=(a_{1}+i\ves_1(\parti_{1})\frac{d_{3}}{b}) v_{i\ves_1+(j\pm
1)\ves_2},
\end{equation}
\begin{equation}\label{equ:3.58}
  x^{\pm \ves_1}\parti_{2}.v_{i\ves_1+j\ves_2}=(a_2+ j\ves_2(\parti_{2})\frac{d_{3}}{a_{3}}) v_{(i \pm 1) \ves_1+j\ves_2}
\end{equation}
for $i,  j \in \mathbb{Z}$ again by Lemma \ref{le:5.3}. Moreover,
$a_{3}^{2}=b^{2}=d_{3}^{2}$.\vspace{0.2cm}

Among the constants $a_{2}$, $d_2$, $a_{3}$, $b$ and $d_{3}$, we have the following relations.

\vspace{0.2cm}

{\it Claim 2}. {\it $a_{2}=d_2$ and $a_{3}=b=d_{3}$.}\vspace{0.2cm}

Let $i$ be a nonzero integer such that
$a_{1}+i\ves_1(\parti_{1})\frac{d_{2}}{a_2}\not=0$ and
$a_{1}+i\ves_1(\parti_{1})\frac{d_{3}}{b}\not=0$. From (\ref{equ:3.53})
and (\ref{equ:3.57}) we deduce
\begin{equation}
x^{-\ves_3}\parti_{1}.x^{\ves_3}\parti_{1}.v_{i\ves_1}=(a_{1}+i\ves_1(\parti_{1})\frac{d_{2}}{a_2})^2 v_{i\ves_1}\not=0
\end{equation}
and
\begin{equation}
x^{-\ves_2}\parti_{1}.x^{\ves_2}\parti_{1}.v_{i\ves_1}=(a_{1}+i\ves_1(\parti_{1})\frac{d_{3}}{b})^2 v_{i\ves_1}\not=0.
\end{equation}
Then Lemma
\ref{le:5.3} implies
$(a_{1}+i\ves_1(\parti_{1})\frac{d_{2}}{a_2})^{2}=(a_{1}+i\ves_1(\parti_{1})\frac{d_{3}}{b})^{2}$,
which gives $\frac{d_{2}}{a_2}=\frac{d_{3}}{b}$. Similarly, by
(\ref{equ:3.34}), (\ref{equ:3.58}) and Lemma \ref{le:5.3}, we find $(a_2+
j\ves_2(\parti_{2}))^2=(a_2+ j\ves_2(\parti_{2})\frac{d_{3}}{a_{3}})^2$
for some $j\not=0$, which implies $a_{3}=d_{3}$.

Moreover, $a_{2}^{2}=d_{2}^{2}$ indicates $\frac{d_{2}}{a_2}=\pm 1$.
To prove the claim, it suffices to show $\frac{d_{2}}{a_2}=1$.
Suppose $\frac{d_{2}}{a_2}=-1$, and then we will see this leads to a
contradiction.

Since $a_{1}+\ves_1(\parti_{1})\not=0$ or
$a_{1}-\ves_1(\parti_{1})\not=0$, we may assume
$a_{1}-\ves_1(\parti_{1})\not=0$; the other case can be proved
similarly. By (\ref{equ:3.53}), (\ref{equ:3.57}) and the assumption $\frac{d_{2}}{a_2}=\frac{d_{3}}{b}=-1$, we have
\begin{equation}
x^{-\ves_3}\parti_{1}.x^{\ves_3}\parti_{1}.v_{\ves_1}=(a_{1}-\ves_1(\parti_{1}))^{2}v_{\ves_1}\not=0,
\end{equation}
\begin{equation}
x^{-\ves_2}\parti_{1}.x^{\ves_2}\parti_{1}.v_{\ves_1}=(a_{1}-\ves_1(\parti_{1}))^{2}v_{\ves_1}\not=0.
\end{equation}
Take $\{v_{\ves_1+k\ves_3} \mid k \in \mathbb{Z}\}$ and
$\{v_{\ves_1+j\ves_2} \mid j \in \mathbb{Z}\}$ as they were in
(\ref{equ:3.51}) and (\ref{equ:3.55}), respectively. According to
Lemma \ref{le:5.3}, we can choose $\{0\not=v_{\ves_1 + j \ves_2 +
k\ves_3} \in M_{ \ves_1 + j \ves_2+k \ves_3}\mid j,k \in
\mathbb{Z}\backslash\{0\}\}$ such that
\begin{equation}\label{equ:3.62}
 x^{\pm
\ves_2}\parti_{1}.v_{\ves_1 + j \ves_2+k \ves_3}=(a_{1}-\ves_1(\parti_{1})) v_{\ves_1+(j \pm 1) \ves_2+k\ves_3},\ \ \ \ \ \ \
\end{equation}
\begin{equation}\label{equ:3.61}
x^{ \pm \ves_3}\parti_{1}.v_{\ves_1 + j \ves_2+k \ves_3}=(a_{1}-\ves_1(\parti_{1}))
v_{\ves_1+ j \ves_2+(k \pm 1)\ves_3}\ \ \ \ \ \ \
\end{equation}
for $j,k\in\mbb{Z}$. Moreover, Lemma \ref{le:5.3}, (\ref{equ:3.51}) and (\ref{equ:3.55}) give
\begin{equation}\label{equ:3.63}
x^{\ves_2+\ves_3}\parti_{1}.v_{\ves_1 + j \ves_2+k \ves_3}=s(a_{1}-\ves_1(\parti_{1}))
v_{\ves_1+(j+1) \ves_2+(k+1)\ves_3},
\end{equation}
\begin{equation}\label{equ:3.64}
 x^{ \pm \ves_3}\parti_{2}.v_{\ves_1 + j \ves_2+k \ves_3}=(a_2+sj
\ves_2(\parti_{2})) v_{\ves_1+ j \ves_2+(k \pm 1)\ves_3},\ \ \ \
\end{equation}
\begin{equation}\label{equ:3.65}
  x^{\pm \ves_2}\parti_{3}.v_{\ves_1 + j \ves_2+k \ves_3}=(a_{3}+ sk\ves_3(
\parti_{3})) v_{\ves_1+(j \pm 1)\ves_2+k\ves_3}\ \ \ \ \
\end{equation}
for $j,k \in\mbb{Z}$, where $s\in \{ 1, -1 \}$ is to be specified.
Substituting (\ref{equ:3.34}), (\ref{equ:3.58}), (\ref{equ:3.64})
and $a_{3}=d_{3}$ into
$[x^{\ves_1}\parti_{2},x^{\ves_3}\parti_{2}].v_{\ves_2}=0 $, and making use of (\ref{eq:1}), we get
\begin{equation}\label{equ:3.66}
x^{\ves_1}\parti_{2}.v_{\ves_2+\ves_3}=(a_2+s\ves_2(
\parti_{2}))v_{\ves_1+\ves_2+\ves_3}.
\end{equation}
Using (\ref{equ:3.34}), (\ref{equ:3.58}), (\ref{equ:3.64}),
(\ref{equ:3.66}) and $a_{3}=d_{3}$ in $
[x^{\ves_1}\parti_{2},x^{-\ves_3}\parti_{2}].v_{\ves_2+\ves_3}=0 $, we find
\begin{equation}
(a_2+\ves_2(
\parti_{2}))^{2}=(a_2+s\ves_2(
\parti_{2}))^{2},
\end{equation}
which implies $s=1$. Then (\ref{equ:3.31}), (\ref{equ:3.58}), (\ref{equ:3.61}), (\ref{equ:3.66}) and $a_{3}=d_{3}$ give
\begin{equation}
x^{\ves_1+\ves_3}\parti_{2}.v_{\ves_2}=\frac{1}{\ves_1(\parti_{1})}[x^{\ves_3}\parti_{1},x^{\ves_1}\parti_{2}].v_{\ves_2}=-(a_2+
\ves_2(\parti_{2}))v_{\ves_1+\ves_2+\ves_3}.
\end{equation}
By $
[x^{\ves_1+\ves_3}\parti_{2},x^{\ves_2}\parti_{1}].v_{0}=[x^{\ves_1}\parti_{2},x^{\ves_2+\ves_3}\parti_{1}].v_{0}$, we get
\begin{equation}
a_{1}\ves_2(\parti_{2})+a_2\ves_1(\parti_{1})=0.
\end{equation}
Applying it to
$[x^{-\ves_1}\parti_{2},[x^{\ves_1}\parti_{2},x^{\ves_2}\parti_{1}]].v_{0}=\ves_2^{2}(\parti_{2})x^{\ves_2}\parti_{1}.v_{0}$,
we obtain
\begin{equation}
a_2\ves_1(\parti_{1})\ves_2(\parti_{2})=0,
\end{equation}
which is absurd. So we must have $
\frac{d_{2}}{a_2}=\frac{d_{3}}{b}=1$, from which it follows
\begin{equation}\label{equ:3.72}
a_{2}=d_2,\quad a_{3}=b=d_{3}.
\end{equation}
Therefore Claim 2 holds.  \vspace{0.2cm}

{\it Step 3. Extension to $\{v_{i\ves_1+j\ves_2+k\ves_3}\mid
i,j,k \in \mathbb{Z}\}$ and general action of
${\cal X}$.}\vspace{0.2cm}

In Step 1 and Step 2 we have properly chosen $\{v_{j\ves_2+k\ves_3}, v_{i\ves_1+k\ves_3},
v_{i\ves_1+j\ves_2}\mid i,j,k \in \mathbb{Z}\}$ and determined part of the action of
${\cal X}$ (see (\ref{equ:3.32})--(\ref{equ:3.35}), (\ref{equ:3.51})--(\ref{equ:3.58}) and Claim 2).
Now we start to extend it to the general action of ${\cal X}$.

As (\ref{equ:3.53}), (\ref{equ:3.57}) and Claim 2 show
\begin{equation}
x^{-\ves_3}\parti_{1}.x^{\ves_3}\parti_{1}.v_{i\ves_1}=(a_{1}+i\ves_1(\parti_{1}))^2 v_{i\ves_1}
\end{equation}
and
\begin{equation}
x^{-\ves_2}\parti_{1}.x^{\ves_2}\parti_{1}.v_{i\ves_1}=(a_{1}+i\ves_1(\parti_{1}))^2 v_{i\ves_1},
\end{equation}
we can apply Lemma \ref{le:5.3} to the action of $x^{\pm\ves_2}\parti_{1},x^{\pm\ves_2}\parti_{3},$
$x^{\pm\ves_3}\parti_{1}$ and $x^{\pm\ves_3}\parti_{2}$ on the subspaces
\begin{equation}
\bigoplus_{j,k\in\mbb{Z}}M_{i\ves_1+j\ves_2+k\ves_3} \qquad \textrm{ for }i\in\mbb{Z}\backslash\{0\} \textrm{ such that } a_{1}+i\ves_1(\parti_{1})\not=0.
\end{equation}
Seeing that there exists at most one $i\in\mathbb{Z}$ such that
$a_{1}+i\ves_1(\parti_{1})=0$, we denote it by $i'$. Obviously $i'\not=0$. Fix any $i\in\mbb{Z}\backslash\{0,i'\}$. Then Lemma \ref{le:5.3} enables us
to choose $\{0\not=v_{i\ves_1+ j \ves_2+k\ves_3 } \in M_{ i\ves_1+ j
\ves_2+k\ves_3 } \mid j,k \in \mathbb{Z}\backslash\{0\}\}$ such that
\begin{equation}\label{equ:4.79}
 x^{\pm
\ves_2}\parti_{1}.v_{i\ves_1+ j \ves_2+k\ves_3 }=(a_{1}+i\ves_1(\parti_{1}))
v_{i\ves_1+(j \pm 1) \ves_2+k\ves_3},
\end{equation}
\begin{equation}\label{equ:4.78}
x^{ \pm \ves_3}\parti_{1}.v_{i\ves_1+ j \ves_2+k\ves_3}=(a_{1}+i\ves_1(\parti_{1})) v_{ i\ves_1+ j \ves_2+(k \pm 1)\ves_3}
\end{equation}
for $j,k \in \mathbb{Z}$. Lemma \ref{le:5.3} further implies
\begin{equation}\label{equ:4.80}
x^{\ves_2+\ves_3}\parti_{1}.v_{i\ves_1+ j \ves_2+k\ves_3}=s_{i}(a_{1}+i\ves_1(\parti_{1})) v_{ i\ves_1+(j+1) \ves_2+(k+1)\ves_3},
\end{equation}
\begin{equation}\label{equ:4.81}
 x^{ \pm \ves_3}\parti_{2}.v_{i\ves_1+ j \ves_2+k\ves_3 }=(a_2+s_{i}j
\ves_2(\parti_{2})) v_{i\ves_1+ j \ves_2+(k \pm
1)\ves_3},\ \ \ \ \ \
\end{equation}
\begin{equation}\label{equ:4.82}
  x^{\pm \ves_2}\parti_{3}.v_{i\ves_1+ j \ves_2+k\ves_3 }=(a_{3}+ s_{i}k
\ves_3(\parti_{3})) v_{i\ves_1+(j \pm 1) \ves_2+k\ves_3}\ \ \ \ \ \
\end{equation}
for $j,k \in \mathbb{Z}$, where $s_{i}\in \{ 1,  -1 \}$ is to be determined.

Set $s_{0}=1$. Then (\ref{equ:3.32})--(\ref{equ:3.35}) coincide with (\ref{equ:4.79})--(\ref{equ:4.82}).

If $i'$ does not exist, (\ref{equ:4.79})--(\ref{equ:4.82}) together
with (\ref{equ:3.32})--(\ref{equ:3.35}) show the general action of
$x^{\pm\ves_2}\parti_{1}$, $x^{\pm\ves_2}\parti_{3}$,
$x^{\pm\ves_3}\parti_{1}$ and $x^{\pm\ves_3}\parti_{2}$. To complete
this step, we only need to determine $s_i$ for
$i\in\mbb{Z}\backslash\{0\}$ and to specify the general action of
$x^{\pm\ves_1}\parti_{2}$ and $x^{\pm\ves_1}\parti_{3}$.

If $i'$ does exist, except for determining $s_i$ for
$i\in\mbb{Z}\backslash\{0, i'\}$, we need to determine
$\{v_{i'\ves_1+ j \ves_2+k\ves_3} \mid
j,k\in\mbb{Z}\backslash\{0\}\}$ and to specify the action of
$x^{\pm\ves_2}\parti_{1}$, $x^{\pm\ves_2}\parti_{3}$,
$x^{\pm\ves_3}\parti_{1}$ and $x^{\pm\ves_3}\parti_{2}$ on the
subspace $\bigoplus_{j,k\in\mbb{Z}}M_{i'\ves_1+j\ves_2+k\ves_3}$.
Moreover, the general action of $x^{\pm\ves_1}\parti_{2}$ and
$x^{\pm\ves_1}\parti_{3}$ needs to be determined.

We give the details in the following.
\vspace{0.1cm}

\emph{Case 1. $i'$ does not exist.}\vspace{0.1cm}

To begin with, we want to determine $s_i$ for $i \in\mathbb{Z}\backslash\{0\}$. Let $p$ be any fixed
integer. Applying (\ref{equ:3.58}) and (\ref{equ:4.81}) to $
[x^{\ves_1}\parti_{2},x^{\ves_3}\parti_{2}].v_{ p\ves_1+\ves_2}=0 $, we have
\begin{equation}\label{equ:8.30}
x^{\ves_1}\parti_{2}.v_{ p\ves_1+\ves_2+\ves_3}=\frac{(a_2+\ves_2(\parti_{2}))(a_2+s_{p+1}
\ves_2(\parti_{2}))}{a_2+ s_{p}\ves_2(\parti_{2})}v_{ (p+1)\ves_1+\ves_2+\ves_3}.
\end{equation}
Moreover, substituting (\ref{equ:3.58}), (\ref{equ:4.81}) and (\ref{equ:8.30}) into $
[x^{\ves_1}\parti_{2},x^{-\ves_3}\parti_{2}].v_{p\ves_1+\ves_2+\ves_3}=0 $, and making use of (\ref{eq:1}),
we obtain
\begin{equation}\label{equ:4.83}
(a_2+ s_{p}\ves_2(\parti_{2}))^{2}=(a_2+s_{p+1}
\ves_2(\parti_{2}))^{2},
\end{equation}
which implies $s_{p}=s_{p+1}$. Since $p\in\mathbb{Z}$ is arbitrary
and $s_{0}=1$,  induction shows
\begin{equation}\label{eq:3}
s_{i}=1 \quad \textrm{for all } i\in\mathbb{Z}.
\end{equation}
So (\ref{equ:4.79})--(\ref{equ:4.82}) together with
(\ref{equ:3.32})--(\ref{equ:3.35}) show
\begin{equation}\label{e:3}
 x^{\pm
\ves_2}\parti_{1}.v_{i\ves_1+ j \ves_2+k\ves_3 }=(a_{1}+i\ves_1(\parti_{1}))
v_{i\ves_1+(j \pm 1) \ves_2+k\ves_3},
\end{equation}
\begin{equation}\label{e:4}
x^{ \pm \ves_3}\parti_{1}.v_{i\ves_1+ j \ves_2+k\ves_3}=(a_{1}+i\ves_1(\parti_{1})) v_{ i\ves_1+ j \ves_2+(k \pm 1)\ves_3},
\end{equation}
\begin{equation}\label{e:5}
 x^{ \pm \ves_3}\parti_{2}.v_{i\ves_1+ j \ves_2+k\ves_3 }=(a_2+j
\ves_2(\parti_{2})) v_{i\ves_1+ j \ves_2+(k \pm
1)\ves_3},
\end{equation}
\begin{equation}\label{e:6}
  x^{\pm \ves_2}\parti_{3}.v_{i\ves_1+ j \ves_2+k\ves_3 }=(a_{3}+ k
\ves_3(\parti_{3})) v_{i\ves_1+(j \pm 1) \ves_2+k\ves_3}
\end{equation}
for $i,j,k \in \mathbb{Z}$.

Next we want to derive the general action of $x^{\pm\ves_1}\parti_{2}$ and $x^{\pm \ves_1}\parti_{3}$.

We first consider determining the action of $x^{\pm\ves_1}\parti_{2}$.
Recall
that (\ref{equ:3.58}) and Claim 2 give
\begin{equation}\label{e:9}
x^{\pm\ves_1}\parti_{2}.v_{i\ves_1+
j\ves_2}=(a_2+j\ves_2(\parti_{2}))v_{(i\pm 1)\ves_1+ j\ves_2} \ \textrm{ for } i,j\in\mbb{Z}.
\end{equation}
Since
\begin{equation}\label{equ:8.32}
x^{\pm\ves_1}\parti_{2}.x^{\pm\ves_3}\parti_{2}.v_{ i\ves_1+ j \ves_2+k\ves_3 }=x^{\pm\ves_3}\parti_{2}.x^{\pm\ves_1}\parti_{2}.v_{ i\ves_1+ j \ves_2+k\ves_3 },
 \end{equation}
equation (\ref{e:5}) and induction on $k$ lead to
\begin{equation}\label{equ:3.80}
x^{\pm \ves_1}\parti_{2}.v_{ i\ves_1+ j \ves_2+k\ves_3 }=(a_2+j
\ves_2(\parti_{2}))v_{(i\pm1)\ves_1+ j \ves_2+k\ves_3}
\end{equation}
for $i ,k \in \mathbb{Z}$ and $j\in\mathbb{Z}$ such that $a_2+j\ves_2(\parti_{2})\not=0$. As there exists at most one $j\in\mathbb{Z}$ such that $a_2+j
\ves_2(\parti_{2})=0$, we denote it by $j'$. If $j'$ exists, then (\ref{e:5}) implies
\begin{equation}
x^{\pm \ves_1}\parti_{2}.v_{ i\ves_1+ j'\ves_2+k\ves_3}=\pm\frac{1}{\ves_1(\parti_{1})\ves_3(\parti_{3})}[x^{-\ves_3}\parti_{1},[x^{\pm
\ves_1}\parti_{3},x^{\ves_3}\parti_{2}]].v_{ i\ves_1+ j' \ves_2+k\ves_3}=0
\end{equation}
for $i,k \in \mathbb{Z}$, which coincides with (\ref{equ:3.80}). So
no matter such $j'$ exists or not, we have
\begin{equation}\label{eq:4}
x^{\pm \ves_1}\parti_{2}.v_{ i\ves_1+ j \ves_2+k\ves_3 }=(a_2+j
\ves_2(\parti_{2}))v_{(i\pm1)\ves_1+ j \ves_2+k\ves_3} \ \textrm{ for } i ,j,k \in \mathbb{Z}.
\end{equation}

It holds similarly that
\begin{equation}\label{equ:3.81}
x^{\pm \ves_1}\parti_{3}.v_{ i\ves_1+ j \ves_2+k\ves_3 }=(a_{3}+k
\ves_3(\parti_{3}))v_{(i\pm1)\ves_1+ j \ves_2+k\ves_3} \ \textrm{ for }i,j,k \in \mathbb{Z}.
\end{equation}

To sum up, (\ref{e:3})--(\ref{e:6}) together with (\ref{eq:4}) and
(\ref{equ:3.81}) show the action of ${\cal X}$ on $M$. Define
$\mu\in D ^{\ast} $ by $\mu(\parti_{1})=a_{1}$,
$\mu(\parti_{2})=a_2$ and $\mu(\parti_{3})=a_{3}$. Then we can write the action of ${\cal X}$ uniformly by
\begin{equation}\label{e:8}
x^{\es}\parti.v_{ i\ves_1+ j \ves_2+k\ves_3 }=(\mu+i\ves_1+ j \ves_2+k\ves_3)(\parti)v_{i\ves_1+ j \ves_2+k\ves_3+\es}
\end{equation}
for $i,j,k\in\mbb{Z}$, $\es \in \{ \pm\ves_1, \pm\ves_2,\pm\ves_3
\}$ and $\parti \in \kn{\es}\cap\{\parti_{1}, \parti_{2},
\parti_{3}\}$.
Since ${\cal X}$
generates ${\cal S}(\G, D)$ (c.f. (\ref{eq:31})), from (\ref{e:8}) we deduce
\begin{equation}
M \simeq \mathscr{M}_{\mu}.
\end{equation} \vspace{0.1cm}

\emph{Case 2. $i'$
does exist.}\vspace{0.1cm}

Recall that $j'$ denotes the possible integer satisfying $a_2+j'
\ves_2(\parti_{2})=0$. Likewise, we denote by $k'$ the possible integer satisfying $a_{3}+k'
\ves_3(\parti_{3})=0$.\vspace{0.2cm}

{\it Firstly, we want to determine $s_i$ for $i
\in\mathbb{Z}\backslash\{0,i'\}$ and obtain:}

1) the action of $x^{\pm\ves_2}\parti_{1},x^{\pm\ves_2}\parti_{3},$
$x^{\pm\ves_3}\parti_{1}$ and $x^{\pm\ves_3}\parti_{2}$ on $\bigoplus_{i,j,k\in\mbb{Z};i\not=i'}M_{i\ves_1+j\ves_2+k\ves_3}$,\vspace{0.1cm}

 2) the action of  $x^{\ves_1}\parti_{2}$ and $x^{\ves_1}\parti_{3}$ on $\bigoplus_{i,j,k\in\mbb{Z};i\not=i',i'-1}M_{i\ves_1+j\ves_2+k\ves_3}$,\vspace{0.1cm}

 3)  the action of $x^{-\ves_1}\parti_{2}$ and $x^{-\ves_1}\parti_{3}$ on $\bigoplus_{i,j,k\in\mbb{Z};i\not=i',i'+1}M_{i\ves_1+j\ves_2+k\ves_3}$.\vspace{0.2cm}

Since (\ref{equ:8.30}) and (\ref{equ:4.83}) also
hold for $p\in \mathbb{Z}\backslash\{i',\, i'-1\}$ in this case, we have $s_{p}=s_{p+1}$ for
$p\in \mathbb{Z}\backslash\{i',\, i'-1\}$. Moreover, we shall show $s_{i'-1}=s_{i'+1}$. As (\ref{equ:3.57}), (\ref{equ:3.58}) and
Claim 2 indicate
\begin{eqnarray}
  x^{2 \ves_1}\parti_{2}.v_{i\ves_1+j\ves_2}&=&\frac{1}{4\ves_1^{2}(\parti_{1})\ves_2(\parti_{2})}[[x^{\ves_1}\parti_{2},[x^{\ves_1}\parti_{2},x^{\ves_2}\parti_{1}]],x^{-\ves_2}\parti_{1}].v_{i\ves_1+j\ves_2}\nonumber\\
  &=&(a_2+ j\ves_2(\parti_{2})) v_{(i +2) \ves_1+j\ves_2}\quad \textrm{ for }i, j \in \mathbb{Z},\label{equ:6.1}
\end{eqnarray}
applying it and (\ref{equ:4.81}) to $
[x^{2\ves_1}\parti_{2},x^{\ves_3}\parti_{2}].v_{ (i'-1)\ves_1+\ves_2}=0$, we
get
\begin{equation}\label{equ:8.31}
x^{2\ves_1}\parti_{2}.v_{ (i'-1)\ves_1+\ves_2+\ves_3}=\frac{(a_2+\ves_2(\parti_{2}))(a_2+s_{i'+1}
\ves_2(\parti_{2}))}{a_2+ s_{i'-1}\ves_2(\parti_{2})}v_{ (i'+1)\ves_1+\ves_2+\ves_3}.
\end{equation}
 Then inserting (\ref{equ:4.81}), (\ref{equ:6.1}) and (\ref{equ:8.31})
into
$[x^{2\ves_1}\parti_{2},x^{-\ves_3}\parti_{2}].v_{(i'-1)\ves_1+\ves_2+\ves_3}=0 $, and making use of (\ref{eq:1}),
we obtain
\begin{equation}\label{equ:6.2}
(a_2+ s_{i'-1}\ves_2(\parti_{2}))^{2}=(a_2+s_{i'+1}
\ves_2(\parti_{2}))^{2},
\end{equation}
which implies $s_{i'-1}=s_{i'+1}$. Thus this and
$s_{p}=s_{p+1}$ for $p\in \mathbb{Z}\backslash\{i',\, i'-1\}$
indicate
\begin{equation}
s_{i}=s_{0}=1 \textrm{ for } i\in\mbb{Z}\backslash\{i'\}.
\end{equation}
So (\ref{equ:4.79})--(\ref{equ:4.82}) together with
(\ref{equ:3.32})--(\ref{equ:3.35}) show
\begin{equation}\label{e:13}
 x^{\pm
\ves_2}\parti_{1}.v_{i\ves_1+ j \ves_2+k\ves_3 }=(a_{1}+i\ves_1(\parti_{1}))
v_{i\ves_1+(j \pm 1) \ves_2+k\ves_3},
\end{equation}
\begin{equation}\label{e:11}
x^{ \pm \ves_3}\parti_{1}.v_{i\ves_1+ j \ves_2+k\ves_3}=(a_{1}+i\ves_1(\parti_{1})) v_{ i\ves_1+ j \ves_2+(k \pm 1)\ves_3},
\end{equation}
\begin{equation}\label{e:10}
 x^{ \pm \ves_3}\parti_{2}.v_{i\ves_1+ j \ves_2+k\ves_3 }=(a_2+j
\ves_2(\parti_{2})) v_{i\ves_1+ j \ves_2+(k \pm
1)\ves_3},
\end{equation}
\begin{equation}\label{e:12}
  x^{\pm \ves_2}\parti_{3}.v_{i\ves_1+ j \ves_2+k\ves_3 }=(a_{3}+ k
\ves_3(\parti_{3})) v_{i\ves_1+(j \pm 1) \ves_2+k\ves_3}
\end{equation}
for $i,j,k \in \mathbb{Z}$ with $i\not=i'$.

Similar arguments as those from (\ref{e:9}) to
(\ref{equ:3.81}) show
\begin{equation}\label{equ:8.33}
x^{ \ves_1}\parti_{2}.v_{ i\ves_1+ j \ves_2+k\ves_3 }=(a_2+j
\ves_2(\parti_{2}))v_{(i+1)\ves_1+ j \ves_2+k\ves_3}\ \textrm{ for }i,j,k \in \mathbb{Z},\,i\not=i',i'-1,
\end{equation}
\begin{equation}\label{equ:8.35}
x^{ \ves_1}\parti_{3}.v_{ i\ves_1+ j \ves_2+k\ves_3 }=(a_{3}+k
\ves_3(\parti_{3}))v_{(i+1)\ves_1+ j \ves_2+k\ves_3}\ \textrm{ for }i,j,k \in \mathbb{Z},\,i\not=i',i'-1,
\end{equation}
\begin{equation}\label{equ:8.34}
x^{ -\ves_1}\parti_{2}.v_{ i\ves_1+ j \ves_2+k\ves_3 }=(a_2+j
\ves_2(\parti_{2}))v_{(i-1)\ves_1+ j \ves_2+k\ves_3}\ \textrm{ for }i,j,k \in \mathbb{Z},\,i\not=i',i'+1,
\end{equation}
\begin{equation}\label{equ:8.36}
x^{-\ves_1}\parti_{3}.v_{ i\ves_1+ j \ves_2+k\ves_3 }=(a_{3}+k
\ves_3(\parti_{3}))v_{(i-1)\ves_1+ j \ves_2+k\ves_3}\ \textrm{ for }i,j,k \in \mathbb{Z},\,i\not=i',i'+1.
\end{equation}
\vspace{0.1cm}

{\it Secondly, we want to define $\{v_{i'\ves_1+ j \ves_2+k\ves_3}\mid j\in\mbb{Z}\backslash\{0,j'\}, k\in\mbb{Z}\backslash\{0\}\}$, and derive}

 1) $\textrm{the action of } x^{\pm\ves_3}\parti_{1} \textrm{ and } x^{\pm\ves_3}\parti_{2} \textrm{ on } \bigoplus_{j,k\in\mbb{Z}; j\not=j'}M_{i'\ves_1+j\ves_2+k\ves_3}$,\vspace{0.1cm}

 2) $\textrm{the action of } x^{\ves_1}\parti_{2} \textrm{ and } x^{\ves_1}\parti_{3} \textrm{ on } \bigoplus_{j,k\in\mbb{Z};j\not=j'}M_{(i'-1)\ves_1+j\ves_2+k\ves_3}$,\vspace{0.1cm}

 3) $\textrm{the action of } x^{-\ves_1}\parti_{2} \textrm{ and } x^{-\ves_1}\parti_{3} \textrm{ on } \bigoplus_{j,k\in\mbb{Z};j\not=j'}M_{i'\ves_1+j\ves_2+k\ves_3}$,\vspace{0.1cm}

 4) $\textrm{the action of } x^{\ves_2}\parti_{1} \textrm{ and } x^{\ves_2}\parti_{3} \textrm{ on } \bigoplus_{j,k\in\mbb{Z}; j\not=j',j'-1}M_{i'\ves_1+j\ves_2+k\ves_3}$,\vspace{0.1cm}

 5) $\textrm{the action of } x^{-\ves_2}\parti_{1} \textrm{ and } x^{-\ves_2}\parti_{3} \textrm{ on } \bigoplus_{j,k\in\mbb{Z}; j\not=j',j'+1}M_{i'\ves_1+j\ves_2+k\ves_3}$.\vspace{0.2cm}

To define the $v$'s, we first need to prove $x^{\ves_1}\parti_{2}.v_{(i'-1)\ves_1+ j
\ves_2+k\ves_3}\not=0$ for $j, k\in\mbb{Z}$ with $j\not=j'$. Fix any
$j\in\mathbb{Z}\backslash\{j'\}$. From (\ref{equ:3.58}) and Claim 2 we know that
\begin{equation}\label{eq:5}
x^{\ves_1}\parti_{2}.v_{(i'-1)\ves_1+ j\ves_2}=(a_2+j\ves_2(\parti_{2}))v_{i'\ves_1+ j\ves_2}\not=0.
 \end{equation}
 Since (\ref{e:10}) leads to
 \begin{eqnarray}
x^{\mp\ves_3}\parti_{2}.x^{\ves_1}\parti_{2}.v_{(i'-1)\ves_1+ j\ves_2+(k\pm1)\ves_3}&=&x^{\ves_1}\parti_{2}.(x^{\mp\ves_3}\parti_{2}.v_{(i'-1)\ves_1+ j\ves_2+(k\pm1)\ves_3})\nonumber\\
&=& (a_2+j\ves_2(\parti_{2}))x^{\ves_1}\parti_{2}.v_{(i'-1)\ves_1+ j\ves_2+k\ves_3},
\end{eqnarray}
we have $x^{\ves_1}\parti_{2}.v_{(i'-1)\ves_1+ j
\ves_2+k\ves_3}\not=0$ for $k\in\mbb{Z}$ by (\ref{eq:5}) and induction on $k$. So we get
\begin{equation}
x^{\ves_1}\parti_{2}.v_{(i'-1)\ves_1+ j
\ves_2+k\ves_3}\not=0\ \ \textrm{ for }j, k\in\mbb{Z},\, j\not=j'.
\end{equation}

Now we can define $\{v_{i'\ves_1+ j \ves_2+k\ves_3}\mid j\in\mbb{Z}\backslash\{0,j'\}, k\in\mbb{Z}\backslash\{0\}\}$ by
\begin{equation}\label{equ:9.4}
x^{\ves_1}\parti_{2}.v_{(i'-1)\ves_1+ j \ves_2+k\ves_3}=(a_2+j\ves_2(\parti_{2})) v_{i'\ves_1+ j \ves_2+k\ves_3}.
\end{equation}
This together with (\ref{equ:3.51}) and (\ref{equ:3.58}) indicate
\begin{equation}\label{equ:9.2}
x^{\ves_1}\parti_{2}.v_{(i'-1)\ves_1+ j \ves_2+k\ves_3}=(a_2+j\ves_2(\parti_{2}))
v_{i'\ves_1+ j \ves_2+k\ves_3} \quad\textrm{ for }j,k\in
\mathbb{Z},\,j\not=j'.
\end{equation}

Next we shall derive the other action in 1)--5).

 From applying (\ref{e:10}) and (\ref{equ:9.2}) to
 \begin{equation}
x^{ \pm
\ves_3}\parti_{2}.x^{ \ves_1}\parti_{2}.v_{(i'-1)\ves_1+ j \ves_2+k\ves_3}=x^{
\ves_1}\parti_{2}.x^{ \pm \ves_3}\parti_{2}.v_{(i'-1)\ves_1+ j \ves_2+k\ves_3},
\end{equation}
it can be immediately deduced
\begin{equation}\label{equ:8.4}
x^{ \pm \ves_3}\parti_{2}.v_{i'\ves_1+ j \ves_2+k\ves_3}=(a_2+j \ves_2(\parti_{2}))
v_{i'\ves_1+ j \ves_2+(k \pm 1)\ves_3}\quad  \textrm{ for }j,k\in
\mathbb{Z},\, j\not=j'.
\end{equation}

Recall that (\ref{equ:3.58}) and Claim 2 give $x^{-\ves_1}\parti_{2}.v_{i'\ves_1+ j\ves_2}=(a_2+j\ves_2(\parti_{2}))v_{(i'-1)\ves_1+ j\ves_2}$. Since
\begin{equation}
x^{-\ves_1}\parti_{2}.x^{ \pm\ves_3}\parti_{2}.v_{i'\ves_1+ j
\ves_2+k\ves_3} =x^{
\pm\ves_3}\parti_{2}.x^{-\ves_1}\parti_{2}.v_{i'\ves_1+ j
\ves_2+k\ves_3},
 \end{equation}
we get
\begin{equation}\label{equ:10.2}
x^{-\ves_1}\parti_{2}.v_{i'\ves_1+ j \ves_2+k\ves_3}=(a_2+j \ves_2(\parti_{2}))
v_{(i'-1)\ves_1+ j \ves_2+k \ves_3}\quad  \textrm{ for }j,k\in \mathbb{Z},\,
j\not=j'
\end{equation}
by (\ref{e:10}), (\ref{equ:8.4}) and induction on $k$.

As (\ref{e:11}) and (\ref{equ:8.33}) show
\begin{eqnarray}
x^{\ves_1+\ves_3}\parti_{2}.v_{(i'-2)\ves_1+ j \ves_2+k \ves_3}&=&\frac{1}{\ves_1(\parti_1)}[x^{\ves_3}\parti_{1},x^{\ves_1}\parti_{2}].v_{(i'-2)\ves_1+ j \ves_2+k \ves_3}\nonumber\\
&=&(a_2+j
\ves_2(\parti_{2}))v_{(i'-1)\ves_1+ j \ves_2+(k+1)\ves_3}\label{equ:9.3}
\end{eqnarray}
for $j,k\in\mbb{Z}$, using it and (\ref{equ:8.33}), (\ref{equ:9.2})
in $[x^{\ves_1+\ves_3}\parti_{2},x^{\ves_1}\parti_{2}].v_{(i'-2)\ves_1+ j
\ves_2+k \ves_3}=0$, we get
\begin{equation}\label{equ:8.1}
x^{\ves_1+\ves_3}\parti_{2}.v_{(i'-1)\ves_1+ j \ves_2+k \ves_3}=(a_2+j
\ves_2(\parti_{2}))v_{i'\ves_1+ j \ves_2+(k+1)\ves_3}\qquad   \textrm{ for
}j,k\in \mathbb{Z},\, j\not=j'.
\end{equation}

Inserting (\ref{e:11}), (\ref{equ:9.2}) and (\ref{equ:8.1})
into
 \begin{equation}
 [x^{\ves_3}\parti_{1},x^{\ves_1}\parti_{2}].v_{(i'-1)\ves_1+ j\ves_2+k\ves_3}=
 \ves_1(\parti_{1})x^{\ves_1+\ves_3}\parti_{2}.v_{(i'-1)\ves_1+j \ves_2+k\ves_3}
 \end{equation}
 and
 \begin{equation}
[x^{-\ves_3}\parti_{1},x^{\ves_1+\ves_3}\parti_{2}].v_{(i'-1)\ves_1+ j \ves_2+(k-1)\ves_3}=\ves_1(\parti_{1})x^{\ves_1}\parti_{2}.v_{(i'-1)\ves_1+ j \ves_2+(k-1)\ves_3}
\end{equation}
respectively, we have
\begin{equation}\label{equ:10.4}
x^{\pm\ves_3}\parti_{1}.v_{i'\ves_1+ j \ves_2+k\ves_3}=(a_{1}+i'
\ves_1(\parti_{1}))v_{i'\ves_1+ j \ves_2+(k\pm 1)\ves_3}=0  \ \textrm{ for }j,k\in \mathbb{Z},\,  j\not=j'.
\end{equation}

From (\ref{equ:3.55}) and Claim 2 we see that $x^{\ves_1}\parti_{3}.v_{(i'-1)\ves_1+ j \ves_2}=a_{3}v_{i'\ves_1+ j
\ves_2}$ for $j\in \mathbb{Z}$. Inserting (\ref{e:10}),
(\ref{equ:8.4}) and (\ref{equ:8.1}) into
\begin{equation}
[x^{\ves_1}\parti_{3},x^{\ves_3}\parti_{2}].v_{(i'-1)\ves_1+
j\ves_2+k\ves_3}=\ves_3(\parti_{3})x^{\ves_1+\ves_3}\parti_{2}.v_{(i'-1)\ves_1+
j\ves_2+k\ves_3},
\end{equation}
we find
\begin{equation}\label{equ:8.2}
x^{\ves_1}\parti_{3}.v_{(i'-1)\ves_1+ j \ves_2+k\ves_3}=(a_{3}+k
\ves_3(\parti_{3}))v_{i'\ves_1+ j \ves_2+k\ves_3}\qquad  \textrm{ for }j,k\in
\mathbb{Z},\, j\not=j'
\end{equation}
by induction on $k$.

Applying (\ref{e:10}), (\ref{equ:8.36}), (\ref{equ:9.3}) and
(\ref{equ:8.1}) to
\begin{equation}
[x^{-\ves_1}\parti_{3},x^{\ves_1+\ves_3}\parti_{2}].v_{(i'-1)\ves_1+ j \ves_2+(k-1)\ves_3}
=\ves_3(\parti_{3})x^{\ves_3}\parti_{2}.v_{(i'-1)\ves_1+ j \ves_2+(k-1)\ves_3},
\end{equation}
we obtain
\begin{equation}\label{equ:10.3}
x^{-\ves_1}\parti_{3}.v_{i'\ves_1+ j \ves_2+k\ves_3}=(a_{3}+k
\ves_3(\parti_{3}))v_{(i'-1)\ves_1+ j \ves_2+k\ves_3}\quad  \textrm{ for }
j,k\in \mathbb{Z},\, j\not=j'.
\end{equation}

Using (\ref{equ:10.4}) in
\begin{equation}
x^{\pm\ves_2}\parti_{1}.v_{i'\ves_1+ j \ves_2+k\ves_3}=\pm\frac{1}{\ves_2(\parti_{2})\ves_3(\parti_{3})}[[x^{\pm\ves_2}\parti_{3},x^{\ves_3}\parti_{2}],x^{-\ves_3}\parti_{1}].v_{i'\ves_1+ j \ves_2+k\ves_3},
\end{equation}
we have
\begin{equation}\label{e:1.1}
x^{\ves_2}\parti_{1}.v_{i'\ves_1+ j \ves_2+k\ves_3}=0=(a_{1}+i'
\ves_1(\parti_{1}))v_{i'\ves_1+ (j+1) \ves_2+k\ves_3} \quad \textrm{ for }j,k\in
\mathbb{Z},\, j\not=j', j'-1,
\end{equation}
\begin{equation}\label{e:1.2}
x^{-\ves_2}\parti_{1}.v_{i'\ves_1+ j \ves_2+k\ves_3}=0=(a_{1}+i'
\ves_1(\parti_{1}))v_{i'\ves_1+ (j+1) \ves_2+k\ves_3} \quad  \textrm{ for
}j,k\in \mathbb{Z},\, j\not=j', j'+1.
\end{equation}

From using (\ref{e:12}) and (\ref{equ:8.2}) in
\begin{equation}
x^{\pm\ves_2}\parti_{3}.x^{\ves_1}\parti_{3}.v_{(i'-1)\ves_1+ j \ves_2+k\ves_3}=x^{\ves_1}\parti_{3}.x^{\pm\ves_2}\parti_{3}.v_{(i'-1)\ves_1+ j \ves_2+k\ves_3},
\end{equation}
it can be deduced
\begin{equation}\label{equ:8.7}
x^{\ves_2}\parti_{3}.v_{i'\ves_1+ j \ves_2+k\ves_3}=(a_{3}+k
\ves_3(\parti_{3}))v_{i'\ves_1+ (j+1) \ves_2+k\ves_3}\   \textrm{ for }
j\in \mbb{Z}\backslash\{j', j'-1\},\, k\in \mbb{Z}\backslash\{k'\},
\end{equation}
\begin{equation}\label{equ:10.1}
x^{-\ves_2}\parti_{3}.v_{i'\ves_1+ j \ves_2+k\ves_3}=(a_{3}+k
\ves_3(\parti_{3}))v_{i'\ves_1+ (j-1) \ves_2+k\ves_3}\   \textrm{ for }
j\in \mbb{Z}\backslash\{j', j'+1\}, \, k\in \mbb{Z}\backslash\{k'\}.
\end{equation}
On the other hand, from applying (\ref{equ:8.4}), (\ref{equ:8.7}) and (\ref{equ:10.1}) to
\begin{equation}
x^{\pm\ves_2}\parti_{3}.v_{i'\ves_1+ j \ves_2+k'\ves_3}=\frac{1}{\ves_2^{2}(\parti_{2})}[[x^{\pm\ves_2}\parti_{3},x^{\ves_3}\parti_{2}],x^{-\ves_3}\parti_{2}].v_{i'\ves_1+ j \ves_2+k'\ves_3},
\end{equation}
it can be derived
\begin{equation}\label{equ:10.5}
x^{\ves_2}\parti_{3}.v_{i'\ves_1+ j \ves_2+k'\ves_3}=0=(a_{3}+k'
\ves_3(\parti_{3}))v_{i'\ves_1+ (j+1) \ves_2+k'\ves_3}\   \textrm{ for }
j\in \mbb{Z}\backslash\{j', j'-1\},
\end{equation}
\begin{equation}\label{equ:10.6}
x^{-\ves_2}\parti_{3}.v_{i'\ves_1+ j \ves_2+k'\ves_3}=0=(a_{3}+k'
\ves_3(\parti_{3}))v_{i'\ves_1+ (j-1) \ves_2+k'\ves_3}\  \textrm{ for }
j\in \mbb{Z}\backslash\{j', j'+1\},
\end{equation}
which coincide with (\ref{equ:8.7}) and (\ref{equ:10.1}),
respectively.\vspace{0.3cm}

{\it Thirdly, under the condition that $j'$ does exist, we define $\{v_{i'\ves_1+ j' \ves_2+k\ves_3}\mid k\in\mbb{Z}\backslash\{0,k'\}\}$, and determine}

  1) $\textrm{the action of } x^{\ves_1}\parti_{2} \textrm{ and } x^{\ves_1}\parti_{3} \textrm{ on } \bigoplus_{k\in\mathbb{Z}\backslash\{k'\}}M_{(i'-1)\ves_1+j'\ves_2+k\ves_3}$,\vspace{0.1cm}

  2) $\textrm{the action of } x^{-\ves_1}\parti_{2} \textrm{ and } x^{-\ves_1}\parti_{3} \textrm{ on } \bigoplus_{k\in\mathbb{Z}\backslash\{k'\}}M_{i'\ves_1+j'\ves_2+k\ves_3}$,\vspace{0.1cm}

  3) $\textrm{the action of } x^{\ves_2}\parti_{1} \textrm{ and } x^{\ves_2}\parti_{3} \textrm{ on } \bigoplus_{k\in\mathbb{Z}\backslash\{k'\}}M_{i'\ves_1+j\ves_2+k\ves_3} \textrm{ for } j=j',j'-1$,\vspace{0.1cm}

  4) $\textrm{the action of } x^{-\ves_2}\parti_{1} \textrm{ and } x^{-\ves_2}\parti_{3} \textrm{ on } \bigoplus_{k\in\mathbb{Z}\backslash\{k'\}}M_{i'\ves_1+j\ves_2+k\ves_3} \textrm{ for } j=j',j'+1$,\vspace{0.1cm}

 5) $\textrm{the action of } x^{\ves_3}\parti_{1} \textrm{ and } x^{\ves_3}\parti_{2} \textrm{ on } \bigoplus_{k\in\mathbb{Z}\backslash\{k',k'-1\}}M_{i'\ves_1+j'\ves_2+k\ves_3}$,\vspace{0.1cm}

 6) $\textrm{the action of } x^{-\ves_3}\parti_{1} \textrm{ and } x^{-\ves_3}\parti_{2} \textrm{ on } \bigoplus_{k\in\mathbb{Z}\backslash\{k',k'+1\}}M_{i'\ves_1+j'\ves_2+k\ves_3}$. \vspace{0.2cm}

If $j'$ does not exist, this part can be skipped.

Under the condition that $j'$ exists, we want to define $v_{i'\ves_1+ j'\ves_2+k\ves_3}$  for $k\in\mathbb{Z}\backslash \{0, k'\}$. As (\ref{e:12}) and (\ref{equ:8.2}) show
\begin{equation}
x^{-\ves_2}\parti_{3}.x^{\ves_1}\parti_{3}.v_{(i'-1)\ves_1+ j'\ves_2+k\ves_3}=x^{\ves_1}\parti_{3}.x^{-\ves_2}\parti_{3}.v_{(i'-1)\ves_1+j' \ves_2+k\ves_3}\not=0 \, \textrm{ for }k\in\mathbb{Z}\backslash \{ k'\},
\end{equation}
we have
\begin{equation}
x^{\ves_1}\parti_{3}.v_{(i'-1)\ves_1+ j'\ves_2+k\ves_3}\not=0 \ \textrm{ for }
 k\in\mathbb{Z}\backslash \{ k'\}.
 \end{equation}
So we can define $v_{i'\ves_1+ j'\ves_2+k\ves_3}$'s by
\begin{equation}\label{eq:6}
x^{\ves_1}\parti_{3}.v_{(i'-1)\ves_1+ j'\ves_2+k\ves_3}=(a_{3}+k\ves_3(\parti_{3})) v_{i'\ves_1+ j' \ves_2+k\ves_3} \ \textrm{ for }k\in\mathbb{Z}\backslash \{0, k'\}.
\end{equation}
This and (\ref{equ:3.55}) together with Claim 2 show
\begin{equation}\label{equ:10.7}
x^{\ves_1}\parti_{3}.v_{(i'-1)\ves_1+ j'\ves_2+k\ves_3}=(a_{3}+k\ves_3(\parti_{3})) v_{i'\ves_1+ j' \ves_2+k\ves_3} \  \textrm{ for } k\in\mathbb{Z}\backslash\{k'\}.
\end{equation}

Next we shall derive the action of the other elements.

Applying
(\ref{e:12}), (\ref{equ:8.2}) and (\ref{equ:10.7}) to
\begin{equation}
x^{\pm\ves_2}\parti_{3}.x^{\ves_1}\parti_{3}.v_{(i'-1)\ves_1+
j\ves_2+k\ves_3}=x^{\ves_1}\parti_{3}.x^{\pm\ves_2}\parti_{3}.v_{(i'-1)\ves_1+
j\ves_2+k\ves_3},
\end{equation}
we obtain
\begin{equation}\label{equ:8.6}
x^{\ves_2}\parti_{3}.v_{i'\ves_1+ j\ves_2+k\ves_3}=(a_{3}+k
\ves_3(\parti_{3})) v_{i'\ves_1+ (j +1) \ves_2+k\ves_3}  \textrm{ for } j\in\mathbb{Z}\backslash \{j',j'-1\},\, k\in\mathbb{Z}\backslash \{ k'\},
\end{equation}
\begin{equation}\label{equ:10.9}
x^{-\ves_2}\parti_{3}.v_{i'\ves_1+ j\ves_2+k\ves_3}=(a_{3}+k
\ves_3(\parti_{3})) v_{i'\ves_1+ (j -1) \ves_2+k\ves_3}  \textrm{ for } j\in\mathbb{Z}\backslash \{j',j'+1\},\, k\in\mathbb{Z}\backslash \{ k'\}.
\end{equation}

From using (\ref{e:12}), (\ref{equ:10.3}) and (\ref{equ:8.6})
in
\begin{equation}
x^{-\ves_1}\parti_{3}.x^{\ves_2}\parti_{3}.v_{i'\ves_1+ (j'-1)\ves_2+k\ves_3}=x^{\ves_2}\parti_{3}.x^{-\ves_1}\parti_{3}.v_{i'\ves_1+ (j'-1)\ves_2+k\ves_3},
\end{equation}
it can be deduced
\begin{equation}\label{equ:10.8}
x^{-\ves_1}\parti_{3}.v_{i'\ves_1+ j'\ves_2+k\ves_3}=(a_{3}+k \ves_3(\parti_{3}))
v_{(i'-1)\ves_1+ j' \ves_2+k\ves_3}\  \textrm{ for }k\in\mathbb{Z}\backslash \{ k'\}.
\end{equation}

Since (\ref{e:13}) and (\ref{equ:8.35}) indicate
\begin{eqnarray}
x^{\ves_1+\ves_2}\parti_{3}.v_{ (i'-2)\ves_1+ j\ves_2+k\ves_3}&=&\frac{1}{\ves_1(\parti_{1})}[x^{\ves_2}\parti_{1},x^{\ves_1}\parti_{3}].v_{(i'-2)\ves_1+ j\ves_2+k\ves_3}\nonumber\\
&=&(a_{3}+k \ves_3(\parti_{3})) v_{(i'-1)\ves_1+ (j+1) \ves_2+k\ves_3} \label{equ:11.3}
\end{eqnarray}
for $j,k\in\mbb{Z}$, applying it and (\ref{equ:8.35}), (\ref{equ:8.2}), (\ref{equ:10.7}) to
\begin{equation}\label{equ:11.4}
x^{\ves_1+\ves_2}\parti_{3}.x^{\ves_1}\parti_{3}.v_{ (i'-2)\ves_1+ j\ves_2+k\ves_3}=x^{\ves_1}\parti_{3}.x^{\ves_1+\ves_2}\parti_{3}.v_{ (i'-2)\ves_1+ j\ves_2+k\ves_3},
\end{equation}
we can derive
\begin{equation}\label{equ:11.5}
x^{\ves_1+\ves_2}\parti_{3}.v_{ (i'-1)\ves_1+ j\ves_2+k\ves_3}=(a_{3}+k
\ves_3(\parti_{3})) v_{i'\ves_1+ (j+1) \ves_2+k\ves_3}\  \textrm{ for }j,k\in\mbb{Z},\,k\not=k'.
\end{equation}

Then using (\ref{e:12}), (\ref{equ:9.2}), (\ref{equ:8.6}) and
(\ref{equ:11.5}) in
\begin{equation}
[x^{\ves_1}\parti_{2},x^{\ves_2}\parti_{3}].v_{(i'-1)\ves_1+ (j'-1)\ves_2+k\ves_3}=\ves_2(\parti_{2})x^{\ves_1+\ves_2}\parti_{3}.v_{(i'-1)\ves_1+ (j'-1)\ves_2+k\ves_3},
\end{equation}
we get
\begin{equation}
x^{\ves_1}\parti_{2}.v_{ (i'-1)\ves_1+ j'\ves_2+k\ves_3}=(a_{2}+j'
\ves_2(\parti_{2}))v_{ i'\ves_1+ j'\ves_2+k\ves_3}=0 \ \textrm{ for }
k\in\mathbb{Z}\backslash \{k'\}.
\end{equation}
Moreover, using (\ref{e:12}), (\ref{equ:8.34}), (\ref{equ:11.3})
and (\ref{equ:11.5}) in
\begin{equation}
[x^{-\ves_1}\parti_{2},x^{\ves_1+\ves_2}\parti_{3}].v_{(i'-1)\ves_1+ (j'-1)\ves_2+k\ves_3}=\ves_2(\parti_{2})x^{\ves_2}\parti_{3}.v_{(i'-1)\ves_1+ (j'-1)\ves_2+k\ves_3},
\end{equation}
we obtain
\begin{equation}
x^{-\ves_1}\parti_{2}.v_{ i'\ves_1+ j'\ves_2+k\ves_3}=(a_{2}+j'
\ves_2(\parti_{2}))v_{ (i'-1)\ves_1+ j'\ves_2+k\ves_3}=0 \ \textrm{ for }
k\in\mathbb{Z}\backslash \{k'\}.
\end{equation}

From inserting (\ref{e:13}), (\ref{equ:8.2}), (\ref{equ:10.7})
 and (\ref{equ:11.5}) into
\begin{equation}
[x^{\ves_2}\parti_{1},x^{\ves_1}\parti_{3}].v_{(i'-1)\ves_1+ j\ves_2+k\ves_3}=\ves_1(\parti_{1})x^{\ves_1+\ves_2}\parti_{3}.v_{(i'-1)\ves_1+ j\ves_2+k\ves_3} ,
\end{equation}
\begin{equation}
[x^{-\ves_2}\parti_{1},x^{\ves_1+\ves_2}\parti_{3}].v_{(i'-1)\ves_1+ j\ves_2+k\ves_3}=\ves_1(\parti_{1})x^{\ves_1}\parti_{3}.v_{(i'-1)\ves_1+ j\ves_2+k\ves_3},
\end{equation}
 it can be respectively deduced
\begin{equation}
x^{\ves_2}\parti_{1}.v_{ i'\ves_1+ j\ves_2+k\ves_3}=(a_{1}+i'
\ves_1(\parti_{1}))v_{ i'\ves_1+ (j+1)\ves_2+k\ves_3}=0
\textrm{ for }j\in\{j',j'-1\},\, k\not=k',
\end{equation}
\begin{equation}
x^{-\ves_2}\parti_{1}.v_{ i'\ves_1+ j\ves_2+k\ves_3}=(a_{1}+i'
\ves_1(\parti_{1}))v_{ i'\ves_1+ (j-1)\ves_2+k\ves_3}=0   \textrm{ for
}j\in\{j',j'+1\},\, k\not=k'.
\end{equation}

Substituting (\ref{e:10}), (\ref{equ:10.7}) and (\ref{equ:10.8})
into
\begin{equation}
[x^{-\ves_1}\parti_{3},[x^{\ves_1}\parti_{3},x^{\pm\ves_3}\parti_{2}]].v_{(i'-1)\ves_1+ j'\ves_2+k\ves_3}=\ves_3^{2}(\parti_{3})x^{\pm\ves_3}\parti_{2}.v_{(i'-1)\ves_1+ j'\ves_2+k\ves_3},
\end{equation}
 we get
\begin{equation}\label{equ:11.1}
x^{\ves_3}\parti_{2}.v_{ i'\ves_1+ j'  \ves_2+k\ves_3}=0=(a_2+j'
\ves_2(\parti_{2}))v_{i'\ves_1+ j' \ves_2+(k+1)\ves_3}\  \textrm{ for }
k\in\mathbb{Z}\backslash \{ k', k'-1\},
\end{equation}
\begin{equation}\label{equ:11.2}
x^{-\ves_3}\parti_{2}.v_{ i'\ves_1+ j'  \ves_2+k\ves_3}=0=(a_2+j'
\ves_2(\parti_{2}))v_{i'\ves_1+ j' \ves_2+(k-1)\ves_3}\   \textrm{ for }
k\in\mathbb{Z}\backslash \{ k', k'+1\}.
\end{equation}

Applying (\ref{equ:10.4}), (\ref{equ:8.6}) and (\ref{equ:10.9})
to
\begin{equation}
[x^{-\ves_2}\parti_{3},[x^{\ves_2}\parti_{3},x^{\pm\ves_3}\parti_{1}]].v_{i'\ves_1+ (j'-1)\ves_2+k\ves_3}=\ves_3^{2}(\parti_{3})x^{\pm\ves_3}\parti_{1}.v_{i'\ves_1+ (j'-1)\ves_2+k\ves_3},
\end{equation}
 we obtain
\begin{equation}\label{equ:4.109}
 x^{\ves_3}\parti_{1}.v_{ i'\ves_1+ j'  \ves_2+k\ves_3}=0=(a_{1}+i' \ves_1(\parti_{1}))
v_{i'\ves_1+ j' \ves_2+(k+1)\ves_3}\   \textrm{ for }
k\in\mathbb{Z}\backslash \{ k', k'-1\},
\end{equation}
\begin{equation}\label{equ:4.110}
 x^{-\ves_3}\parti_{1}.v_{ i'\ves_1+ j'  \ves_2+k\ves_3}=0=(a_{1}+i' \ves_1(\parti_{1}))
v_{i'\ves_1+ j' \ves_2+(k-1)\ves_3}\   \textrm{ for }
k\in\mathbb{Z}\backslash \{ k', k'+1\}.
\end{equation}\vspace{0.1cm}

{\it Fourthly, we want to derive}

  1)  $\textrm{the action of } x^{\ves_1}\parti_{2} \textrm{ and } x^{\ves_1}\parti_{3} \textrm{ on } \bigoplus_{(j',k')\not=(j,k) \in \mbb{Z}^2}M_{i'\ves_1+j\ves_2+k\ves_3}$, \vspace{0.1cm}

  2) $\textrm{the action of } x^{-\ves_1}\parti_{2} \textrm{ and } x^{-\ves_1}\parti_{3} \textrm{ on } \bigoplus_{(j',k')\not=(j,k) \in \mbb{Z}^2}M_{(i'+1)\ves_1+j\ves_2+k\ves_3}$.\vspace{0.2cm}

Taking $i=i'+1$ in (\ref{e:13})--(\ref{e:12}), we have
\begin{equation}
x^{ \pm \ves_3}\parti_{1}.v_{ (i'+1)\ves_1+ j
\ves_2+k\ves_3}=\ves_1(\parti_{1}) v_{(i'+1)\ves_1+ j \ves_2+(k \pm 1)\ves_3},
\end{equation}
\begin{equation}
 x^{\pm
\ves_2}\parti_{1}.v_{ (i'+1)\ves_1+ j \ves_2+k\ves_3}=\ves_1(\parti_{1})
v_{ (i'+1)\ves_1+(j \pm 1) \ves_2+k\ves_3},
\end{equation}
\begin{equation}\label{equ:8.3}
 x^{ \pm \ves_3}\parti_{2}.v_{(i'+1)\ves_1+ j \ves_2+k\ves_3}=(a_2+j
\ves_2(\parti_{2})) v_{(i'+1)\ves_1+ j \ves_2+(k \pm 1)\ves_3},
\end{equation}
\begin{equation}\label{equ:8.5}
  x^{\pm \ves_2}\parti_{3}.v_{ (i'+1)\ves_1+ j \ves_2+k\ves_3}=(a_{3}+ k
\ves_3(\parti_{3})) v_{ (i'+1)\ves_1+(j \pm 1) \ves_2+k\ves_3}
\end{equation}
for $j,k \in \mathbb{Z}$.

Recall that (\ref{equ:3.58}) and Claim 2 give $x^{ \ves_1}\parti_{2}.v_{i'\ves_1+ j \ves_2}=(a_2+j
\ves_2(\parti_{2})) v_{(i'+1)\ves_1+ j \ves_2}$ for $j\in\mbb{Z}$. Since
\begin{equation}
x^{ \ves_1}\parti_{2}.x^{\pm\ves_3}\parti_{2}.v_{i'\ves_1+ j \ves_2+k\ves_3}=x^{\pm\ves_3}\parti_{2}.x^{ \ves_1}\parti_{2}.v_{i'\ves_1+ j \ves_2+k\ves_3},
\end{equation}
we obtain
\begin{equation}\label{eq:7}
x^{ \ves_1}\parti_{2}.v_{  i'\ves_1+ j \ves_2+k\ves_3}=(a_2+j \ves_2(\parti_{2}))
v_{ (i'+1)\ves_1+ j \ves_2+k\ves_3}\quad \textrm{ for } j,k \in \mathbb{Z},\,
j\not=j'
\end{equation}
by (\ref{equ:8.4}), (\ref{equ:8.3}) and induction on $k$. On the other hand, since (\ref{equ:8.6}), (\ref{equ:10.9}), (\ref{equ:8.5}) and (\ref{eq:7}) lead to
\begin{eqnarray}
0&=&[x^{ -\ves_2}\parti_{3},[x^{
\ves_2}\parti_{3},x^{ \ves_1}\parti_{2}]].v_{  i'\ves_1 +j' \ves_2+k\ves_3}\nonumber\\
&=&x^{ -\ves_2}\parti_{3}.x^{
\ves_2}\parti_{3}.x^{ \ves_1}\parti_{2}.v_{  i'\ves_1 +j' \ves_2+k\ves_3}+x^{ \ves_1}\parti_{2}.x^{ \ves_2}\parti_{3}.x^{-\ves_2}\parti_{3}.v_{  i'\ves_1 +j' \ves_2+k\ves_3}
\end{eqnarray}
for $k\in\mbb{Z}\backslash\{ k'\}$, we find
\begin{equation}
x^{ \ves_1}\parti_{2}.v_{ i'\ves_1+ j' \ves_2+k\ves_3}=0=(a_2+j'
\ves_2(\parti_{2})) v_{ (i'+1)\ves_1+ j' \ves_2+k\ves_3}\quad \textrm{for } k\in\mbb{Z}\backslash\{ k'\}
\end{equation}
by (\ref{equ:8.6}), (\ref{equ:10.9}) and (\ref{equ:8.5}).

Similarly, in analogy with (\ref{eq:7}), we get
\begin{equation}
x^{ -\ves_1}\parti_{2}.v_{ (i'+1)\ves_1+ j \ves_2+k\ves_3}=(a_2+j
\ves_2(\parti_{2})) v_{ i'\ves_1+ j \ves_2+k\ves_3}\quad \textrm{ for } j,k \in
\mathbb{Z},\, j\not=j'.
\end{equation}
By the fact $[x^{ -\ves_2}\parti_{3},[x^{ -\ves_1}\parti_{2},x^{
\ves_2}\parti_{3}]].v_{ (i'+1)\ves_1 +j' \ves_2+k\ves_3}=0$, we derive
\begin{equation}
x^{ -\ves_1}\parti_{2}.v_{ (i'+1)\ves_1+ j' \ves_2+k\ves_3}=0=(a_2+j'
\ves_2(\parti_{2})) v_{ i'\ves_1+ j' \ves_2+k\ves_3} \quad \textrm{for
} k\in\mbb{Z}\backslash\{ k'\}.
\end{equation}

Recall that (\ref{equ:3.54}) and Claim 2 give $x^{\ves_1}\parti_{3}.v_{ i'\ves_1+k\ves_3}=(a_{3}+k\ves_3(\parti_{3})) v_{ (i'+1)\ves_1+k\ves_3}$ for $k\in\mbb{Z}$. Since
\begin{equation}
x^{ \ves_1}\parti_{3}.x^{ \pm\ves_2}\parti_{3}.v_{
 i'\ves_1+j\ves_2+k\ves_3}=x^{ \pm\ves_2}\parti_{3}.x^{ \ves_1}\parti_{3}.v_{
 i'\ves_1+j\ves_2+k\ves_3},
\end{equation}
we obtain
\begin{equation}\label{eq:8}
x^{\ves_1}\parti_{3}.v_{ i'\ves_1+ j\ves_2+k\ves_3}=(a_{3}+k \ves_3(\parti_{3}))
v_{ (i'+1)\ves_1+ j \ves_2+k\ves_3}\quad \textrm{ for } j,k \in \mathbb{Z},\,
k\not=k'
\end{equation}
by (\ref{equ:8.7}), (\ref{equ:10.1}), (\ref{equ:8.6}),
(\ref{equ:10.9}), (\ref{equ:8.5}) and induction on $j$. Since (\ref{equ:8.4}), (\ref{equ:8.3}) and (\ref{eq:8}) lead to
\begin{eqnarray}
0&=&[x^{ -\ves_3}\parti_{2},[x^{
\ves_3}\parti_{2},x^{ \ves_1}\parti_{3}]].v_{ i'\ves_1 +j \ves_2+k' \ves_3 }\nonumber\\
&=&x^{ -\ves_3}\parti_{2}.x^{\ves_3}\parti_{2}.x^{ \ves_1}\parti_{3}.v_{ i'\ves_1 +j \ves_2+k' \ves_3 }+x^{ \ves_1}\parti_{3}.x^{\ves_3}\parti_{2}.x^{-\ves_3}\parti_{2}.v_{ i'\ves_1 +j \ves_2+k' \ves_3 }
\end{eqnarray}
for $j\in\mbb{Z}\backslash\{j'\}$, it can be deduced from (\ref{equ:8.4}) and (\ref{equ:8.3}) that
\begin{equation}
x^{ \ves_1}\parti_{3}.v_{i'\ves_1+ j \ves_2+k' \ves_3 }=0=(a_{3}+k'\ves_3(\parti_{3}))
v_{ (i'+1)\ves_1+ j \ves_2+k'\ves_3} \  \textrm{ for } j\in\mbb{Z}\backslash\{j'\}.
\end{equation}

Similarly, in analogy with (\ref{eq:8}), we get
\begin{equation}
x^{-\ves_1}\parti_{3}.v_{ (i'+1)\ves_1+ j\ves_2+k\ves_3}=(a_{3}+k
\ves_3(\parti_{3})) v_{ i'\ves_1+ j \ves_2+k\ves_3}\quad\textrm{ for } j,k \in
\mathbb{Z},\, k\not=k',
\end{equation}
and by the fact $[x^{ -\ves_3}\parti_{2},[x^{ -\ves_1}\parti_{3},x^{
\ves_3}\parti_{2}]].v_{(i'+1)\ves_1 +j \ves_2+k' \ves_3 }=0$, we deduce
\begin{equation}
x^{ -\ves_1}\parti_{3}.v_{(i'+1)\ves_1+ j \ves_2 +k' \ves_3 }=0=(a_{3}+k'
\ves_3(\parti_{3})) v_{ i'\ves_1+ j \ves_2+k'\ves_3} \quad
\textrm{for } j \in\mbb{Z}\backslash\{j'\}.
\end{equation}\vspace{0.1cm}

{\it Finally, we completely specify the action of ${\cal X}$ on $M$ and give the conclusion of this case.}\vspace{0.2cm}

By now, we have determined $\{v_{i\ves_1+ j \ves_2 +k \ves_3}\mid (i',j',k')\not=(i,j,k)\in\mbb{Z}^3 \}$ in this case.
Define $\mu\in D ^{\ast} $ by $\mu(\parti_{1})=a_{1}$,
$\mu(\parti_{2})=a_2$ and $\mu(\parti_{3})=a_{3}$. Write $\zeta=i'\ves_1 +j' \ves_2+k'
\ves_3$ for short. Summing up the acting relations mentioned above, we can write them uniformly by
\begin{equation}\label{e:7}
x^{\es}\parti.v_{\be}=(\be+\mu)(\parti)v_{\be+\es}
\end{equation}
for $\es \in \{ \pm\ves_1, \pm\ves_2,\pm\ves_3\}$, $\parti \in \kn{\es}\cap\{\parti_{1}, \parti_{2},
\parti_{3}\}$ and $\be\in\G\backslash\{\zeta, \,\zeta-\es\}$.

If $j'$ or $k'$ does not exist, i.e., $\zeta$ does not exist, then (\ref{e:7}) shows the general action of ${\cal X}$. Since ${\cal X}$ generates ${\cal S}(\G, D)$ (c.f. (\ref{eq:31})), we deduce
\begin{equation}
M \simeq \mathscr{M}_{\mu}.
\end{equation}

If both $j'$ and $k'$ exist, then $\zeta=-\mu$ and we have a few acting relations left to be specified, say,
\begin{equation}\label{eq:15}
\textrm{how }x^{\es}\parti \textrm{ act on } M_{\zeta}\textrm{ and on } M_{\zeta-\es} \textrm{ for }\es \in \{ \pm\ves_1, \pm\ves_2,\pm\ves_3
\}\textrm{ and }\parti \in \kn{\es}\cap\{\parti_{1}, \parti_{2},
\parti_{3}\}.
\end{equation}
The action in (\ref{eq:15}) is divided into three subcases.

The first subcase is that $x^{\es}\parti.M_{\zeta}=\{0\}$ and
$x^{\es}\parti.v_{\zeta-\es}=0$ for all $\es \in \{ \pm\ves_1, \pm\ves_2,\pm\ves_3
\}$ and $\parti \in \kn{\es}\cap\{\parti_{1}, \parti_{2},
\parti_{3}\}$. This and (\ref{e:7}) give the general action of ${\cal X}$. Since ${\cal X}$ generates ${\cal S}(\G, D)$ (c.f. (\ref{eq:31})), we deduce
\begin{equation}
M \simeq
\mathscr{M}_{\mu}.\end{equation}

The second subcase is that $x^{\es}\parti.M_{\zeta}\not=\{0\}$ for some $\es \in
\{\pm\ves_1, \pm\ves_2, \pm\ves_3\}$ and some $\parti \in
\kn{\es}\cap\{\parti_{1},
\parti_{2}, \parti_{3}\}$. Assume
$x^{\ves_3}\parti_{1}.M_{\zeta}\not=\{0\}$; the other cases can be proved
similarly. Pick a nonzero vector  $v_{\zeta}$ of $M_{\zeta}$. Write
\begin{equation}\label{e:15}
x^{\ves_3}\parti_{1}.v_{\zeta}=\hat{a}_{1}v_{\zeta+\ves_3},
\ \
x^{\ves_1}\parti_{2}.v_{\zeta}=\hat{a}_{2}v_{\zeta+\ves_1},
\ \
x^{\ves_2}\parti_{3}.v_{\zeta}=\hat{a}_{3}v_{\zeta+\ves_2}
\end{equation}
with $\hat{a}_{1}, \hat{a}_{2}, \hat{a}_{3} \in \mathbb{F}$, where $\hat{a}_{1}\not=0$.

By (\ref{e:7}), we have
\begin{equation}
x^{\ves_3}\parti_{1}.(x^{-\ves_3}\parti_{l}.v_{\zeta+\ves_3})=x^{-\ves_3}\parti_{l}.(x^{\ves_3}\parti_{1}.v_{\zeta+\ves_3})=0 \ \textrm{ for }l \in \{1, 2\},
\end{equation}
which implies
\begin{equation}\label{e:14}
x^{-\ves_3}\parti_{l}.v_{\zeta+\ves_3}=0 \quad \textrm{for }l \in \{1, 2\}.
\end{equation}
Moreover, (\ref{e:7}) yields
\begin{equation}
x^{\ves_3}\parti_{1}.(x^{\pm\ves_2}\parti_{1}.v_{\zeta\mp\ves_2})
=x^{\pm\ves_2}\parti_{1}.(x^{\ves_3}\parti_{1}.v_{\zeta\mp\ves_2})=0,
\end{equation}
which shows
\begin{equation}
x^{\pm\ves_2}\parti_{1}.v_{\zeta\mp\ves_2}=0.
\end{equation}

Since (\ref{e:14}) leads to
\begin{equation}
[x^{-\ves_1}\parti_{2}, [x^{-\ves_3}\parti_{1},
x^{\ves_1}\parti_{3}]].v_{\zeta+\ves_3}=\ves_3(\parti_{3})\ves_1(\parti_{1})x^{-\ves_3}\parti_{2}.v_{\zeta+\ves_3}=0,
\end{equation}
from (\ref{e:7}) and (\ref{e:14}) we derive
\begin{eqnarray}
x^{-\ves_1}\parti_{2}.v_{\zeta+\ves_1} & = & \frac{1}{\ves_3(\parti_{3})\ves_1(\parti_{1})}x^{-\ves_1}\parti_{2}.(x^{-\ves_3}\parti_{1}.x^{\ves_1}\parti_{3} .v_{\zeta+\ves_3})\nonumber\\
& = &\frac{1}{\ves_3(\parti_{3})\ves_1(\parti_{1})}[x^{-\ves_1}\parti_{2}, [x^{-\ves_3}\parti_{1},
x^{\ves_1}\parti_{3}]].v_{\zeta+\ves_3}=0.
\end{eqnarray}
Similarly, from
\begin{equation}
[x^{\ves_1}\parti_{2}, [x^{-\ves_3}\parti_{1},
x^{-\ves_1}\parti_{3}]].v_{\zeta+\ves_3}=-\ves_3(\parti_{3})\ves_1(\parti_{1})x^{-\ves_3}\parti_{2}.v_{\zeta+\ves_3}=0
\end{equation}
it can be deduced
\begin{equation}
x^{\ves_1}\parti_{2}.v_{\zeta-\ves_1}=0.
\end{equation}

Since
\begin{eqnarray}
& & x^{\ves_2}\parti_{1}.v_{\zeta}\nonumber\\
&=&\frac{1}{\ves_3(\parti_{3})\ves_2(\parti_{2})}[x^{-\ves_3}\parti_{2}, [x^{\ves_2}\parti_{3}, x^{\ves_3}\parti_{1}]].v_{\zeta}\nonumber\\
&=&\frac{1}{\ves_3(\parti_{3})\ves_2(\parti_{2})}(\hat{a}_{1}\ves_3(\parti_{3})\ves_2(\parti_{2})v_{\zeta+\ves_2}-x^{\ves_2}\parti_{3}.x^{-\ves_3}\parti_{2}.(x^{\ves_3}\parti_{1}.v_{\zeta})+x^{\ves_3}\parti_{1}.(x^{\ves_2}\parti_{3}.x^{-\ves_3}\parti_{2}.v_{\zeta}))\nonumber\\
&=&\hat{a}_{1}v_{\zeta+\ves_2}\not=0\label{equ:11.6}
\end{eqnarray}
and
\begin{equation}
x^{\ves_2}\parti_{1}.(x^{\pm\ves_2}\parti_{3}.v_{\zeta\mp\ves_2})=x^{\pm\ves_2}\parti_{3}.(x^{\ves_2}\parti_{1}.v_{\zeta\mp\ves_2})=0,
\end{equation}
 we have
 \begin{equation}
 x^{\pm\ves_2}\parti_{3}.v_{\zeta\mp\ves_2}=0.
 \end{equation}
Moreover,
\begin{equation}\label{e:16}
x^{-\ves_3}\parti_{1}.v_{\zeta}=-\frac{1}{\ves_3(\parti_{3})\ves_2(\parti_{2})}[x^{-\ves_2}\parti_{3}, [x^{-\ves_3}\parti_{2}, x^{\ves_2}\parti_{1}]].v_{\zeta}=\hat{a}_{1}v_{\zeta-\ves_3}\not=0
\end{equation}
and
\begin{equation}
x^{-\ves_3}\parti_{1}.(x^{\ves_3}\parti_{l}.v_{\zeta-\ves_3})=x^{\ves_3}\parti_{l}.(x^{-\ves_3}\parti_{1}.v_{\zeta-\ves_3})=0 \ \textrm{
for } l \in \{1, 2\}
\end{equation}
indicate
\begin{equation}
x^{\ves_3}\parti_{l}.v_{\zeta-\ves_3}=0 \ \textrm{ for }l \in \{1, 2\}.
\end{equation}
From
\begin{equation}
[x^{\pm\ves_1}\parti_{3}, [x^{-\ves_2}\parti_{1},
x^{\mp\ves_1}\parti_{2}]].v_{\zeta+\ves_2}=\mp\ves_2(\parti_{2})\ves_1(\parti_{1})x^{-\ves_2}\parti_{3}.v_{\zeta+\ves_2}=0
\end{equation}
it can be derived
\begin{equation}
x^{\pm\ves_1}\parti_{3}.v_{\zeta\mp\ves_1}=0.
\end{equation}

So we have
\begin{equation}\label{eq:9}
x^{\es}\parti.v_{\zeta-\es}=0 \ \textrm{ for all }\es \in \{ \pm\ves_1, \pm\ves_2,\pm\ves_3 \} \textrm{ and } \parti \in
\kn{\es}\cap\{\parti_{1}, \parti_{2}, \parti_{3}\}.
\end{equation}

Thereby  we get
\begin{eqnarray}
& &x^{-\ves_2}\parti_{1}.v_{\zeta}=-\frac{1}{\ves_3(\parti_{3})\ves_2(\parti_{2})}[x^{-\ves_3}\parti_{2}, [x^{-\ves_2}\parti_{3}, x^{\ves_3}\parti_{1}]].v_{\zeta}=\hat{a}_{1}v_{\zeta-\ves_2}, \label{e:17}\\
& &x^{\pm\ves_3}\parti_{2}.v_{\zeta}=\pm\frac{1}{\ves_3(\parti_{3})\ves_1(\parti_{1})}[x^{-\ves_1}\parti_{3}, [x^{\pm\ves_3}\parti_{1}, x^{\ves_1}\parti_{2}]].v_{\zeta}=\hat{a}_{2}v_{\zeta\pm\ves_3},\\
& &x^{-\ves_1}\parti_{2}.v_{\zeta}=-\frac{1}{\ves_3(\parti_{3})\ves_1(\parti_{1})}[x^{-\ves_3}\parti_{1}, [x^{-\ves_1}\parti_{3}, x^{\ves_3}\parti_{2}]].v_{\zeta}=\hat{a}_{2}v_{\zeta-\ves_1},\\
& &x^{\pm\ves_1}\parti_{3}.v_{\zeta}=\pm\frac{1}{\ves_2(\parti_{2})\ves_1(\parti_{1})}[x^{-\ves_2}\parti_{1}, [x^{\pm\ves_1}\parti_{2}, x^{\ves_2}\parti_{3}]].v_{\zeta}=\hat{a}_{3}v_{\zeta\pm\ves_1},\\
& &x^{-\ves_2}\parti_{3}.v_{\zeta}=-\frac{1}{\ves_2(\parti_{2})\ves_1(\parti_{1})}[x^{-\ves_1}\parti_{2}, [x^{-\ves_2}\parti_{1}, x^{\ves_1}\parti_{3}]].v_{\zeta}=\hat{a}_{3}v_{\zeta-\ves_2}.\label{eq:10}
\end{eqnarray}

Define $\eta \in D^{\ast}$ by $\eta(\parti_{1})=\hat{a}_{1}$,
$\eta(\parti_{2})=\hat{a}_{2}$ and $\eta(\parti_{3})=\hat{a}_{3}$. Then we can write (\ref{e:15}), (\ref{equ:11.6}), (\ref{e:16}) and  (\ref{e:17})--(\ref{eq:10}) uniformly as
 \begin{equation}\label{e:18}
x^{\es}\parti.v_{\zeta}=\eta(\parti)v_{\zeta+\es} \textrm{ for all }\es \in \{
\pm\ves_1, \pm\ves_2,\pm\ves_3 \} \textrm{ and }\parti \in \kn{\es}\cap\{\parti_{1},
\parti_{2}, \parti_{3}\}.
\end{equation}
 Relations (\ref{eq:9}) and (\ref{e:18}) show (\ref{eq:15}). So we have completely specified the general action of ${\cal X}$ (c.f. (\ref{e:7}), (\ref{eq:9}), (\ref{e:18})). Since ${\cal X}$ generates ${\cal S}(\G, D)$ (c.f. (\ref{eq:31})), we deduce
\begin{equation}
 M\simeq \mathscr{A}_{\mu,\eta} \simeq \mathscr{A}_{0,\eta}.
\end{equation}

The third subcase is that
\begin{equation}\label{eq:11}
x^{\es}\parti.M_{\zeta}=\{0\} \textrm{ for all }\es \in \{
\pm\ves_1, \pm\ves_2,\pm\ves_3 \} \textrm{ and }\parti \in \kn{\es}\cap\{\parti_{1},
\parti_{2}, \parti_{3}\},
\end{equation}
while there exist some $\es' \in \{ \pm\ves_1, \pm\ves_2,\pm\ves_3 \}$ and $\parti' \in \kn{\es'}\cap\{\parti_{1},
\parti_{2}, \parti_{3}\}$ such that
\begin{equation}\label{eq:12}
x^{\es'}\parti'.v_{\zeta-\es'} \not=0.
\end{equation}
Pick a nonzero vector $v_{\zeta}$ of $M_{\zeta}$. Write
\begin{equation}\label{e:19}
x^{\ves_3}\parti_{1}.v_{\zeta-\ves_3}=\hat{a}_{1}'v_{\zeta},
\
x^{\ves_1}\parti_{2}.v_{\zeta-\ves_1}=\hat{a}_{2}'v_{\zeta},
\
x^{\ves_2}\parti_{3}.v_{\zeta-\ves_2}=\hat{a}_{3}'v_{\zeta}
\end{equation}
with $\hat{a}_{1}',\; \hat{a}_{2}',\; \hat{a}_{3}' \in \mathbb{F}$. Then we have
\begin{eqnarray}
& &x^{\pm\ves_2}\parti_{1}.v_{\zeta\mp\ves_2}=\pm\frac{1}{\ves_3(\parti_{3})\ves_2(\parti_{2})}[x^{-\ves_3}\parti_{2}, [x^{\pm\ves_2}\parti_{3}, x^{\ves_3}\parti_{1}]].v_{\zeta\mp\ves_2}=\hat{a}_{1}'v_{\zeta},\label{eq:13}\\
& &x^{-\ves_3}\parti_{1}.v_{\zeta+\ves_3}=-\frac{1}{\ves_3(\parti_{3})\ves_2(\parti_{2})}[x^{-\ves_2}\parti_{3}, [x^{-\ves_3}\parti_{2}, x^{\ves_2}\parti_{1}]].v_{\zeta+\ves_3}=\hat{a}_{1}'v_{\zeta},\\
& &x^{\pm\ves_3}\parti_{2}.v_{\zeta\mp\ves_3}=\pm\frac{1}{\ves_3(\parti_{3})\ves_1(\parti_{1})}[x^{-\ves_1}\parti_{3}, [x^{\pm\ves_3}\parti_{1}, x^{\ves_1}\parti_{2}]].v_{\zeta\mp\ves_3}=\hat{a}_{2}'v_{\zeta},\\
& &x^{-\ves_1}\parti_{2}.v_{\zeta+\ves_1}=-\frac{1}{\ves_3(\parti_{3})\ves_1(\parti_{1})}[x^{-\ves_3}\parti_{1}, [x^{-\ves_1}\parti_{3}, x^{\ves_3}\parti_{2}]].v_{\zeta+\ves_1}=\hat{a}_{2}'v_{\zeta},\\
& &x^{\pm\ves_1}\parti_{3}.v_{\zeta\mp\ves_1}=\pm\frac{1}{\ves_2(\parti_{2})\ves_1(\parti_{1})}[x^{-\ves_2}\parti_{1}, [x^{\pm\ves_1}\parti_{2}, x^{\ves_2}\parti_{3}]].v_{\zeta\mp\ves_1}=\hat{a}_{3}'v_{\zeta},\\
& &x^{-\ves_2}\parti_{3}.v_{\zeta+\ves_2}=-\frac{1}{\ves_2(\parti_{2})\ves_1(\parti_{1})}[x^{-\ves_1}\parti_{2}, [x^{-\ves_2}\parti_{1}, x^{\ves_1}\parti_{3}]].v_{\zeta+\ves_2}=\hat{a}_{3}'v_{\zeta}.\label{eq:14}
\end{eqnarray}
Define $\eta \in D^{\ast}$ by $\eta(\parti_{1})=\hat{a}_{1}'$,
$\eta(\parti_{2})=\hat{a}_{2}'$ and $\eta(\parti_{3})=\hat{a}_{3}'$. Then we can write (\ref{e:19})--(\ref{eq:14}) uniformly as
 \begin{equation}\label{e:20}
x^{\es}\parti.v_{\zeta-\es}=\eta(\parti)v_{\zeta} \ \textrm{ for all }\es \in \{
\pm\ves_1, \pm\ves_2,\pm\ves_3 \} \textrm{ and }\parti \in \kn{\es}\cap\{\parti_{1},
\parti_{2}, \parti_{3}\}.
\end{equation}
Relations (\ref{eq:11}) and (\ref{e:20}) show (\ref{eq:15}). So we have completely specified the general action of ${\cal X}$ (c.f. (\ref{e:7}), (\ref{eq:11}), (\ref{e:20})).
 Since ${\cal X}$ generates ${\cal S}(\G, D)$ (c.f. (\ref{eq:31})), we deduce
\begin{equation}
M\simeq \mathscr{B}_{\mu,\eta}\simeq \mathscr{B}_{0,\eta}.
\end{equation}

Thus we complete the proof of Lemma 3.3. $\qquad\Box$

\vspace{0.2cm}

\begin{lemma}\label{le:5.5}
Suppose the conditions of Lemma \ref{le:5.4} fail. Moreover, if there exists some $\nu \in \G$ such that
\begin{eqnarray}
& & x^{-\ves_3}\parti_{1}.x^{\ves_3}\parti_{1}.w_{\nu} =0 \textrm{ and }
x^{-\ves_3}\parti_{2}.x^{\ves_3}\parti_{2}.w_{\nu} =0,\nonumber\\
& \textrm{or,} & x^{-\ves_1}\parti_{2}.x^{\ves_1}\parti_{2}.w_{\nu} =0 \textrm{ and }
x^{-\ves_1}\parti_{3}.x^{\ves_1}\parti_{3}.w_{\nu} =0,\nonumber\\
& \textrm{or,} & x^{-\ves_2}\parti_{1}.x^{\ves_2}\parti_{1}.w_{\nu} =0 \textrm { and }
x^{-\ves_2}\parti_{3}.x^{\ves_2}\parti_{3}.w_{\nu} =0,\nonumber
\end{eqnarray}
then $M\simeq \bigoplus
_{\theta \in \G} \mathbb{F}w_{\theta}$, where each component is a
trivial submodule of ${\cal S}(\G, D)$.
\end{lemma}

\noindent{\bf Proof}. Suppose there exists some $\nu \in \G$
such that $x^{-\ves_3}\parti_{1}.x^{\ves_3}\parti_{1}.w_{\nu} =0$ and
$x^{-\ves_3}\parti_{2}.x^{\ves_3}\parti_{2}.w_{\nu} =0$; the other cases
can be proved similarly. By a  translation of the indices if
necessary, we may assume $\nu=0$. In other words,
\begin{equation}\label{eq:16}
x^{-\ves_3}\parti_{1}.x^{\ves_3}\parti_{1}.w_{0} =0 \textrm{ and } x^{-\ves_3}\parti_{2}.x^{\ves_3}\parti_{2}.w_{0} =0.
\end{equation}
Lemma \ref{le:5.2} then implies
\begin{equation}
x^{-\ves_3}\parti_{1}.x^{\ves_3}\parti_{1}.w_{k\ves_3} =0 \textrm{ and } x^{-\ves_3}\parti_{2}.x^{\ves_3}\parti_{2}.w_{k\ves_3} =0 \textrm{ for } k\in\mbb{Z}.
\end{equation}

We begin with the following claim:\vspace{0.3cm}

{\it Claim}. {\it $x^{-\ves_2}\parti_{1}.x^{\ves_2}\parti_{1}.w_{0}=0$. Moreover, either $x^{\ves_2}\parti_{1}.w_{0}=0$ or $x^{-\ves_3}\parti_{2}.x^{\ves_3}\parti_{2}. w_{\ves_2}\not=0$.}\vspace{0.2cm}

If
$x^{\ves_2}\parti_{1}.w_{0}=0$, this claim holds trivially.

Suppose
$x^{\ves_2}\parti_{1}.w_{0}\not=0$. We proceed our proof in three cases.\vspace{0.2cm}

{\it Case 1. $x^{\ves_3}\parti_{1}.w_{0} =0$.}\vspace{0.2cm}

Since $x^{\ves_2}\parti_{1}.w_{0}\not=0$ and
$x^{\ves_3}\parti_{1}.(x^{\ves_2}\parti_{1}.w_{0})
=x^{\ves_2}\parti_{1}.x^{\ves_3}\parti_{1}.w_{0} =0$, we have
$x^{\ves_3}\parti_{1}.w_{\ves_2} =0$, which further implies
\begin{eqnarray}
x^{\ves_2+\ves_3}\parti_{1}.w_{0} &= & \frac{1}{\ves_3(\parti_{3})}[x^{\ves_2}\parti_{3},x^{\ves_3}\parti_{1}].w_{0}\nonumber\\
&=&\frac{1}{\ves_3(\parti_{3})}(x^{\ves_2}\parti_{3}.x^{\ves_3}\parti_{1}.w_{0}-x^{\ves_3}\parti_{1}.(x^{\ves_2}\parti_{3}.w_{0})) =0.
\end{eqnarray}
Thus
\begin{equation}
x^{\ves_2+\ves_3}\parti_{1}.x^{-\ves_3}\parti_{2}.w_{0}=[x^{\ves_2+\ves_3}\parti_{1},x^{-\ves_3}\parti_{2}].w_{0} =-\ves_2(\parti_{2})x^{\ves_2}\parti_{1}.w_{0}\not=0,
 \end{equation}
which indicates
\begin{equation}\label{equ:4.148}
x^{-\ves_3}\parti_{2}.w_{0}\not=0 \textrm{ and } x^{\ves_2+\ves_3}\parti_{1}.w_{-\ves_3}\not=0.
\end{equation}
Since (\ref{eq:16}) shows
\begin{equation}
x^{\ves_3}\parti_{2}.(x^{-\ves_3}\parti_{2}.w_{0})
=x^{-\ves_3}\parti_{2}.x^{\ves_3}\parti_{2}.w_{0} =0,
\end{equation}
we get $x^{\ves_3}\parti_{2}.w_{-\ves_3} =0$ by (\ref{equ:4.148}). Then
\begin{equation}
x^{\ves_3}\parti_{2}.x^{\ves_2}\parti_{1}.w_{-\ves_3}=[x^{\ves_3}\parti_{2},x^{\ves_2}\parti_{1}].w_{-\ves_3} =\ves_2(\parti_{2})x^{\ves_2+\ves_3}\parti_{1}.w_{-\ves_3}\not=0,
\end{equation}
 which implies $x^{\ves_2}\parti_{1}.w_{-\ves_3}\not=0$.

 If $x^{\ves_3}\parti_{1}.w_{-\ves_3}=0$,
 then
 \begin{equation}x^{\ves_3}\parti_{1}.(x^{\ves_2}\parti_{1}.w_{-\ves_3})=x^{\ves_2}\parti_{1}.x^{\ves_3}\parti_{1}.w_{-\ves_3}=0,
 \end{equation}
 which leads to $x^{\ves_3}\parti_{1}.w_{\ves_2-\ves_3}=0$. Thus
 \begin{eqnarray}
 x^{\ves_2+\ves_3}\parti_{1}.w_{-\ves_3}&=&\frac{1}{\ves_3(\parti_{3})}[x^{\ves_2}\parti_{3},x^{\ves_3}\parti_{1}].w_{-\ves_3}\nonumber\\
 &=&\frac{1}{\ves_3(\parti_{3})}(x^{\ves_2}\parti_{3}.x^{\ves_3}\parti_{1}.w_{-\ves_3} - x^{\ves_3}\parti_{1}.(x^{\ves_2}\parti_{3}.w_{-\ves_3}))=0,
  \end{eqnarray}
  which contradicts (\ref{equ:4.148}).

Assume $x^{\ves_3}\parti_{1}.w_{-\ves_3}\not=0$. Then (\ref{equ:4.148}) and
  $x^{\ves_3}\parti_{1}.w_{0} =0$ give
\begin{equation}
0=x^{-\ves_3}\parti_{2}.x^{\ves_3}\parti_{1}.w_{0}=x^{\ves_3}\parti_{1}.(x^{-\ves_3}\parti_{2}.w_{0})\not=0,
\end{equation}
which is a contradiction. \vspace{0.2cm}

{\it Case 2. $x^{-\ves_3}\parti_{1}.w_{0} =0$.}\vspace{0.2cm}

Replacing
$\ves_3$ by $-\ves_3$ in the arguments of Case 1, we similarly get
a contradiction.\vspace{0.2cm}

{\it Case 3. $x^{\ves_3}\parti_{1}.w_{0}\not =0$ and
$x^{-\ves_3}\parti_{1}.w_{0} \not=0$.}\vspace{0.2cm}

Equation (\ref{eq:16}) shows
\begin{equation}
x^{-\ves_3}\parti_{1}.(x^{\ves_3}\parti_{1}.w_{0})=x^{\ves_3}\parti_{1}.(x^{-\ves_3}\parti_{1}.w_{0})=0,
\end{equation}
which indicates
$x^{-\ves_3}\parti_{1}.w_{\ves_3} =0$ and $x^{\ves_3}\parti_{1}.w_{-\ves_3} =0$.

If $x^{\ves_2}\parti_{1}.w_{-\ves_3}\not =0$, in analogy with Case 1, $x^{\ves_3}\parti_{1}.w_{-\ves_3} =0$ leads to a contradiction.

If $x^{\ves_2}\parti_{1}.w_{\ves_3}\not =0$, in analogy with Case 2, $x^{-\ves_3}\parti_{1}.w_{\ves_3} =0$ leads to a contradiction.

Assume
$x^{\ves_2}\parti_{1}.w_{\pm\ves_3} =0$. Then we have
\begin{equation}\label{eq:17}
x^{\ves_3}\parti_{1}.(x^{\ves_2}\parti_{1}.w_{0})
=x^{\ves_2}\parti_{1}.(x^{\ves_3}\parti_{1}.w_{0}) =0.
\end{equation}
Since $x^{\ves_2}\parti_{1}.w_{0} \not=0$, (\ref{eq:17}) shows $x^{\ves_3}\parti_{1}.w_{\ves_2}
=0$, which further shows
\begin{equation}\label{eq:18}
x^{\ves_3}\parti_{1}.(x^{-\ves_2}\parti_{1}.w_{\ves_2})=x^{-\ves_2}\parti_{1}.x^{\ves_3}\parti_{1}.w_{\ves_2}=0.
\end{equation}
Since $x^{\ves_3}\parti_{1}.w_{0}\not =0$, we get
$x^{-\ves_2}\parti_{1}.w_{\ves_2}=0$ by (\ref{eq:18}). So we have
\begin{equation}
x^{-\ves_2}\parti_{1}.x^{\ves_2}\parti_{1}.w_{0}=0,
\end{equation}
which completes the proof of the first statement of the claim.
Moreover,
\begin{equation} 0\not=\ves_2^{2}(\parti_2)x^{\ves_2}\parti_{1}.w_{0}=[
x^{\ves_3}\parti_{2},[x^{-\ves_3}\parti_{2}, x^{\ves_2}\parti_{1}]].w_{0}=
x^{\ves_3}\parti_{2}.x^{-\ves_3}\parti_{2}. x^{\ves_2}\parti_{1}.w_{0}
\end{equation}
gives
\begin{equation}\label{equ:4.154}
x^{-\ves_3}\parti_{2}.x^{\ves_3}\parti_{2}. w_{\ves_2}\not=0.
\end{equation}

Thus this claim follows.\vspace{0.3cm}

By symmetry, we can similarly prove
$x^{-\ves_1}\parti_{2}.x^{\ves_1}\parti_{2}.w_{0} =0$.

So we conclude that $x^{-\ves_3}\parti_{1}.x^{\ves_3}\parti_{1}.w_{0} =0$ and $x^{-\ves_3}\parti_{2}.x^{\ves_3}\parti_{2}.w_{0} =0$ lead to
\begin{equation}\label{eq:19}
x^{-\ves_2}\parti_{1}.x^{\ves_2}\parti_{1}.w_{0}=0 \textrm{ and } x^{-\ves_1}\parti_{2}.x^{\ves_1}\parti_{2}.w_{0} =0.
\end{equation}

On the other hand, as the conditions of Lemma \ref{le:5.4} fail, we must have
\begin{equation}
x^{-\ves_1}\parti_{3}.x^{\ves_1}\parti_{3}.w_{0} =0 \textrm{ or }
x^{-\ves_2}\parti_{3}.x^{\ves_2}\parti_{3}.w_{0} =0.
\end{equation}
Without loss of generality, we assume
$x^{-\ves_2}\parti_{3}.x^{\ves_2}\parti_{3}.w_{0} =0$. Observe that the above claim gives
$x^{-\ves_2}\parti_{1}.x^{\ves_2}\parti_{1}.w_{0} =0$. Thus, in analogy with (\ref{eq:19}), from $x^{-\ves_2}\parti_{3}.x^{\ves_2}\parti_{3}.w_{0}=0$ and $x^{-\ves_2}\parti_{1}.x^{\ves_2}\parti_{1}.w_{0} =0$ we can deduce
\begin{equation}
x^{-\ves_1}\parti_{3}.x^{\ves_1}\parti_{3}.w_{0}
=0.
\end{equation}

Therefore we obtain
\begin{equation}
x^{-\ves}\parti.x^{\ves}\parti. w_{0}=0 \textrm{ for all }\ves \in \{ \pm\ves_1, \pm\ves_2, \pm\ves_3 \} \textrm{ and } \parti \in
\kn{\ves}\cap\{\parti_{1}, \parti_{2}, \parti_{3}\}.
\end{equation}
Lemma \ref{le:5.2} and the above discussion then further imply
\begin{equation}\label{equ:3.153}
x^{-\ves}\parti.x^{\ves}\parti. w_{\theta}=0 \ \textrm{ for all }\ves \in \{\pm\ves_1, \pm\ves_2, \pm\ves_3 \} , \; \parti \in
\kn{\ves}\cap\{\parti_{1}, \parti_{2}, \parti_{3}\} \textrm{ and } \theta\in \G.
\end{equation}

Since the above claim states either $x^{\ves_2}\parti_{1}.w_{0}=0$ or $x^{-\ves_3}\parti_{2}.x^{\ves_3}\parti_{2}. w_{\ves_2}\not=0$, (\ref{equ:3.153}) confirms
\begin{equation}
x^{\ves_2}\parti_{1}.w_{0}=0.
\end{equation}
Likewise, it can be proved
\begin{equation}
x^{\ves}\parti. w_{0}=0\;\; \textrm{ for all }\ves \in \{\pm\ves_1,
\pm\ves_2, \pm\ves_3 \}\;\;\mbox{and}\;\; \parti \in
\kn{\ves}\cap\{\parti_{1},
\parti_{2}, \parti_{3}\}.
\end{equation}
Moreover, it can be similarly proved
\begin{equation}\label{eq:20}
x^{\ves}\parti. w_{\theta}=0\;\; \textrm{ for all }\ves \in \{\pm\ves_1,
\pm\ves_2, \pm\ves_3 \}\;\;\parti \in
\kn{\ves}\cap\{\parti_{1},
\parti_{2}, \parti_{3}\}\;\;\mbox{and}\;\;  \theta\in \G.
\end{equation}

Therefore, (\ref{eq:20}) shows $M\simeq \bigoplus
_{\theta \in \G} \mathbb{F}w_{\theta}$, where each component is a
trivial submodule of ${\cal S}(\G, D)$.$\qquad\Box$\vspace{0.2cm}

\begin{lemma}\label{le:5.6}
 The conditions of Lemma \ref{le:5.4} and that of Lemma \ref{le:5.5} enumerate all the possibilities.
 \end{lemma}

\noindent{\bf Proof}. Suppose otherwise, namely, both the conditions of Lemma \ref{le:5.4} and that of Lemma \ref{le:5.5} fail. Then this will lead to a contradiction. The following proof is divided into two cases.\vspace{0.2cm}

{\it Case 1. $x^{-\ves_3}\parti_{1}.x^{\ves_3}\parti_{1}.w_{0} =0$.}

Since both the conditions of Lemma \ref{le:5.4} and that of Lemma \ref{le:5.5} fail, we get
\begin{eqnarray}
& & x^{-\ves_3}\parti_{1}.x^{\ves_3}\parti_{1}.w_{0} =0 , \
x^{-\ves_3}\parti_{2}.x^{\ves_3}\parti_{2}.w_{0} \not=0,\
 x^{-\ves_1}\parti_{2}.x^{\ves_1}\parti_{2}.w_{0} =0 , \\
&  &x^{-\ves_1}\parti_{3}.x^{\ves_1}\parti_{3}.w_{0} \not=0,\
x^{-\ves_2}\parti_{3}.x^{\ves_2}\parti_{3}.w_{0} =0 ,\
x^{-\ves_2}\parti_{1}.x^{\ves_2}\parti_{1}.w_{0} \not=0.
\end{eqnarray}
Lemma \ref{le:5.2} then implies, for any $\theta\in \G$,
\begin{eqnarray}
& & x^{-\ves_3}\parti_{1}.x^{\ves_3}\parti_{1}.w_{\theta} =0 , \
x^{-\ves_3}\parti_{2}.x^{\ves_3}\parti_{2}.w_{\theta} \not=0,\
x^{-\ves_1}\parti_{2}.x^{\ves_1}\parti_{2}.w_{\theta} =0 , \\
&  & x^{-\ves_1}\parti_{3}.x^{\ves_1}\parti_{3}.w_{\theta} \not=0,\
x^{-\ves_2}\parti_{3}.x^{\ves_2}\parti_{3}.w_{\theta} =0,\
x^{-\ves_2}\parti_{1}.x^{\ves_2}\parti_{1}.w_{\theta} \not=0 .
\end{eqnarray}

Since $x^{-\ves_2}\parti_{1}.x^{\ves_2}\parti_{1}.w_{\theta} \not=0$ for any $\theta\in \G$, we get $x^{\pm\ves_2}\parti_{1}.w_{\theta} \not=0$ for any $\theta\in \G$. Similarly, $x^{\pm\ves_3}\parti_{2}.w_{\theta} \not=0$ for any $\theta\in \G$.

If $x^{\ves_3}\parti_{1}.w_{0} =0$, by
\begin{equation}\label{e:21}
x^{\ves_3}\parti_{1}.(x^{\pm\ves_2}\parti_{1}.w_{j\ves_2})=x^{\pm\ves_2}\parti_{1}.x^{\ves_3}\parti_{1}.w_{j\ves_2}
 \end{equation}
 and  induction on $j$, we get $x^{\ves_3}\parti_{1}.w_{j\ves_2} =0$ for all $j \in \mathbb{Z}$. Then by the fact
 \begin{equation}\label{e:22}
 x^{\ves_3}\parti_{1}.(x^{\pm\ves_3}\parti_{2}.w_{j\ves_2+k\ves_3})=x^{\pm\ves_3}\parti_{2}.x^{\ves_3}\parti_{1}.w_{j\ves_2+k\ves_3}
  \end{equation}
  and  induction on $k$, we obtain $x^{\ves_3}\parti_{1}.w_{j\ves_2+k\ves_3} =0$ for all $j,k \in \mathbb{Z}$. Therefore,
\begin{equation}
0\not=\ves_3(\parti_{3})\ves_2(\parti_{2})x^{\ves_2}\parti_{1}.w_{0} =[[x^{\ves_3}\parti_{1}, x^{\ves_2}\parti_{3}], x^{-\ves_3}\parti_{2}].w_{0}=0
\end{equation}
leads to a contradiction.

If $x^{\ves_3}\parti_{1}.w_{0} \not=0$, then we get
$x^{-\ves_3}\parti_{1}.w_{\ves_3} =0$. By similar arguments as those from (\ref{e:21}) to (\ref{e:22}), we can prove
\begin{equation}
x^{-\ves_3}\parti_{1}.w_{j\ves_2+k\ves_3} =0 \textrm{ for all } j,k \in \mathbb{Z}.
\end{equation}
So it can be deduced
\begin{equation}
0\not=-\ves_3(\parti_{3})\ves_2(\parti_{2})x^{\ves_2}\parti_{1}.w_{0} =[[x^{-\ves_3}\parti_{1}, x^{\ves_2}\parti_{3}], x^{\ves_3}\parti_{2}].w_{0}=0,
\end{equation}
which is a contradiction.
\vspace{0.2cm}

{\it Case 2. $x^{-\ves_3}\parti_{1}.x^{\ves_3}\parti_{1}.w_{0} \not=0$.}\vspace{0.2cm}

Since both the conditions of Lemma \ref{le:5.4} and that of Lemma \ref{le:5.5} fail, we get
\begin{eqnarray}
& & x^{-\ves_3}\parti_{1}.x^{\ves_3}\parti_{1}.w_{0} \not=0 , \ x^{-\ves_2}\parti_{1}.x^{\ves_2}\parti_{1}.w_{0} =0,\
x^{-\ves_2}\parti_{3}.x^{\ves_2}\parti_{3}.w_{0} \not=0 , \\
& &x^{-\ves_1}\parti_{3}.x^{\ves_1}\parti_{3}.w_{0} =0,\
x^{-\ves_1}\parti_{2}.x^{\ves_1}\parti_{2}.w_{0} \not=0,\ x^{-\ves_3}\parti_{2}.x^{\ves_3}\parti_{2}.w_{0} =0.
\end{eqnarray}
Lemma \ref{le:5.2} then implies, for any $\theta\in \G$,
\begin{eqnarray}
& & x^{-\ves_3}\parti_{1}.x^{\ves_3}\parti_{1}.w_{\theta} \not=0 , \ x^{-\ves_2}\parti_{1}.x^{\ves_2}\parti_{1}.w_{\theta} =0,\
x^{-\ves_2}\parti_{3}.x^{\ves_2}\parti_{3}.w_{\theta} \not=0 , \\
& &x^{-\ves_1}\parti_{3}.x^{\ves_1}\parti_{3}.w_{\theta} =0,\
x^{-\ves_1}\parti_{2}.x^{\ves_1}\parti_{2}.w_{\theta} \not=0,\ x^{-\ves_3}\parti_{2}.x^{\ves_3}\parti_{2}.w_{0} =0.
\end{eqnarray}

Since
\begin{equation}
[[x^{\pm\ves_2}\parti_{1}, x^{\ves_3}\parti_{2}], x^{\mp\ves_2}\parti_{3}].w_0=\pm\ves_3(\parti_{3})\ves_2(\parti_{2})x^{\ves_3}\parti_{1}.w_0,
\end{equation}
similarly, in analogy with Case 1, we will get a contradiction. We omit the details.
$\qquad\Box$ \vspace{0.2cm}

In summary, Lemmas \ref{le:5.4}--\ref{le:5.6}, together with Lemma \ref{le:3.1} and Theorem \ref{th:3.2}, give: \vspace{0.2cm}

\begin{lemma}\label{le:5.7}
If $\mbox{\it dim}\; D=3$ and $\G \simeq \mathbb{Z}^{3}$, then
Theorem \ref{th:4.1} holds.
\end{lemma}

\section{General case of $\dim D=3$}

In this section, we deal with the general case of $\dim D=3$, that
is, we want to prove:

\begin{lemma}\label{le:4.3}
If $\mbox{\it dim}\; D =3$, then Theorem \ref{th:4.1} holds.
\end{lemma}

\noindent{\bf Proof}. Assume that $M=\bigoplus_{\theta\in
\G}M_{\theta}$ is a $\G$-graded ${\cal S}(\G, D)$-module with $\dim
M_{\theta}=1$ for each $\theta\in \G$. Suppose $\G'$ is an arbitrary
subgroup of
 $\G$ such that $\bigcap_{\al \in \G'}\kn \al = \{0\}$. Then ${\cal S}(\G', D)$ is not only a subalgebra of ${\cal S}(\G, D)$, but also a simple generalized
 divergence-free Lie algebra. Define
\begin{equation}
M(\nu,\G')=\bigoplus_{\theta\in \G'} M_{\nu+\theta}\quad \textrm{
for any } \nu \in \G.
\end{equation}
Then $M(\nu,\G')$ is a $\G'$-graded ${\cal S}(\G',
D)$-submodule.

By Zorn's Lemma and Lemma \ref{le:5.7}, there exists a maximal
subgroup $\G_{0}$ of $\G$ such that $\bigcap_{\al \in \G_{0}}\kn \al
= \{0\}$ and the following condition holds:\vspace{0.2cm}

{\bf (C1)} Let $\{\theta_{r} \mid r \in I\}$ be the set of all
representatives of cosets of $\G_{0}$ in $\G$. For any $r\in I$,
$M(\theta_{r},\G_{0})$ is isomorphic to:
\begin{eqnarray}
&(\rmnum{1})&  \mathscr{M}_{\mu_{r}}({\cal S}(\G_{0}, D)) \textrm{ for some }  \mu_{r}\in D^{\ast}\backslash \G_{0}; \\
&(\rmnum{2})&  \mathscr{M}_{\mu_{r}}({\cal S}(\G_{0}, D)) \textrm{ for some }  \mu_{r}\in \G_{0}; \\
&(\rmnum{3})&  \mathscr{A}_{\mu_{r},\eta_{r}}({\cal S}(\G_{0}, D)) \textrm{ for some }  \mu_{r}\in \G_{0} \textrm{ and } \eta_{r}\in D^{\ast}\backslash\{0\}; \\
&(\rmnum{4})&  \mathscr{B}_{\mu_{r},\eta_{r}}({\cal S}(\G_{0}, D)) \textrm{ for some }  \mu_{r}\in \G_{0} \textrm{ and } \eta_{r}\in D^{\ast}\backslash\{0\}; \\
&(\rmnum{5})&  \bigoplus_{\nu\in \G_{0}}\mathbb{F}w_{\nu}, \textrm{
where each component is a trivial submodule of }
{\cal S}( \G_{0},D ).
   \end{eqnarray}

To prove this lemma, it suffices to show $\G_{0}=\G$. Suppose
$\G_{0}\not=\G$. We will see that this leads to a contradiction.
Picking $\theta_{1} \in \G\backslash \G_{0}$, we set
\begin{equation}\label{e:23}
\G_{1}=\G_{0}+\mathbb{Z}\theta_{1}.
\end{equation}
Let
\begin{equation}
\{\rho_{\lmd} \mid \lmd\in J\} \textrm{ be the
set of all representatives of cosets of } \G_{1} \textrm{ in } \G.
\end{equation}
Our intention is to prove that $\G_1$ satisfies condition (C1), which results in the contradiction. Let
\begin{equation}\label{eq:36}
\{k\theta_{1} \mid k\in K\} \textrm{ be the set of all
representatives of cosets of } \G_{0} \textrm{ in }\G_{1},\textrm{ where } K\subseteq \mathbb{Z}.
\end{equation}
Then (\ref{e:23}) shows $\#(K)>1$.

Fix $\lmd\in J$ hereafter. We want to see how ${\cal S}(\G_1, D)$ act on $M(\rho_{\lmd},\G_1)$ and therefore deduce $\G_1$ satisfies condition (C1).

Condition (C1) implies:\vspace{0.2cm}

{\bf (C2)} For any $k\in \mathbb{Z}$, $M(\rho_{\lmd}+k \theta_{1},\G_{0})$
is isomorphic to
\begin{eqnarray}
&(\rmnum{1})& \mathscr{M}_{\mu_{k}'}({\cal S}(\G_{0}, D)) \textrm{ for some }  \mu_{k}'\in D^{\ast}\backslash \G_{0}; \\
&(\rmnum{2})& \mathscr{M}_{\mu_{k}'}({\cal S}(\G_{0}, D)) \textrm{ for some }  \mu_{k}'\in \G_{0}; \\
&(\rmnum{3})& \mathscr{A}_{\mu_{k}',\eta_{k}'}({\cal S}(\G_{0}, D)) \textrm{ for some }  \mu_{k}'\in \G_{0} \textrm{ and }  \eta_{k}'\in D^{\ast}\backslash\{0\}; \\
&(\rmnum{4})&  \mathscr{B}_{\mu_{k}',\eta_{k}'}({\cal S}(\G_{0}, D)) \textrm{ for some }  \mu_{k}'\in \G_{0} \textrm{ and }  \eta_{k}'\in D^{\ast}\backslash\{0\}; \\
&(\rmnum{5})& \bigoplus_{\nu\in \G_{0}}\mathbb{F}w_{\nu}, \textrm{
where each component is a trivial submodule of }
{\cal S}( \G_{0},D ).
\end{eqnarray}
In the first four cases of (C2), $\mu_{k}'$'s are respectively
chosen so that there exist nonzero $v_{\rho_{\lmd}+k\theta_{1}+\be}
\in M_{\rho_{\lmd}+k\theta_{1}+\be}$ with $\be\in \G_{0}$ such that
\begin{equation}\label{equ:4.162}
x^{\al}\parti.v_{\rho_{\lmd}+k\theta_{1}+\be}=(\be+\mu_{k}')(\parti)v_{\rho_{\lmd}+k\theta_{1}+\be+\al}
\end{equation}
for $\al\in \G_{0}\backslash\{0\}$, $\parti \in \kn{\al}$ and
$\be\in \G_{0}$ with $\be+\mu_{k}'
\not=0\not=\be+\mu_{k}'+\al$.

Choose $\es_{1},\es_{2} \in \G_{0}\backslash\{0\}$ such that
$\kn{\theta_{1}}\cap \kn{\es_{1}}\cap \kn{\es_{2}}=\{0\}$. Set
\begin{equation}\label{eq:44}
G_{0}=\mathbb{Z}\theta_{1}+\mathbb{Z}\es_{1}+\mathbb{Z}\es_{2}.
\end{equation}
Then ${\cal S}( G_{0},D )$ is not only a subalgebra of ${\cal S}( \G,D )$, but also a simple generalized divergence-free
Lie algebra. By Lemma \ref{le:5.7}, we know that:\vspace{0.2cm}

{\bf (C3)} For any $\be \in \G_{0}$, $M(\rho_{\lmd}+\be,G_{0})$ is isomorphic to
\begin{eqnarray}
   &(\rmnum{1})& \mathscr{M}_{\mu_{\be}''}({\cal S}( G_{0},D )) \textrm{ for some } \mu_{\be}''\in D^{\ast}\backslash G_{0}; \\
   &(\rmnum{2})& \mathscr{M}_{\mu_{\be}''}({\cal S}( G_{0},D )) \textrm{ for some }  \mu_{\be}''\in G_{0}; \\
    &(\rmnum{3})& \mathscr{A}_{\mu_{\be}'',\eta_{\be}''}({\cal S}( G_{0},D )) \textrm{ for some }  \mu_{\be}''\in G_{0} \textrm{ and } \eta_{\be}''\in D^{\ast}\backslash\{0\}; \\
    &(\rmnum{4})& \mathscr{B}_{\mu_{\be}'',\eta_{\be}''}({\cal S}( G_{0},D )) \textrm{ for some }  \mu_{\be}''\in G_{0} \textrm{ and } \eta_{\be}''\in D^{\ast}\backslash\{0\}; \\
    &(\rmnum{5})& \bigoplus_{\nu\in G_{0}}\mathbb{F}w_{\nu},  \textrm{ where each component is a trivial submodule of } {\cal S}( G_{0},D ).
   \end{eqnarray}
In the first four cases of (C3), $\mu_{\be}''$'s are respectively chosen so that, there exist nonzero
$v_{\rho_{\lmd}+\be+\nu} \in
M_{\rho_{\lmd}+\be+\nu}$ with $\nu\in G_{0}$ such that
\begin{equation}\label{eq:33}
x^{\al}\parti.v_{\rho_{\lmd}+\be+\nu}=(\nu+\mu_{\be}'')(\parti)v_{\rho_{\lmd}+\be+\nu+\al}
\end{equation}
for $\al\in G_{0}\backslash\{0\}$, $\parti \in \kn{\al}$ and
$\nu\in G_{0}$ with $\nu+\mu_{\be}''
\not=0\not=\nu+\mu_{\be}''+\al$.

Our aim now is to show that $\G_1$ satisfies condition (C1). The
idea is to sew the ${\cal S}( \G_{0},D )$-submodules $M(\rho_{\lmd}+k\theta_{1},\G_{0})$
for $k\in K$ together via the action of ${\cal S}( G_{0},D )$. Before getting to that, we give some
observations.

For any $\be'\in \G_{1}\backslash \G_{0}$, we can choose
$\theta_{1}'\in (\theta_{1}+\G_{0})\cap G_{0}$ and
 $\es_{1}',\es_{2}' \in \G_{0}\backslash\{0\}$ such that
 \begin{equation}\label{eq:47}
 \kn \theta_{1}' \cap \kn \es_{1}' \cap \kn \es_{2}'=\{0\} \textrm{ and } \be'\in G_{1}=\mathbb{Z}\theta_{1}'+\mathbb{Z}\es_{1}'+\mathbb{Z}\es_{2}'.
 \end{equation}
Explicitly, we write
 $\be'=k\theta_{1}+\be$ with $k\in K\backslash\{0\}$ and $\be\in \G_{0}$. If $\be=0$, we take $\theta_{1}'=\theta_{1}$, and choose $\es_{1}', \es_{2}' \in \G_{0}\backslash\{0\}$ such
    that $\kn{\theta_{1}'}\cap \kn{\es_{1}'}\cap \kn{\es_{2}'}=\{0\}$. In the case $\be\not=0$ and $\kn{\theta_{1}}\not=\kn{\be}$,
 we take $\theta_{1}'=\theta_{1}$, $\es_{1}'=\be$, and choose $\es_{2}' \in \G_{0}\backslash\{0\}$ such
   that $\kn{\theta_{1}'}\cap \kn{\es_{1}'}\cap \kn{\es_{2}'}=\{0\}$.
     When $\kn{\theta_{1}}=\kn{\be}$, by the fact
     $\kn(\theta_{1}+\es_{1})\not=\kn(\be-k\es_{1})$ (c.f. (\ref{eq:44})),  we can take $\theta_{1}'=\theta_{1}+\es_{1}$,
    $\es_{1}'=\be-k\es_{1}$, and choose $\es_{2}' \in \G_{0}\backslash\{0\}$ such that
    $\kn{\theta_{1}'}\cap \kn{\es_{1}'}\cap \kn{\es_{2}'}=\{0\}$.

Notice that, for the above defined $G_1$, ${\cal S}( G_1,D )$ is not only a subalgebra of ${\cal S}( \G,D )$, but also a simple generalized divergence-free
Lie algebra.

Now we proceed our analysis through sewing the ${\cal S}( \G_{0},D )$-submodules
$M(\rho_{\lmd}+k\theta_{1},\G_{0})$ for $k\in K$ together. The proof
is divided into five cases.

\vspace{0.3cm}

\noindent{\it{Case 1.}} \emph{There exists some $k_{0} \in K$
such that $M(\rho_{\lmd}+k_{0}\theta_{1},\G_{0})\simeq \bigoplus_{\nu\in
\G_{0}}\mathbb{F}w_{\nu}$, where each component is a trivial
submodule of ${\cal S}( \G_{0},D )$.} \vspace{0.2cm}

Recall that $G_{0}=\mathbb{Z}\theta_{1}+\mathbb{Z}\es_{1}+\mathbb{Z}\es_{2}$, where $\es_{1},\es_{2} \in \G_{0}\backslash\{0\}$. Since $M(\rho_{\lmd}+k_{0}\theta_{1},\G_{0})\simeq \bigoplus_{\nu\in
\G_{0}}\mathbb{F}w_{\nu}$, where each component is a trivial
submodule of ${\cal S}( \G_{0},D )$, we have
\begin{equation}\label{equ:8.8}
x^{\es_{s}}\parti.M_{\rho_{\lmd}+k_{0}\theta_{1}+\be+l\es_{1}+m\es_{2}}=\{0\}\quad \forall \be \in \G_{0}, \;  l, m\in \mathbb{Z},\; \parti \in \kn{\es_{s}} \textrm{ and }s\in\{1,2\}.
\end{equation}

Pick $\parti' \in \kn\es_{1}\backslash\kn\es_{2}$. In (C3), if $M(\rho_{\lmd}+\be,G_{0})$ for some $\be \in \G_{0}$ is isomorphic to one of the first four cases, then there exist $l',m'\in\mbb{Z}$ such that
\begin{equation}
x^{\es_{1}}\parti'.M_{\rho_{\lmd}+\be+k_{0}\theta_{1}+l'\es_{1}+m'\es_{2}}\not=\{0\},
\end{equation}
 which contradicts (\ref{equ:8.8}). So for any $\be \in \G_{0}$,
\begin{equation}\label{equ:4.1164}
M(\rho_{\lmd}+\be,G_{0})\simeq \bigoplus_{\nu\in G_{0}}\mathbb{F}w_{\nu},
\end{equation}
where each component is a trivial submodule of ${\cal S}( G_{0},D )$. Then we have
\begin{equation}\label{equ:8.9}
x^{\es_{s}}\parti.M_{\rho_{\lmd}+\be+k
\theta_{1}+l\es_{1}+m\es_{2}}=\{0\},\quad \forall
\be \in \G_{0},\; k,l,m\in \mathbb{Z},\; \parti \in \kn{\es_{s}} \textrm{ and } s\in\{1,2\}.
\end{equation}

Pick $\parti' \in \kn\es_{1}\backslash\kn\es_{2}$. Observe that in (C2),
for any fixed $k \in K$, if $M(\rho_{\lmd}+k \theta_{1},\G_{0})$ is isomorphic to one of the first four cases, then
there exist $\be\in \G_{0}$ and $l',m'\in\mbb{Z}$ such that
\begin{equation}
x^{\es_{1}}\parti'.M_{\rho_{\lmd}+\be+k\theta_{1}+l'\es_{1}+m'\es_{2}}\not=\{0\},
\end{equation}
which contradicts (\ref{equ:8.9}). So for any $k\in
K$,
\begin{equation}\label{eq:21}
M(\rho_{\lmd}+k\theta_{1},\G_{0})\simeq \bigoplus_{\nu\in
\G_{0}}\mathbb{F}w_{\nu},
\end{equation}
where each component is a
trivial submodule of ${\cal S}( \G_{0},D )$. In other words, from (\ref{eq:36}) we see that
\begin{equation}\label{equ:4.1165}
x^{\al}\parti.M(\rho_{\lmd},\G_{1})=\{0\},\quad \forall \al\in
\G_{0}\backslash \{0\},\, \parti \in \kn{\al} .
\end{equation}

Recall that for any $\be'\in \G_{1}\backslash \G_{0}$, we can choose
$\theta_{1}'\in (\theta_{1}+\G_{0})\cap G_{0}$ and $\es_{1}',\es_{2}' \in
\G_{0}\backslash\{0\}$ such that $\kn \theta_{1}' \cap \kn \es_{1}'
\cap \kn \es_{2}'=\{0\}$ and $\be'\in
G_{1}=\mathbb{Z}\theta_{1}'+\mathbb{Z}\es_{1}'+\mathbb{Z}\es_{2}'$ (c.f. (\ref{eq:47})). Moreover, we note that $\{x^{\pm\al}\parti\mid \al\in\{\theta_{1}',\es_{1}',\es_{2}'\},\, \parti \in \kn{\al}\}$ generates the simple generalized divergence-free Lie algebra ${\cal S}( G_{1},D )$. As (\ref{equ:4.1164}) and
(\ref{equ:4.1165}) show
\begin{equation}\label{equ:4.165}
x^{\pm\al}\parti.M(\rho_{\lmd},\G_{1})=\{0\}\quad \textrm{ for }
\al\in\{\theta_{1}',\es_{1}',\es_{2}'\} \textrm{ and } \parti \in \kn{\al} ,
\end{equation}
we can deduce
\begin{equation}
x^{\tau}\parti.M(\rho_{\lmd},\G_{1})=\{0\}\quad \textrm{ for }\tau\in G_1\backslash\{0\} \textrm{ and }
\parti \in \kn{\tau} .
\end{equation}
In particular,
\begin{equation}\label{equ:4.11165}
x^{\be'}\parti.M(\rho_{\lmd},\G_{1})=\{0\}\quad \textrm{ for }
\parti \in \kn{\be'} .
\end{equation}

 Since $\be'\in \G_{1}\backslash \G_{0}$ is arbitrary,
 (\ref{equ:4.1165}) and (\ref{equ:4.11165}) yield
\begin{equation}
M(\rho_{\lmd},\G_{1})\simeq \bigoplus_{\nu\in \G_{1}}\mathbb{F}w_{\nu}, \textrm{
where each component is a trivial submodule of } {\cal S}( \G_{1},D ) .
\end{equation}

\vspace{0.2cm}

\noindent{\it{Case 2.}} \emph{There exists some $k_{0} \in K$ such
that
\begin{equation}\label{equ:4.168}
M(\rho_{\lmd}+k_{0}\theta_{1},\G_{0})\simeq \mathscr{A}_{\mu'_{k_{0}},\eta'_{k_{0}}}({\cal S}( \G_{0},D)) \textrm{ for some }\mu_{k_{0}}'\in \G_{0} \textrm{ and } \eta_{k_0}'\in D^{\ast}\backslash\{0\}.
\end{equation}}

To begin with, we want to confirm which cases the submodules in (C3) are isomorphic to. Write $M_{\theta}=\mbb{F}w'_{\theta}$ for $\theta\in\G$.

Since $M(\rho_{\lmd}+k_{0}\theta_{1},\G_{0})\simeq \mathscr{A}_{\mu'_{k_{0}},\eta'_{k_{0}}}({\cal S}( \G_{0},D))$, (\ref{equ:4.162}) shows
\begin{eqnarray}\label{equ:4.169}
 & & x^{-\es_{s}}\parti.x^{\es_{s}}\parti.w'_{\rho_{\lmd}+k_{0}\theta_{1}+\be+l\es_{1}+m\es_{2}}\nonumber\\
 &=&(\be+l\es_{1}+m\es_{2}+\mu_{k_{0}}')^{2}(\parti)w'_{\rho_{\lmd}+k_{0}\theta_{1}+\be+l\es_{1}+m\es_{2}}
\end{eqnarray}
for $\be \in \G_{0}$ and $l, m\in \mathbb{Z}$ such that $\be+l\es_{1}+m\es_{2}+\mu_{k_{0}}'\not=0\not=\be+l\es_{1}+m\es_{2}+\mu_{k_{0}}'+\es_s$, where $\parti \in \kn{\es_{s}}$ and $s\in\{1,2\}$.
Pick $\parti' \in
\kn\es_{1}\backslash\kn\es_{2}$. Fix any $\be_{0} \in \G_{0}$. From (\ref{equ:4.169}) we see that there
exist $l',m'\in \mathbb{Z}$ such that
\begin{equation}
x^{-\es_{1}}\parti'.x^{\es_{1}}\parti'.w'_{\rho_{\lmd}+k_{0}\theta_{1}+\be_{0}+l'\es_{1}+m'\es_{2}}\not=0,
\end{equation}
which then implies, in (C3),
\begin{equation}\label{eq:22}
 M(\rho_{\lmd}+\be_{0},G_{0}) \textrm{ can only be isomorphic to one of the first four cases}.
 \end{equation}
So by (\ref{eq:33}) we have
\begin{eqnarray}\label{equ:4.170}
 & & x^{-\es_{s}}\parti.x^{\es_{s}}\parti.w'_{\rho_{\lmd}+\be_{0}+k_{0}\theta_{1}+l\es_{1}+m\es_{2}}\nonumber\\
 &=&(k_{0}\theta_{1}+l\es_{1}+m\es_{2}+\mu_{\be_{0}}'')^{2}(\parti)w'_{\rho_{\lmd}+\be_{0}+k_{0}\theta_{1}+l'\es_{1}+m'\es_{2}}
\end{eqnarray}
for $l,m\in \mathbb{Z}$ such that
$k_{0}\theta_{1}+l\es_{1}+m\es_{2}+\mu_{\be_{0}}''\not=0 \not= k_{0}\theta_{1}+l\es_{1}+m\es_{2}+\mu_{\be_{0}}''+\es_s$, where $\parti \in \kn{\es_{s}}$ and $s\in\{1,2\}$. Comparing (\ref{equ:4.170}) with (\ref{equ:4.169}), we see that
\begin{equation}\label{eq:34}
(k_{0}\theta_{1}+l_s\es_{1}+m_s\es_{2}+\mu_{\be_{0}}'')^{2}(\parti)=(\be_{0}+l_s\es_{1}+m_s\es_{2}+\mu_{k_{0}}')^{2}(\parti)
\end{equation}
for $l_s,m_s\in \mathbb{Z}$ such that
$\be_0+l_s\es_{1}+m_s\es_{2}+\mu_{k_{0}}'\not=0\not
=\be_0+l_s\es_{1}+m_s\es_{2}+\mu_{k_{0}}'+\es_s$ and
$k_{0}\theta_{1}+l_s\es_{1}+m_s\es_{2}+\mu_{\be_{0}}''\not=0 \not=
k_{0}\theta_{1}+l_s\es_{1}+m_s\es_{2}+\mu_{\be_{0}}''+\es_s$, where
$\parti \in \kn{\es_{s}}$ and $s\in\{1,2\}$. Moreover, since there
exist infinite pairs of such $l_s$ and $m_s$ for either
$s\in\{1,2\}$,  we can deduce from (\ref{eq:34}) that
\begin{equation}\label{equ:6.3}
(k_{0}\theta_{1}+\mu_{\be_{0}}'')(\parti)=(\be_{0}+\mu_{k_{0}}')(\parti)\quad\textrm
{ for }\parti \in (\kn \es_{1}\cup \kn \es_{2})\backslash(\kn
\es_{1}\cap \kn \es_{2}),
\end{equation}
\begin{equation}\label{equ:8.10}
(k_{0}\theta_{1}+\mu_{\be_{0}}'')^{2}(\parti'')=(\be_{0}+\mu_{k_{0}}')^{2}(\parti'')\quad\textrm
{ for }\parti'' \in \kn \es_{1}\cap \kn \es_{2}.
\end{equation}
Suppose that
$(k_{0}\theta_{1}+\mu_{\be_{0}}'')(\parti'')\not=(\be_{0}+\mu_{k_{0}}')(\parti'')$
for some $\parti'' \in \kn \es_{1}\cap \kn \es_{2}$. Then
(\ref{equ:8.10}) gives rise to
\begin{equation}\label{eq:24}
(k_{0}\theta_{1}+\mu_{\be_{0}}'')(\parti'')=-(\be_{0}+\mu_{k_{0}}')(\parti'')\not=0.
\end{equation}
We will see this leads to a contradiction. Recall that $\parti'\in \kn \es_{1}\backslash\kn \es_{2}$. Choose $l_1,m_1\in \mathbb{Z}$ such that $\be_0+l_1\es_{1}+m_1\es_{2}+\mu_{k_{0}}'\not=0\not=\be_0+(l_1+1)\es_{1}+m_1\es_{2}+\mu_{k_{0}}'$,
$k_{0}\theta_{1}+l_1\es_{1}+m_1\es_{2}+\mu_{\be_{0}}''\not=0 \not= k_{0}\theta_{1}+(l_1+1)\es_{1}+m_1\es_{2}+\mu_{\be_{0}}''$ and
\begin{equation}\label{eq:23}
(\be_{0}+l_1\es_{1}+m_1\es_{2}+\mu_{k_{0}}')(\parti')\not=0.
\end{equation}
Observe that (\ref{equ:4.168}) and (\ref{eq:22}) enable us to choose nonzero vectors
$u_{\rho_{\lmd}+k_{0}\theta_{1}+\be_{0}+(l_1+1)\es_{1}+m_1\es_{2}}$
and
$u'_{\rho_{\lmd}+k_{0}\theta_{1}+\be_{0}+(l_1+1)\es_{1}+m_1\es_{2}}$
in $M_{\rho_{\lmd}+k_{0}\theta_{1}+\be_{0}+(l_1+1)\es_{1}+m_1\es_{2}}$, respectively,
such that,
\begin{equation}\label{equ:6.4}
  x^{\es_{1}}\parti.w'_{\rho_{\lmd}+k_{0}\theta_{1}+\be_{0}+l_1\es_{1}+m_1\es_{2}}
 =(\be_{0}+l_1\es_{1}+m_1\es_{2}+\mu_{k_{0}}')(\parti)u_{\rho_{\lmd}+k_{0}\theta_{1}+\be_{0}+(l_1+1)\es_{1}+m_1\es_{2}},
\end{equation}
\begin{equation}\label{equ:6.5}
 x^{\es_{1}}\parti.w'_{\rho_{\lmd}+k_{0}\theta_{1}+\be_{0}+l_1\es_{1}+m_1\es_{2}}
 =(k_{0}\theta_{1}+l_1\es_{1}+m_1\es_{2}+\mu_{\be_{0}}'')(\parti)u'_{\rho_{\lmd}+k_{0}\theta_{1}
 +\be_{0}+(l_1+1)\es_{1}+m_1\es_{2}}
\end{equation}
for $\parti \in \kn \es_{1}$. Since (\ref{equ:6.3}) and (\ref{eq:23}) indicate
\begin{equation}
(k_{0}\theta_{1}+l_1\es_{1}+m_1\es_{2}+\mu_{\be_{0}}'')(\parti')=(\be_{0}+l_1\es_{1}+m_1\es_{2}+\mu_{k_{0}}')(\parti')\not=0,
\end{equation}
taking $\parti=\parti'$ in (\ref{equ:6.4}) and (\ref{equ:6.5}), we get
\begin{equation}\label{eq:25}
u_{\rho_{\lmd}+k_{0}\theta_{1}+\be_{0}+(l_1+1)\es_{1}+m_1\es_{2}}=u'_{\rho_{\lmd}+k_{0}
\theta_{1}+\be_{0}+(l_1+1)\es_{1}+m_1\es_{2}}.
\end{equation}
On the other hand, (\ref{eq:24}) shows
\begin{equation}
(k_{0}\theta_{1}+l_1\es_{1}+m_1\es_{2}+\mu_{\be_{0}}'')(\parti'')=
-(\be_{0}+l_1\es_{1}+m_1\es_{2}+\mu_{k_{0}}')(\parti'')\not=0.
\end{equation}
Taking $\parti=\parti''$ in (\ref{equ:6.4}) and (\ref{equ:6.5}), we
get
\begin{equation}
u_{\rho_{\lmd}+k_{0}\theta_{1}+\be_{0}+(l_1+1)\es_{1}+m_1\es_{2}}=-u'_{\rho_{\lmd}+k_{0}\theta_{1}+\be_{0}+(l_1+1)\es_{1}+m_1\es_{2}},
\end{equation}
which contradicts (\ref{eq:25}). So we must have
\begin{equation}\label{eq:26}
(k_{0}\theta_{1}+\mu_{\be_{0}}'')(\parti'')=(\be_{0}+\mu_{k_{0}}')(\parti'')\quad
\textrm{ for }\parti'' \in \kn \es_{1}\cap \kn \es_{2}.
\end{equation}
As $\kn \es_{1}+
\kn \es_{2}=D$, (\ref{equ:6.3}) and (\ref{eq:26}) lead to
\begin{equation}\label{equ:4.171}
\mu_{\be_{0}}''-\be_{0}=\mu_{k_{0}}'-k_{0}\theta_{1}\quad \textrm { for any }\be_{0} \in \G_{0}.
\end{equation}

Now we can confirm which cases the submodules in (C3) are isomorphic
to. Since $-\mu_{k_{0}}'\in \G_{0}$,  we deduce from
(\ref{equ:4.171}) that
\begin{equation}\label{e:24}
\mu_{-\mu_{k_{0}}'}''=\mu_{k_{0}}'-k_{0}\theta_{1}-\mu_{k_{0}}'=-k_{0}\theta_{1}\in
G_{0}.
\end{equation}
Moreover, the fact
$M(\rho_{\lmd}+k_{0}\theta_{1},\G_{0})\simeq
\mathscr{A}_{\mu'_{k_{0}},\eta'_{k_{0}}}({\cal S}( \G_{0},D))$ gives
\begin{equation}
x^{\es_{s}}\parti.M_{\rho_{\lmd}+k_{0}\theta_{1}-\mu_{k_{0}}'-\es_{s}}=\{0\}\quad
\textrm{ for all } \parti \in \kn{\es_{s}}\textrm{ and } s\in\{1,2\},
\end{equation}
\begin{equation}\label{e:25}
x^{\es_{s}}\parti.M_{\rho_{\lmd}+k_{0}\theta_{1}-\mu_{k_{0}}'}\not=\{0\}\quad
\textrm{ for some } s\in\{1,2\}\textrm{ and } \parti \in \kn{\es_{s}}.
\end{equation}
So (C3), (\ref{eq:22}), (\ref{e:24})--(\ref{e:25}) and $\kn \es_{1}+
\kn \es_{2}=D$ show
\begin{equation}\label{equ:8.11}
M(\rho_{\lmd}-\mu_{k_{0}}',G_{0})\simeq \mathscr{A}_{\mu_{-\mu_{k_{0}}'}'',\eta}({\cal S}( G_{0},D))=\mathscr{A}_{-k_{0}\theta_{1},\eta}({\cal S}( G_{0},D)) \ \textrm{ for some } \eta \in D^\ast\backslash\{0\}.
\end{equation}
For any $\be\in \G_{0}$, if $\mu_{\be}''\in G_{0}$, we have
\begin{equation}
\be=-\mu_{k_{0}}'+\mu_{\be}''+k_{0}\theta_{1}\in -\mu_{k_{0}}'+G_{0}
\end{equation}
by (\ref{equ:4.171}). Conversely,
$\be\in\G_0\backslash(-\mu_{k_{0}}'+G_{0})$ implies
$\mu_{\be}''\not\in G_{0}$. Thus this together with (C3),
(\ref{eq:22}) and (\ref{equ:8.11}) imply
\begin{equation}\label{equ:4.180}
M(\rho_{\lmd}+\be,G_{0})=M(\rho_{\lmd}-\mu_{k_{0}}',G_{0})\simeq
\mathscr{A}_{-k_{0}\theta_{1},\eta}({\cal S}( G_{0},D))\simeq
\mathscr{A}_{\be+\mu_{k_{0}}'
-k_{0}\theta_{1},\eta}({\cal S}( G_{0},D))
\end{equation}
for $ \be\in  \G_{0}\cap(-\mu_{k_{0}}'+G_{0})$, and
\begin{equation}\label{equ:4.181}
M(\rho_{\lmd}+\be,G_{0})\simeq \mathscr{M}_{\mu_{\be}''}({\cal S}( G_{0},D))=\mathscr{M}_{\be+\mu_{k_{0}}'
-k_{0}\theta_{1}}({\cal S}( G_{0},D))\; \textrm{ for }\be\in
\G_{0}\backslash(-\mu_{k_{0}}'+G_{0}).
\end{equation}

Then we want to confirm which cases the submodules in (C2) are
isomorphic to.

Fix $\be \in \G_{0}$ and $k\in \mathbb{Z}$. Note that (\ref{equ:4.180}) and (\ref{equ:4.181}) give
\begin{eqnarray}\label{equ:4.172}
 & & x^{-\es_{s}}\parti.x^{\es_{s}}\parti.w'_{\rho_{\lmd}+k\theta_{1}+\be+l\es_{1}+m\es_{2}}\nonumber\\
 &=&(k\theta_{1}+l\es_{1}+m\es_{2}+\mu_{\be}'')^{2}(\parti)w'_{\rho_{\lmd}+k\theta_{1}+\be+l\es_{1}+m\es_{2}}
\end{eqnarray}
for $l, m \in
\mathbb{Z}$ such that $k\theta_{1}+l\es_{1}+m\es_{2}+\mu_{\be}''\not=0\not=k\theta_{1}+l\es_{1}+m\es_{2}+\mu_{\be}''+\es_s$, where $\parti \in \kn{\es_{s}}$ and $s\in\{1,2\}$. Pick $\parti' \in
\kn\es_{1}\backslash\kn\es_{2}$. Then (\ref{equ:4.172}) shows that there exist $l',m'\in \mathbb{Z}$
such that
\begin{equation}
x^{-\es_{1}}\parti'.x^{\es_{1}}\parti'.w'_{\rho_{\lmd}+k\theta_{1}+\be+l'\es_{1}+m'\es_{2}}\not=0,
\end{equation}
which further implies, in (C2),
\begin{equation}\label{eq:27}
M(\rho_{\lmd}+k\theta_{1},\G_{0}) \textrm{ can only be
isomorphic to one of the first four cases}.
\end{equation}
So (\ref{equ:4.162}) yields
\begin{eqnarray}\label{equ:4.173}
 & & x^{-\es_{s}}\parti.x^{\es_{s}}\parti.w'_{\rho_{\lmd}+k\theta_{1}+\be+l\es_{1}+m\es_{2}}\nonumber\\
 &=&(\be + l \es_{1} + m \es_{2} + \mu_{k}')^{2}(\parti)w'_{\rho_{\lmd}+k\theta_{1}+\be+l\es_{1}+m\es_{2}}
\end{eqnarray}
for $l,m \in\mbb{Z}$
such that $\be + l\es_{1} + m\es_{2} + \mu_{k}'\not=0\not=\be + l\es_{1} + m\es_{2} + \mu_{k}'+\es_s$, where $\parti \in
\kn{\es_{s}}$ and $s\in\{1,2\}$. By the similar arguments as those from
(\ref{equ:4.169}) to (\ref{equ:4.171}), we get
\begin{equation}\label{equ:4.174}
\mu_{k}'-k \theta_{1}=\mu_{\be}''-\be\quad\textrm { for any }k \in
\mbb{Z}.
\end{equation}
So combining (\ref{equ:4.171})
with (\ref{equ:4.174}), we obtain
\begin{equation}\label{equ:4.175}
\mu_{k}'-k\theta_{1}=\mu_{\be}''-\be=\mu_{k_{0}}'-k_{0}\theta_{1}\quad\textrm {
for any }k \in \mathbb{Z},\; \be \in \G_{0}.
\end{equation}

Since $\mu_{k_{0}}'\in \G_{0}$, from (\ref{eq:36}) and
(\ref{equ:4.175}), we deduce
\begin{equation}
\mu_{k}'=(k-k_{0})
\theta_{1}+\mu_{k_{0}}'\not\in \G_{0} \quad \textrm{ for } k \in K\backslash\{k_{0}\}.
\end{equation}
Thus (\ref{eq:27}) and (C2) imply
\begin{equation}\label{equ:4.176}
M(\rho_{\lmd}+k\theta_{1},\G_{0})\simeq
\mathscr{M}_{\mu_{k}'}({\cal S}(
\G_{0},D))=\mathscr{M}_{\mu_{k_{0}}'+(k-k_{0})
\theta_{1}}({\cal S}(
\G_{0},D))\  \textrm{ for }k \in K\backslash\{k_{0}\}.
\end{equation}

Next we sew the ${\cal S}(\G_0,D)$-submodules $M(\rho_{\lmd}+k\theta_{1},\G_{0})$ for
$k\in K$ together, and derive how ${\cal S}(\G_1,D)$ act on
$M(\rho_{\lmd},\G_1)$. We fix some $\be_{0}\in
\G_{0}\backslash(-\mu_{k_{0}}'+G_{0})$ in this part.

For any $k\in K\backslash\{k_0\}$, (\ref{equ:4.176}) enables us to choose
$\{0\not=v_{\rho_{\lmd}+k\theta_{1}+\be}'\in M_{\rho_{\lmd}+k\theta_{1}+\be}\mid \be\in \G_{0}\}$ such that
\begin{equation}\label{equ:4.182}
x^{\al}\parti.v_{\rho_{\lmd}+k\theta_{1}+\be}'=(\be+\mu_{k_{0}}'+(k-k_{0})
\theta_{1}) (\parti)v_{\rho_{\lmd}+k\theta_{1}+\be+\al}'
\end{equation}
for $\al\in \G_{0}\backslash\{0\}$ and $\parti \in \kn{\al}$.
Moreover, by (\ref{equ:4.168}), we can choose
$\{0\not=v_{\rho_{\lmd}+k_0\theta_{1}+\be}'\in
M_{\rho_{\lmd}+k_0\theta_{1}+\be}\mid \be\in \G_{0}\}$ such that
\begin{equation}\label{eq:37}
x^{\al}\parti.v_{\rho_{\lmd}+k_0\theta_{1}+\be}'=(\be+\mu_{k_{0}}') (\parti)v_{\rho_{\lmd}+k_0\theta_{1}+\be+\al}'
\end{equation}
for $\parti \in \kn{\al}$, $\al\in \G_{0}\backslash\{0\}$ and $\be\in\G_0\backslash\{-\mu_{k_{0}}'\}$, and
\begin{equation}\label{eq:35}
x^{\al}\parti.v_{\rho_{\lmd}+k_0\theta_{1}-\mu_{k_{0}}'}'=\eta_{k_{0}}' (\parti)v_{\rho_{\lmd}+k_0\theta_{1}-\mu_{k_{0}}'+\al}'\textrm{ for } \parti \in \kn{\al} \textrm{ and }\al\in \G_{0}\backslash\{0\}.
\end{equation}
On the other hand, since
$\be_{0}\in  \G_{0}\backslash(-\mu_{k_{0}}'+G_{0})$,
(\ref{equ:4.181}) enables us to choose
$\{0\not=v_{\rho_{\lmd}+\be_{0}+\nu}''\in M_{\rho_{\lmd}+\be_{0}+\nu}\mid\nu\in G_{0}\}$ such that
\begin{equation}\label{equ:4.183}
x^{\al}\parti.v_{\rho_{\lmd}+\be_{0}+\nu}''=(\nu+\mu_{k_{0}}'+\be_{0}-k_{0}\theta_{1})
(\parti)v_{\rho_{\lmd}+\be_{0}+\nu+\al}''
\end{equation}
for $\al\in G_{0}\backslash\{0\}$ and $\parti \in \kn{\al}$. Assume
\begin{equation}\label{equ:4.184}
v_{\rho_{\lmd}+\be_{0}+k\theta_{1}+\ga}''=\hat{b}_{k,\ga}v_{\rho_{\lmd}+\be_{0}+k\theta_{1}+\ga}'\quad\textrm{
for } k\in K \textrm{
and } \ga\in \G_{0}\cap G_{0},
\end{equation}
where $\hat{b}_{k,\ga}\in\mbb{F}$.

We now show that $\hat{b}_{k,\ga}=\hat{b}_{k}$ is independent of $\ga$. Since $\be_{0}\in  \G_{0}\backslash(-\mu_{k_{0}}'+G_{0})$, we have
\begin{equation}
\be_{0}+\ga\not=-\mu_{k_{0}}' \quad \textrm{ for any } \ga\in \G_{0}\cap G_{0}.
\end{equation}
So fixing any $k\in K$, from (\ref{equ:4.182}) or (\ref{eq:37}), and from
(\ref{equ:4.183}) we respectively derive
\begin{equation}\label{equ:4.185}
x^{\al}\parti.v_{\rho_{\lmd}+k\theta_{1}+\be_{0}+\ga}'=(\be_{0}+\ga+\mu_{k}')(\parti)v_{\rho_{\lmd}+k\theta_{1}+\be_{0}+\ga+\al}',
\end{equation}
\begin{equation}\label{equ:4.186}
x^{\al}\parti.v_{\rho_{\lmd}+\be_{0}+k\theta_{1}+\ga}''=
(\be_{0}+\ga+\mu_{k}')(\parti)v_{\rho_{\lmd}+\be_{0}+k\theta_{1}+\ga+\al}''
\end{equation}
for $\al\in (\G_{0}\cap G_{0})\backslash\{0\}$, $\parti \in \kn{\al}$ and $\ga\in \G_{0}\cap G_{0}$,
 where $\mu_{k}'=\mu_{k_{0}}'+(k-k_{0})\theta_{1}$. From $\be_{0}\in  \G_{0}\backslash(-\mu_{k_{0}}'+G_{0})$ again, it follows
\begin{equation}\label{eq:32}
\be_{0}+\mu_{k}'+\ga=\be_{0}+\mu_{k_{0}}'+(k-k_{0})\theta_{1}+\ga\not=0\quad \textrm{for any }\ga\in \G_{0}\cap G_{0}.
\end{equation}
We see that $\hat{b}_{k,\ga}=\hat{b}_{k,0}$ in two cases. On one hand, for $\al'\in (\G_{0}\cap G_{0})\backslash\{0\}$ such that $\kn \al' \not=\kn (\be_{0}+\mu'_{k})$, taking $\al=\al'$, $\ga=0$ and $\parti \in \kn\al' \backslash \kn (\be_{0}+\mu'_{k})$ in (\ref{equ:4.185}) and (\ref{equ:4.186}), we obtain
\begin{equation}\label{eq:28}
\hat{b}_{k,0}=\hat{b}_{k,\al'}.
 \end{equation}
Noticing that such $\al'$ does exist (e.g. $\es_1$ or $\es_2$ does the trick), we denote one by $\al'_0$. On the other hand, for $\al''\in (\G_{0}\cap G_{0})\backslash\{0\}$ such that $\kn \al'' =\kn (\be_{0}+\mu'_{k})$, we first have
\begin{equation}
\kn (\al''-\al'_{0})\not =\kn (\be_{0}+\mu'_{k}),
\end{equation}
which implies $\hat{b}_{k,0}=\hat{b}_{k,\al''-\al'_0}$ with $\al'$ replaced by $\al''-\al'_0$ in (\ref{eq:28}).
 Moreover,  it follows from (\ref{eq:32}) that
 \begin{equation}
 \kn \al'' =\kn (\be_{0}+\mu'_{k})=\kn (\be_{0}+\mu'_{k}+\al''),
 \end{equation}
which implies $\kn \al'_{0}\not =\kn (\be_{0}+\mu'_{k}+\al''-\al'_{0})$. Thus, putting $\al=\al'_0$, $\ga=\al''-\al'_{0}$ and taking $\parti \in \kn\al'_{0} \backslash \kn (\be_{0}+\mu'_{k}+\al''-\al'_{0})$ in (\ref{equ:4.185}) and (\ref{equ:4.186}), we obtain
\begin{equation}\label{eq:29}
\hat{b}_{k,\al''}=\hat{b}_{k,\al''-\al'_{0}}=\hat{b}_{k,0}.
\end{equation}
So combining (\ref{eq:28}) with (\ref{eq:29}), we get
\begin{equation}\label{equ:4.190}
\hat{b}_{k,\ga}=\hat{b}_{k,0}=\hat{b}_k \quad \textrm{for } \ga\in \G_{0}\cap G_{0}.
\end{equation}

Recall that $\{k\theta_{1} \mid k\in K\}$ is the set of all
representatives of cosets of $\G_{0}$ in $\G_{1}$ (c.f (\ref{eq:36})). Multiplying
$\{v_{\rho_{\lmd}+k\theta_{1}+\be}'\mid \be \in \G_{0}\}$ by $\hat{b}_{k}$ for each $k\in K$, we get a new set
\begin{equation}
\{0\not=v_{\rho_{\lmd}+\varrho}\in M_{\rho_{\lmd}+\varrho  }\mid \varrho\in
\G_{1}\}.
\end{equation}
We then consider the action of ${\cal S}(\G_1, D)$ on $M(\rho_{\lmd},\G_{1})$.

First, (\ref{equ:4.182}), (\ref{eq:37}) and
(\ref{eq:35}) show
\begin{equation}\label{equ:4.191}
x^{\al}\parti.v_{\rho_{\lmd}+\varrho}=(\varrho+\mu'_{k_{0}} -k_{0}\theta_{1})(\parti)v_{\rho_{\lmd}+\varrho+\al}
\end{equation}
for $\al\in \G_{0}\backslash\{0\}$, $\parti \in \kn{\al}$ and
$\varrho\in \G_{1}\backslash\{k_0\theta_{1}-\mu_{k_{0}}'\}$, and
\begin{equation}\label{eq:38}
x^{\al}\parti.v_{\rho_{\lmd}+k_0\theta_{1}-\mu_{k_{0}}'}=\eta_{k_{0}}' (\parti)v_{\rho_{\lmd}+k_0\theta_{1}-\mu_{k_{0}}'+\al}\ \textrm{ for }\al\in \G_{0}\backslash\{0\}, \; \parti \in \kn{\al},
\end{equation}
while (\ref{equ:4.183}) and (\ref{equ:4.184}) imply
\begin{equation}\label{equ:4.192}
x^{\tau
}\parti'.v_{\rho_{\lmd}+\be_{0}+\nu}=(\nu+\be_{0}+\mu'_{k_{0}}
-k_{0}\theta_{1}) (\parti')v_{\rho_{\lmd}+\be_{0}+\nu+\tau }
\end{equation}
for $\tau\in G_{0}\backslash\{0\}$, $\parti'\in \kn{\tau}$ and $\nu\in G_{0}$.

Second, we determine the action of ${\cal S}( G_{0}, D)$ on $M(\rho_{\lmd},\G_1)$. Since ${\cal X}_0=\{x^{\pm \al}\parti \mid \al\in\{\theta_{1}, \es_{1}, \es_{2}\},\,\parti\in \kn{\al}\}$ generates ${\cal S}( G_{0}, D)$, we only need to determine the action of ${\cal X}_0$ on $M(\rho_{\lmd},\G_1)$. Recall that $\es_1,\es_2\in\G_0$. Then (\ref{equ:4.191}) and (\ref{eq:38}) show
\begin{equation}\label{equ:4.193}
x^{\pm \es_{s}}\parti.v_{\rho_{\lmd}+\varrho}=(\varrho+\mu'_{k_{0}} -k_{0}\theta_{1})(\parti)v_{\rho_{\lmd}+\varrho\pm \es_{s}}
\end{equation}
for $\varrho\in \G_{1}\backslash\{k_0\theta_{1}-\mu_{k_{0}}'\}$, $\parti \in \kn{\es_{s}}$ and $s\in \{1,2\}$, and
\begin{equation}\label{eq:39}
x^{\pm \es_{s}}\parti.v_{\rho_{\lmd}+k_0\theta_{1}-\mu_{k_{0}}'}=\eta_{k_{0}}' (\parti)v_{\rho_{\lmd}+k_0\theta_{1}-\mu_{k_{0}}'\pm \es_{s}}\ \textrm{ for }\parti \in \kn{\es_{s}}\textrm{ and }s\in \{1,2\}.
\end{equation}
Fix any $\varrho \in \G_1\backslash\{k_0\theta_{1}-\mu_{k_{0}}'\}$. As $\G_1=\G_0+\mbb{Z}\theta_1=\G_0+G_0$, we see that
\begin{equation}\label{eq:42}
\varrho=\be+\nu \textrm{ for some } \be\in\G_0 \textrm{ and }\nu\in G_0.
\end{equation}
We then derive how $x^{\pm \theta_{1}}\parti$ act on $v_{\rho_{\lmd}+\varrho}$ in two cases. If $\kn{\theta_{1}}=
\kn(\varrho+\mu'_{k_{0}} -k_{0}\theta_{1})$, (\ref{equ:4.180}) or
(\ref{equ:4.181}) gives
\begin{equation}\label{equ:4.194}
x^{\pm \theta_{1}}\parti.v_{\rho_{\lmd}+\varrho}=0=(\varrho+\mu'_{k_{0}} -k_{0}\theta_{1})(\parti)v_{\rho_{\lmd}+\varrho \pm \theta_{1}} \textrm{ for }\parti \in \kn{\theta_{1}}.
\end{equation}
Otherwise, we have $\textrm{dim}(\kn{\theta_{1}}\cap \kn(\varrho+\mu'_{k_{0}} -k_{0}\theta_{1}))=1$. Recall that $\varrho=\be+\nu$ (c.f. (\ref{eq:42})). Choose $l,\, m \in \mathbb{Z}$ such that
\begin{equation}\label{equ:4.195}
\kn{\theta_{1}}\cap \kn(\varrho+\mu'_{k_{0}} -k_{0}\theta_{1})\cap \kn(\be-\be_0+l\es_{1}+m\es_{2})=\{0\}.
\end{equation}
Pick nonzero $\parti_{1}\in \kn{\theta_{1}}\cap \kn(\varrho+\mu'_{k_{0}} -k_{0}\theta_{1})$ and $\parti_{2}\in \kn{\theta_{1}}\cap \kn(\be-\be_{0}+l\es_{1}+m\es_{2})$, then $\{\parti_{1}, \parti_{2}\}$ forms a basis of $\kn{\theta_{1}}$ and $(\varrho+\mu'_{k_{0}} -k_{0}\theta_{1})(\parti_{2})\not=0$. On one hand, (\ref{equ:4.180}) or (\ref{equ:4.181}) shows
\begin{equation}\label{equ:4.196}
x^{\pm \theta_{1}}\parti_{1}.v_{\rho_{\lmd}+\varrho}=0=(\varrho+\mu'_{k_{0}} -k_{0}\theta_{1})(\parti_{1})v_{\rho_{\lmd}+\varrho \pm \theta_{1}}.
\end{equation}
On the other hand, since (\ref{equ:4.192}) gives the action of $x^{\pm \theta_{1}}\parti_2$ on $M(\rho_{\lmd}+\be_{0},G_0)$, we extend it to the action of $x^{\pm \theta_{1}}\parti_2$ on $v_{\rho_{\lmd}+\varrho}$ with the assistance of (\ref{equ:4.191}). Notice that $\be-\be_0+l\es_{1}+m\es_{2}\in\G_0\backslash\{0\}$ and $M_{\rho_{\lmd}+\be_{0}-l\es_{1}-m\es_{2}+\nu}\subseteq M(\rho_{\lmd}+\be_{0},G_0)$ (c.f. (\ref{eq:42})). The following diagram illustrates the idea:
$${\xymatrix{    M_{\rho_{\lmd}+\varrho} \ar[rr]^-{{x^{\pm \theta_{1}}\parti_2}} & &  M_{\rho_{\lmd}+\varrho \pm \theta_1}  \\
   & & &\\
 M_{\rho_{\lmd}+\be_{0}-l\es_{1}-m\es_{2}+\nu} \ar[uu]^-{{x^{\be-\be_{0}+l\es_{1}+m\es_{2}}\parti_2}} \ar[rr]^-{{x^{\pm \theta_{1}}\parti_2}} & &  M_{\rho_{\lmd}+\be_{0}-l\es_{1}-m\es_2+\nu \pm \theta_1} \ar[uu]_-{{x^{\be-\be_{0}+l\es_{1}+m\es_{2}}\parti_2}} }}$$
 Since $\parti_{2}\in \kn(\be-\be_{0}+l\es_{1}+m\es_{2})$, we have
\begin{equation}\label{eq:43}
(\be_{0}-l\es_{1}-m\es_{2}+\nu+\mu'_{k_{0}} -k_{0}\theta_{1})(\parti_2)=(\varrho+\mu'_{k_{0}} -k_{0}\theta_{1})(\parti_2)\not=0.
 \end{equation}
 Applying (\ref{equ:4.191}) and (\ref{equ:4.192}) to
 \begin{equation}\label{eq:41}
 x^{\pm \theta_{1}}\parti_{2}.x^{\be-\be_{0}+l\es_{1}+m\es_{2}}\parti_{2}.v_{\rho_{\lmd}+\be_{0}-l\es_{1}-m\es_{2}+\nu}=x^{\be-\be_{0}+l\es_{1}+m\es_{2}}\parti_{2}.x^{\pm \theta_{1}}\parti_{2}.v_{\rho_{\lmd}+\be_{0}-l\es_{1}-m\es_{2}+\nu}
 \end{equation}
and making use of (\ref{eq:43}), we get
\begin{equation}\label{equ:4.197}
x^{\pm \theta_{1}}\parti_{2}.v_{\rho_{\lmd}+\varrho}=(\varrho+\mu'_{k_{0}} -k_{0}\theta_{1})(\parti_{2})v_{\rho_{\lmd}+\varrho \pm \theta_{1}}.
\end{equation}
So (\ref{equ:4.194}), (\ref{equ:4.196}) and (\ref{equ:4.197}) show
\begin{equation}\label{eq:40}
x^{\pm \theta_{1}}\parti.v_{\rho_{\lmd}+\varrho}=(\varrho+\mu'_{k_{0}} -k_{0}\theta_{1})(\parti)v_{\rho_{\lmd}+\varrho \pm \theta_{1}} \ \textrm{ for }\varrho\in \G_1\backslash\{k_0\theta_{1}-\mu_{k_{0}}'\},\; \parti \in \kn{\theta_{1}}.
\end{equation}
Afterwards we need to derive the action of $x^{\pm \theta_{1}}\parti$ on $v_{\rho_{\lmd}+k_0\theta_{1}-\mu_{k_{0}}'}$. Recall that $\kn{\theta_{1}}\cap \kn{\es_{1}}\cap \kn{\es_{2}}=\{0\}$ (c.f. (\ref{eq:44})). Pick nonzero $\parti'_{1}\in \kn{\theta_{1}}\cap \kn{\es_1} $ and $\parti'_{2}\in \kn{\theta_{1}}\cap \kn{\es_2}$, then $\{\parti'_{1}, \parti'_{2}\}$ forms a basis of $\kn{\theta_{1}}$. Note that $\es_2(\parti'_1)\not=0$ and $\es_1(\parti'_2)\not=0$. Moreover, we pick $\tilde{\parti_1}\in \kn{\es_1}\backslash \kn{\theta_{1}} $ and $\tilde{\parti_2}\in \kn{\es_2}\backslash \kn{\theta_{1}} $. Then using (\ref{equ:4.193}), (\ref{eq:39}) and (\ref{eq:40}), we get
\begin{eqnarray}
x^{\pm \theta_{1}}\parti'_1.v_{\rho_{\lmd}+k_0\theta_{1}-\mu_{k_{0}}'} & = & \pm\frac{1}{\theta_1(\tilde{\parti_1})\es_1(\parti'_2)}[x^{-\es_1}\tilde{\parti_1},[x^{\pm\theta_{1}}\parti'_2,x^{\es_1}\parti'_1]].v_{\rho_{\lmd}+k_0\theta_{1}-\mu_{k_{0}}'}\nonumber\\
&=& \eta_{k_{0}}' (\parti'_1)v_{\rho_{\lmd}+k_0\theta_{1}-\mu_{k_{0}}'\pm \theta_{1}}.\label{eq:45}
\end{eqnarray}
Similarly, we can prove
\begin{equation}\label{eq:46}
x^{\pm \theta_{1}}\parti'_2.v_{\rho_{\lmd}+k_0\theta_{1}-\mu_{k_{0}}'} =\eta_{k_{0}}' (\parti'_2)v_{\rho_{\lmd}+k_0\theta_{1}-\mu_{k_{0}}'\pm \theta_{1}}.
\end{equation}
Since ${\cal X}_0=\{x^{\pm \al}\parti \mid \al\in\{\theta_{1},
\es_{1}, \es_{2}\},\,\parti\in \kn{\al}\}$ generates ${\cal S}(
G_{0}, D)$, we can deduce from (\ref{equ:4.193}), (\ref{eq:39}),
(\ref{eq:40}), (\ref{eq:45}) and (\ref{eq:46}) that
\begin{equation}\label{equ:4.198}
x^{\tau}\parti.v_{\rho_{\lmd}+\varrho}=(\varrho+\mu'_{k_{0}} -k_{0}\theta_{1})(\parti)v_{\rho_{\lmd}+\varrho+\tau}
\end{equation}
for $\varrho\in \G_{1}\backslash\{k_0\theta_{1}-\mu_{k_{0}}'\}$, $\tau\in G_0\backslash\{0\}$ and $\parti \in \kn{\tau}$, and
\begin{equation}\label{eq:30}
x^{\tau}\parti.v_{\rho_{\lmd}+k_0\theta_{1}-\mu_{k_{0}}'}=\eta_{k_{0}}' (\parti)v_{\rho_{\lmd}+k_0\theta_{1}-\mu_{k_{0}}'+\tau}\ \textrm{ for } \tau\in G_0\backslash\{0\} \textrm{ and }\parti \in \kn{\tau}.
\end{equation}

Third, we determine the action of $x^{\be'}\parti$ on $M(\rho_{\lmd},\G_1)$ for any $\be'\in \G_{1}\backslash \G_{0}$ and $\parti\in\kn \be'$. Fix any $\be'\in \G_{1}\backslash \G_{0}$. Recall that we can choose $\theta_{1}'\in (\theta_{1}+\G_{0})\cap G_{0}$ and $\es_{1}',\es_{2}' \in \G_{0}\backslash\{0\}$ such that $\kn \theta_{1}' \cap \kn \es_{1}' \cap \kn \es_{2}'=\{0\}$ and $\be'\in G_{1}=\mathbb{Z}\theta_{1}'+\mathbb{Z}\es_{1}'+\mathbb{Z}\es_{2}'$ (c.f. (\ref{eq:47})). Observe that the set ${\cal X}_1=\{x^{\pm \al}\parti \mid \al\in\{\theta_{1}',\, \es_{1}',\, \es_{2}'\},\,\parti\in \kn{\al}\}$ generates ${\cal S}( G_{1}, D)$, which is a subalgebra of ${\cal S}( \G,D )$ and also a simple generalized divergence-free
Lie algebra. As (\ref{equ:4.191}), (\ref{eq:38}), (\ref{equ:4.198}) and (\ref{eq:30}) give the action of ${\cal X}_1$ on $M(\rho_{\lmd},\G_1)$, we can derive
\begin{equation}\label{equ:4.199}
x^{\be'}\parti.v_{\rho_{\lmd}+\varrho}=(\varrho+\mu'_{k_{0}} -k_{0}\theta_{1})(\parti)v_{\rho_{\lmd}+\varrho+\be' }\ \textrm{ for } \varrho\in \G_{1}\backslash\{k_0\theta_{1}-\mu_{k_{0}}'\} \textrm{ and } \parti\in \kn{\be'}.
\end{equation}
\begin{equation}\label{eq:48}
x^{\be'}\parti.v_{\rho_{\lmd}+k_0\theta_{1}-\mu_{k_{0}}'}=\eta_{k_{0}}' (\parti)v_{\rho_{\lmd}+k_0\theta_{1}-\mu_{k_{0}}'+\be'}\ \textrm{ for } \parti \in \kn{\be'}.
\end{equation}

To sum up, (\ref{equ:4.191}), (\ref{eq:38}), (\ref{equ:4.199}) and
(\ref{eq:48}) show
\begin{equation}\label{equ:4.200}
M(\rho_{\lmd},\G_{1})\simeq \mathscr{A}_{\mu'_{k_{0}}-k_{0}\theta_{1},\eta_{k_{0}}' }({\cal S}(
\G_{1},D)), \ \textrm{ where } \mu'_{k_{0}}-k_{0}\theta_{1}\in \G_{1}\textrm{ and } \eta'_{k_0}\in D^{\ast}\backslash\{0\}.
\end{equation}

\vspace{0.2cm}

\noindent{\it{Case 3.}} \emph{There exists some $k_{0} \in K$ such
that
\begin{equation}\label{equ:4.201}
M(\rho_{\lmd}+k_{0}\theta_{1},\G_{0})\simeq \mathscr{M}_{\mu_{k_{0}}'}({\cal S}( \G_{0},D )) \textrm{ for some }\mu_{k_{0}}'\in \G_{0}.
\end{equation}}

In analogy with Case 2, we obtain
\begin{equation}\label{equ:4.203}
M(\rho_{\lmd},\G_{1})\simeq \mathscr{M}_{\mu'_{k_{0}}-k_{0}\theta_{1}}({\cal S}(
\G_{1},D)), \textrm{ where } \mu'_{k_{0}}-k_{0}\theta_{1}\in \G_{1}.
\end{equation}\vspace{0.1cm}

\noindent{\it{Case 4.}} \emph{There exists some $k_{0} \in K$ such
that
\begin{equation}\label{equ:4.202}
 M(\rho_{\lmd}+k_{0}\theta_{1},\G_{0})\simeq   \mathscr{B}_{\mu'_{k_{0}},\eta'_{k_{0}}}({\cal S}( \G_{0}, D)) \textrm{ for some }\mu_{k_{0}}'\in \G_{0}\textrm{ and } \eta'_{k_0}\in D^{\ast}\backslash\{0\}.
\end{equation}}

In analogy with Case 2, we obtain
\begin{equation}\label{equ:4.204}
 M(\rho_{\lmd},\G_{1})\simeq   \mathscr{B}_{\mu'_{k_{0}}-k_{0}\theta_{1},\eta'_{k_{0}}}({\cal S}( \G_{1}, D)), \ \textrm{ where } \mu'_{k_{0}}-k_{0}\theta_{1}\in \G_{1}\textrm{ and } \eta'_{k_0}\in D^{\ast}\backslash\{0\}..
\end{equation}

\vspace{0.1cm}

\noindent{\it{Case 5.}} \emph{For each $k \in K$,
\begin{equation}\label{equ:4.205}
M(\rho_{\lmd}+k\theta_{1},\G_{0})\simeq   \mathscr{M}_{\mu'_{k}}({\cal S}(
\G_{0}, D))  \textrm{ with some }\mu'_{k}\in D^{\ast}\backslash
\G_{0}.
\end{equation}}

In analogy with Case 2, we get
\begin{equation}
\mu'_{k}-k\theta_{1}=\mu'_{0}\quad \textrm{for } k\in K,
\end{equation}
and obtain
\begin{equation}\label{equ:4.206}
M(\rho_{\lmd},\G_{1})\simeq \mathscr{M}_{\mu'_{0}}({\cal S}(\G_{1},D)),
\textrm{ where }\mu'_{0}\in D^{\ast}\backslash \G_{1}.
\end{equation}
\vspace{0.2cm}

In summary, we get:

For any $\lmd\in J$, $ M(\rho_{\lmd},\G_{1})$ is isomorphic to
\begin{eqnarray}
  &(\rmnum{1})&\mathscr{M}_{\bar{\mu}_{\lmd}}({\cal S}( \G_{1},D)) \textrm{ for some } \bar{\mu}_{\lmd}\in D^{\ast}\backslash \G_{1}; \\
  &(\rmnum{2})& \mathscr{M}_{\bar{\mu}_{\lmd}}({\cal S}( \G_{1},D))\textrm{ for some } \bar{\mu}_{\lmd}\in \G_{1}; \\
   &(\rmnum{3})& \mathscr{A}_{\bar{\mu}_{\lmd},\bar{\eta}_{\lmd}}({\cal S}( \G_{1},D))\textrm{ for some }  \bar{\mu}_{\lmd}\in \G_{1}\textrm{ and } \bar{\eta}_{\lmd}\in D^{\ast}\backslash\{0\}; \\
   &(\rmnum{4})& \mathscr{B}_{\bar{\mu}_{\lmd},\bar{\eta}_{\lmd}}({\cal S}( \G_{1},D))\textrm{ for some }  \bar{\mu}_{\lmd}\in \G_{1}\textrm{ and } \bar{\eta}_{\lmd}\in D^{\ast}\backslash\{0\}; \\
   &(\rmnum{5})& \bigoplus_{\nu\in \G_{1}}\mathbb{F}w_{\nu},  \textrm{ where each component is a trivial submodule of }{\cal S}( \G_{1},D).
   \end{eqnarray}
Namely, $\G_1$ satisfies condition (C1), which contradicts the maximality of $\G_{0}$. So we must have $\G=\G_{0}$. Therefore, Lemma \ref{le:4.3} follows from (C1), Lemma \ref{le:3.1} and Theorem \ref{th:3.2}.  $\qquad \Box$

\section{The Case $\dim D\geq 4$}

In this section, we shall prove Theorem \ref{th:4.1} under the condition $\dim D\geq 4$. We always assume $\dim D\geq 4$ throughout this section.

Let $M=\bigoplus_{\theta\in \G}M_{\theta}$ be a $\G$-graded
${\cal S}(\G,D)$-module with $\textrm{dim} M_{\theta}=1$ for each $\theta\in \G$. In order to prove Theorem \ref{th:4.1} under the condition $\dim D\geq 4$, we first need to derive all the possible action of ${\cal S}(\G,D)$ on $M$. Picking any three $\mathbb{F}$-linearly independent elements $\parti_{1}, \parti_{2}, \parti_{3}\in D$, we set
\begin{equation}
D_{0}=\mathbb{F}\parti_{1}+\mathbb{F}\parti_{2}+\mathbb{F}\parti_{3}.
\end{equation}
To begin with, we have the following lemma. \vspace{0.2cm}

\begin{lemma}\label{le:6.1} The set
$\{x^{\al}\parti \mid \parti \in \kn \al \cap D_{0}, \al\in \G \textrm{ satisfying } D_{0}\not \subseteq \kn{\al}\}$ act on $M$ in one of the following two ways:

(P1) $x^{\al}\parti. M =\{0\}$
for $\parti \in \kn{\al}\cap D_{0}$ and $\al\in \G$ satisfying $D_{0}\not \subseteq \kn{\al}$;

(P2) there exist
$\mu \in D^{\ast}$ and $\{0\not=v_{\be}\in M_{\be}\mid \be\in \G \textrm{ satisfying }D_{0}\not \subseteq \kn (\be+\mu)\}$
such that
\begin{equation}
x^{\al}\parti.  v_{\be} =(\be+\mu)(\parti)v_{\al+\be}
\end{equation}
for $\parti \in \kn{\al}\cap D_{0}$, $\al\in \G$ satisfying $D_{0}\not \subseteq \kn{\al}$, and $\be \in \G$ satisfying $D_{0}\not \subseteq \kn (\be+\mu)$ and $D_{0}\not \subseteq \kn (\be+\mu+\al)$.
\end{lemma}

\noindent{\bf{Proof.}} Lemma
\ref{le:4.3} and Zorn's Lemma imply that there exists a maximal subgroup $\G_{1}$ of $\G$ such that

(\Rmnum{1}) $\G_{1}$ has a subgroup $\G_{0}$ that satisfies $(\cap_{\al \in \G_{0}}\kn \al)\cap D_{0} =\{0\}$,

 (\Rmnum{2}) $\{x^{\al}\parti \mid \parti \in \kn{\al}\cap D_{0}, \al\in \G_1 \textrm{ satisfying } D_{0}\not \subseteq \kn{\al}\}$ act on $\bigoplus_{\be\in\G_1}M_{\be}$ in one of the following two ways:

$(p1)$ $x^{\al}\parti. M_{\be} =\{0\}$ for $\be \in \G_{1}$, $\parti \in \kn{\al}\cap D_{0}$ and $\al\in \G_{1}$  satisfying $D_{0}\not \subseteq \kn{\al}$;

$(p2)$ there exist $\mu \in D^{\ast}$ and $\{0\not=v_{\be}\in
M_{\be}\mid \be\in \G_{1} \textrm{ satisfying } D_{0}\not \subseteq
\kn (\be+\mu)\}$ such that
\begin{equation}
x^{\al}\parti.  v_{\be} =(\be+\mu)(\parti)v_{\al+\be}
\end{equation}
for $\parti \in \kn{\al}\cap D_{0}$, $\al\in \G_{1}$ satisfying $D_{0}\not \subseteq \kn{\al}$, and $\be \in \G_{1}$ satisfying $D_{0}\not \subseteq \kn (\be+\mu)$ and $D_{0}\not \subseteq \kn (\be+\mu+\al)$.

To prove the lemma, it suffices to show $\G_{1}=\G$. Suppose that
$\G_{1}\not=\G$. We will see that this leads to a contradiction.

Choose $\theta_{1}\in \G\backslash \G_{1}$ such that $D_{0}\not \subseteq \kn{\theta_{1}}$. Let
\begin{equation}
\G_{2}=\G_{1}+\mathbb{Z}\theta_{1}.
\end{equation}
We immediately see that $\G_{2}$ satisfies (\Rmnum{1}). Moreover, we want to prove $\G_{2}$ satisfies (\Rmnum{2}).

Suppose that $\G_{1}$ satisfies $(p2)$ in (\Rmnum{2}); the case in which $\G_{1}$ satisfies $(p1)$ can be proved similarly. Namely, there exist $\mu \in D^{\ast}$
and $\{0\not=v_{\be}\in M_{\be} \mid\be\in \G_{1} \textrm{ satisfying } D_{0}\not \subseteq \kn (\be+\mu)\}$ such that
\begin{equation}\label{equ:4.213}
x^{\al}\parti.  v_{\be} =(\be+\mu)(\parti)v_{\al+\be}
\end{equation}
for $\parti \in \kn{\al}\cap D_{0}$, $\al\in \G_{1}$ satisfying $D_{0}\not \subseteq \kn{\al}$, and $\be \in \G_{1}$ satisfying $D_{0}\not \subseteq \kn (\be+\mu)$ and $D_{0}\not \subseteq \kn (\be+\mu+\al)$. We shall extend this to
\begin{equation}
 \textrm{ the action of }\{x^{\al}\parti \mid \parti \in \kn{\al}\cap D_{0}, \al\in \G_2 \textrm{ satisfying } D_{0}\not \subseteq \kn{\al}\}  \textrm{ on } \bigoplus_{\be\in\G_2}M_{\be}.
\end{equation}

The fact that $\G_{1}$ satisfies (\Rmnum{1}) enables us to choose $\es_{1}, \es_{2} \in \G_{1}$ such that
\begin{equation}\label{eq:49}
\kn{\theta_{1}} \cap \kn{\es_{1}}\cap \kn{\es_{2}} \cap D_{0}=\{0\}.
 \end{equation}
Set
\begin{equation}
G_{0}=\mathbb{Z}\theta_{1} + \mathbb{Z}\es_{1} + \mathbb{Z}\es_{2}.
\end{equation}
Then $(\cap_{\nu \in G_{0}} \kn \nu) \cap D_{0} =\{0\}$, and $G_{0}$ can be viewed as a subgroup of $D_{0}^{\ast}$. Hence we get a simple generalized divergence-free Lie algebra ${\cal S}(G_{0}, D_{0})$, which is also a subalgebra of ${\cal S}(\G,D)$.\vspace{0.2cm}

We then proceed our proof in several steps.\vspace{0.2cm}

{\it Firstly}, we extend the set $\{0\not=v_{\be}\in M_{\be} \mid
\be\in \G_{1}, D_{0}\not \subseteq \kn(\be+\mu)\}$ to $
\{0\not=v_{\be}\in M_{\be} \mid \be\in \G_2,D_{0}\not \subseteq
\kn(\be+\mu)\}, $ and determine the action of ${\cal S}(G_{0},
D_{0})$ on $\bigoplus_{\be\in\G_2}M_{\be}$.\vspace{0.2cm}

Let
\begin{equation}
\{\be_{\vs} \mid \vs\in \bar{I}\} \textrm{ be the
set of all representatives of cosets of } \G_1\cap G_0 \textrm{ in } \G_1.
\end{equation}
Fix any $\vs\in \bar{I}$. Lemma \ref{le:4.3} says that $ M(\be_{\vs},G_{0})$ is isomorphic to:
\begin{eqnarray}
   &(\rmnum{1})& \mathscr{M}_{\mu_{\vs}}({\cal S}( G_{0},D_{0} ))  \textrm{ for some }   \mu_{\vs}\in D_{0}^{\ast}\backslash G_{0};  \label{equ:8.12}\\
   &(\rmnum{2}) & \mathscr{M}_{\mu_{\vs}}({\cal S}( G_{0},D_{0} ))  \textrm{ for some }   \mu_{\vs}\in G_{0}; \\
    &(\rmnum{3})& \mathscr{A}_{\mu_{\vs},\eta_{\vs}}({\cal S}( G_{0},D_{0} ))  \textrm{ for some }   \mu_{\vs}\in G_{0}  \textrm{ and }  \eta_{\vs}\in D_{0}^{\ast}\backslash\{0\}; \\
    &(\rmnum{4})& \mathscr{B}_{\mu_{\vs},\eta_{\vs}}({\cal S}( G_{0},D_{0} ))  \textrm{ for some }   \mu_{\vs}\in G_{0} \textrm{ and } \eta_{\vs}\in D_{0}^{\ast}\backslash\{0\};  \label{eq:50}\\
    &(\rmnum{5})& \bigoplus_{\nu\in G_{0}}\mathbb{F}w_{\nu},\  \textrm{where each component is a trivial submodule of } {\cal S}( G_{0},D_{0} ).
   \end{eqnarray}
In the first four cases, $\mu_{\vs}$'s are respectively chosen so
that there exist nonzero $v'_{\be_{\vs}+\nu} \in M_{\be_{\vs}+\nu}$
with $\nu\in G_{0}$ satisfying $\nu|_{D_0}+\mu_{\vs}\not=0$ such
that
\begin{equation}\label{eq:51}
x^{\al}\parti.v'_{\be_{\vs}+\nu}=(\nu+\mu_{\vs})(\parti)v'_{\be_{\vs}+\nu+\al}
\end{equation}
for $\al\in G_{0}\backslash\{0\}$, $\parti \in \kn{\al}\cap D_{0}$ and
$\nu\in G_{0}$ satisfying $\nu|_{D_0}+\mu_{\vs}\not=0\not=\nu|_{D_0}+\mu_{\vs}+\al|_{D_0}$.

Recall that $\es_1,\es_2 \in\G_1$ (c.f. (\ref{eq:49})). Pick $\parti' \in
(\kn\es_{1}\cap D_0)\backslash(\kn\es_{2}\cap D_0)$.  Then (\ref{equ:4.213}) indicates
\begin{equation}\label{eq:55}
x^{\es_{1}}\parti'.v_{\be_{\vs}+l'\es_{1}+m'\es_{2}}\not=0 \ \textrm{ for some }l',m'\in \mathbb{Z},
\end{equation}
which implies,
\begin{equation}\label{eq:52}
 M(\be_{\vs},G_{0}) \textrm{ can only be isomorphic to one of the first four cases (i.e. ( \ref{equ:8.12})--(\ref{eq:50}))}.
 \end{equation}
So (\ref{eq:51}) shows
\begin{equation}
 x^{-\es_{s}}\parti.x^{\es_{s}}\parti.v_{\be_{\vs}+l\es_{1}+m\es_{2}} = (l\es_{1}+m\es_{2}+\mu_{\vs})^{2}(\parti)v_{\be_{\vs}+l\es_{1}+m\es_{2}}
\end{equation}
for $l,m\in \mathbb{Z}$ such that
$l\es_{1}|_{D_0}+m\es_{2}|_{D_0}+\mu_{\vs}\not=0 \not= l\es_{1}|_{D_0}+m\es_{2}|_{D_0}+\mu_{\vs}+\es_s|_{D_0}$, where $\parti \in \kn{\es_{s}}\cap D_0$ and $s\in\{1,2\}$. On the other hand, (\ref{equ:4.213}) tells
\begin{equation}
 x^{-\es_{s}}\parti.x^{\es_{s}}\parti.v_{\be_{\vs}+l\es_{1}+m\es_{2}} = (\be_{\vs}+l\es_{1}+m\es_{2}+\mu)^{2}(\parti)v_{\be_{\vs}+l\es_{1}+m\es_{2}}
\end{equation}
for $l,m\in \mathbb{Z}$ such that
$D_0\not \subseteq \kn (\be_{\vs}+l\es_{1}+m\es_{2}+\mu)$ and $D_0\not \subseteq  \kn (\be_{\vs}+l\es_{1}+m\es_{2}+\mu+\es_s)$, where $\parti \in \kn{\es_{s}}\cap D_0$ and $s\in\{1,2\}$. By the similar arguments as those from (\ref{eq:34}) to (\ref{equ:4.171}), we can prove
\begin{equation}\label{eq:56}
\mu_{\vs}=(\mu+\be_{\vs})|_{D_{0}} \quad \textrm{ for any }\vs \in \bar{I}.
\end{equation}
So this together with (\ref{eq:51}) and (\ref{eq:52}) enable us to
choose  nonzero $v'_{\be_{\vs}+\nu} \in M_{\be_{\vs}+\nu}$ with
$\nu\in G_{0}$ satisfying $D_{0}\not \subseteq \kn
(\nu+\mu+\be_{\vs})$ such that
\begin{equation}\label{equ:8.13}
x^{\al}\parti.v'_{\be_{\vs}+\nu}=(\nu+\mu+\be_{\vs})(\parti)v'_{\be_{\vs}+\nu+\al}
\end{equation}
for $\parti \in \kn{\al}\cap D_{0}$, $\al\in G_{0}\backslash\{0\}$, and
$\nu\in G_{0}$ satisfying $D_0\not \subseteq \kn (\nu+\mu+\be_{\vs})$ and $D_0\not \subseteq  \kn (\nu+\mu+\be_{\vs}+\al)$. Assume
\begin{equation}\label{eq:53}
v'_{\be_{\vs}+\ga}=\hat{c}_{\vs,\ga}v_{\be_{\vs}+\ga}\quad\textrm{
for }  \ga\in \G_1\cap G_{0} \textrm{ such that } D_0\not \subseteq  \kn (\be_{\vs}+\ga+\mu),
\end{equation}
where $\hat{c}_{\vs,\ga}\in\mbb{F}$ (c.f. (\ref{equ:4.213})). By the similar arguments as those from (\ref{equ:4.184}) to (\ref{equ:4.190}), we can prove that \begin{equation}\label{eq:54}
\hat{c}_{\vs,\ga}=\hat{c}_{\vs} \textrm{ is independent of }\ga.
\end{equation}

Since $\G_{2}=\G_{1}+\mathbb{Z}\theta_{1}=\G_{1}+G_0=\bigoplus_{\vs\in\bar{I}}(\be_{\vs}+G_0)$, multiplying
\begin{equation}
\{v'_{\be_{\vs}+\nu}\mid \nu\in G_0 \textrm{ satisfying }D_{0}\not \subseteq \kn(\be_{\vs}+\nu+\mu) \}
 \end{equation}
 by $\frac{1}{\hat{c}_{\vs}}$ for each $\vs\in\bar{I}$, we get a set
 \begin{equation}
\{0\not=v_{\be}\in M_{\be} \mid \be\in \G_2,D_{0}\not \subseteq \kn(\be+\mu)\},
\end{equation}
which expands $\{0\not=v_{\be}\in M_{\be} \mid \be\in \G_{1},
D_{0}\not \subseteq \kn(\be+\mu)\}$ by (\ref{eq:53}) and (\ref{eq:54}). Then
(\ref{equ:4.213}) still holds and (\ref{equ:8.13}) shows
\begin{equation}\label{equ:4.214}
x^{\al}\parti.  v_{\be}=(\be+\mu)(\parti)v_{\al+\be}
\end{equation}
for $\al\in G_{0}\backslash\{0\}$, $\parti \in \kn{\al}\cap D_{0}$ and $\be \in \G_{2}$ satisfying $D_{0}\not \subseteq \kn (\be+\mu)$ and $D_{0}\not \subseteq \kn (\be+\mu+\al)$.\vspace{0.2cm}

{\it Secondly}, we determine the action of $\{x^{\al}\parti \mid
\parti \in \kn{\al}\cap D_{0}, \al\in \G_1 \textrm{ such that }
D_{0}\not \subseteq \kn{\al} \textrm{ and } \kn{\al}\cap D_{0} \not
= \kn{\theta_{1}}\cap D_{0} \}$ on
$\bigoplus_{\be\in\G_2}M_{\be}$.\vspace{0.2cm}

Fix any $\al' \in \G_{1}$ such that $D_{0}\not \subseteq \kn{\al'}$ and $\kn{\al'}\cap D_{0} \not = \kn{\theta_{1}}\cap D_{0}$. Fix any $\varrho' \in \G_{2}$ such that
\begin{equation}\label{eq:60}
D_{0}\not \subseteq \kn (\varrho'+\mu)  \textrm{ and } D_{0}\not \subseteq \kn (\varrho'+\mu+\al').
\end{equation}
We want to see how $x^{\al'}\parti$ act on $v_{\varrho'}$ for $\parti \in \kn{\al'}\cap D_{0}$.

Since $\kn{\theta_{1}} \cap \kn{\al'}\cap \kn{\es_{1}} \cap D_{0}=\{0\}$ or $\kn{\theta_{1}} \cap \kn{\al'}\cap \kn{\es_{2}} \cap D_{0}=\{0\}$, without loss of generality, we assume that $\kn{\theta_{1}} \cap \kn{\al'}\cap \kn{\es_{1}} \cap D_{0}=\{0\}$. Let
\begin{equation}
G_{1}= \mathbb{Z}\theta_{1} + \mathbb{Z}\al' + \mathbb{Z}\es_1.
\end{equation}
Then $(\cap_{\nu \in G_{1}} \kn \nu)\cap D_{0} =\{0\}$ and $G_{1}$ can be viewed as a subgroup of $D_{0}^{\ast}$. Hence we get a simple generalized divergence-free Lie algebra ${\cal S}(G_1, D_{0})$, which is also a subalgebra of ${\cal S}(\G,D)$. By Lemma \ref{le:4.3}, we know that
$ M(\varrho',G_{1})$ is isomorphic to:
\begin{eqnarray}
   &(\rmnum{1})& \mathscr{M}_{\mu'_{\varrho'}}({\cal S}( G_{1},D_{0} ))\textrm{ for some }  \mu'_{\varrho'}\in D_{0}^{\ast}\backslash G_{1}; \label{equ:8.15}\\
   &(\rmnum{2})& \mathscr{M}_{\mu'_{\varrho'}}({\cal S}( G_{1},D_{0} )) \textrm{ for some }  \mu'_{\varrho'}\in G_{1}; \\
   &(\rmnum{3})& \mathscr{A}_{\mu'_{\varrho'},\eta'_{\varrho'}}({\cal S}( G_{1},D_{0} )) \textrm{ for some } \mu'_{\varrho'}\in G_{1}\textrm{ and } \eta'_{\varrho'}\in D_{0}^{\ast}\backslash\{0\}; \\
   &(\rmnum{4})& \mathscr{B}_{\mu'_{\varrho'},\eta'_{\varrho'}}({\cal S}( G_{1},D_{0} )) \textrm{ for some }  \mu'_{\varrho'}\in G_{1}\textrm{ and } \eta'_{\varrho'}\in D_{0}^{\ast}\backslash\{0\}; \label{eq:59}\\
   &(\rmnum{5})& \bigoplus_{\nu\in G_{1}}\mathbb{F}w_{\nu}, \textrm{ where each component is a trivial submodule of } {\cal S}( G_{1},D_{0} ).
   \end{eqnarray}
In the first four cases (c.f. ( \ref{equ:8.15})--(\ref{eq:59})),
$\mu'_{\varrho'}$'s are respectively chosen so that there exist
nonzero $w'_{\varrho'+\nu} \in M_{\varrho'+\nu}$ with $\nu\in G_{1}$
satisfying $\nu|_{D_0}+\mu'_{\varrho'}\not=0$ such that
\begin{equation}\label{eq:57}
x^{\tau}\parti.w'_{\varrho'+\nu}=(\nu+\mu'_{\varrho'})(\parti)w'_{\varrho'+\nu+\tau}
\end{equation}
for $\tau\in G_{1}\backslash\{0\}$, $\parti \in \kn{\tau}\cap D_{0}$ and
$\nu\in G_{1}$ satisfying $\nu|_{D_0}+\mu'_{\varrho'}\not=0\not=\nu|_{D_0}+\mu'_{\varrho'}+\tau|_{D_0}$. Note that (\ref{equ:4.214}) gives the action of $x^{\al}\parti$ on $\bigoplus_{\be\in\G_2}M_{\be}$ for $\parti\in\kn\al$ and $\al\in\{\pm\theta_1,\pm\es_1\}$. By the similar arguments as those from (\ref{eq:55}) to (\ref{eq:56}), we can prove that
\begin{equation}\label{eq:58}
 M(\varrho',G_1) \textrm{ is isomorphic to one of the first four cases (i.e. ( \ref{equ:8.15})--(\ref{eq:59}))},
 \end{equation}
 and
\begin{equation}\label{e:26}
\mu'_{\varrho'}=(\mu+\varrho')|_{D_{0}}.
\end{equation}

Pick nonzero $\tilde{\parti_{1}} \in \kn{\theta_{1}} \cap
\kn{\al'}\cap D_{0}$ and $\tilde{\parti_{2}} \in
\kn(\theta_{1}+\es_{1}) \cap \kn{\al'}\cap D_{0}$. Then
$\{\tilde{\parti_{1}},\tilde{\parti_{2}}\}$ forms a basis of
$\kn{\al'}\cap D_{0}$. If $(\varrho'+\mu)(\tilde{\parti_{1}})=0$, then
(\ref{eq:60}) and (\ref{eq:57})--(\ref{e:26}) indicate
\begin{equation}\label{eq:61}
x^{\al'}\tilde{\parti_{1}}.  v_{\varrho'}=0 =(\varrho'+\mu)(\tilde{\parti_{1}})v_{\al'+\varrho'}.
\end{equation}
Assume $(\varrho'+\mu)(\tilde{\parti_{1}})\not=0$. Since $\G_{2}=\G_{1}+\mathbb{Z}\theta_{1}$ and $\varrho' \in \G_{2}$, we can write $\varrho'=\be'+k'\theta_{1}$ for some $\be' \in \G_{1}$ and $k' \in\mbb{Z}$. As $\tilde{\parti_{1}} \in \kn{\theta_{1}} \cap \kn{\al'}\cap D_{0}$, we have
\begin{equation}
(\be'+\mu)(\tilde{\parti_{1}})=(\varrho'+\mu)(\tilde{\parti_{1}})\not=0 \textrm{ and }(\al'+\be'+\mu)(\tilde{\parti_{1}})=(\varrho'+\mu)(\tilde{\parti_{1}})\not=0;
\end{equation}
namely, $D_{0}\not \subseteq \kn (\be'+\mu)$ and $D_{0}\not \subseteq \kn (\be'+\mu+\al')$.
So by (\ref{equ:4.213}) and (\ref{equ:4.214}), we have
\begin{equation}\label{equ:4.215}
(\be'+\mu)(\tilde{\parti_{1}})x^{\al'}\tilde{\parti_{1}}.  v_{\varrho'} =x^{\al'}\tilde{\parti_{1}}.x^{k'\theta_1}\tilde{\parti_{1}}.v_{\be'}=x^{k'\theta_1}\tilde{\parti_{1}}.x^{\al'}\tilde{\parti_{1}}.v_{\be'}=(\be'+\mu)^{2}(\tilde{\parti_{1}})v_{\al'+\varrho'},
\end{equation}
which implies
\begin{equation}\label{eq:62}
x^{\al'}\tilde{\parti_{1}}.v_{\varrho'}=(\be'+\mu)(\tilde{\parti_{1}})v_{\al'+\varrho'}=(\varrho'+\mu)(\tilde{\parti_{1}})v_{\al'+\varrho'}.
\end{equation}
Replacing $\tilde{\parti_{1}}$ by $\tilde{\parti_{2}}$ and replacing $\theta_{1}$ by $\theta_{1}+\es_{1}$ in (\ref{eq:61})--(\ref{eq:62}), we can similarly prove
\begin{equation}
x^{\al'}\tilde{\parti_{2}}.  v_{\varrho'} =(\varrho'+\mu)(\tilde{\parti_{2}})v_{\al'+\varrho'}.
\end{equation}
So we obtain
\begin{equation}\label{equ:8.16}
x^{\al'}\parti.  v_{\varrho'} =(\varrho'+\mu)(\parti)v_{\al'+\varrho'}
\end{equation}
for $\parti \in \kn{\al'}\cap D_{0}$, $\al' \in \G_{1}$ satisfying $D_{0}\not \subseteq \kn{\al'}$ and $\kn{\al'}\cap D_{0}\not = \kn{\theta_{1}}\cap D_{0}$, and $\varrho' \in \G_{2}$ satisfying $D_{0}\not \subseteq \kn (\varrho'+\mu)$ and $D_{0}\not \subseteq \kn (\varrho'+\mu+\al')$.\vspace{0.2cm}

 {\it Thirdly}, we determine the action of $\{x^{\al}\parti \mid
\parti \in \kn{\al}\cap D_{0}, \al\in \G_1 \textrm{ such that }
D_{0}\not \subseteq \kn{\al} \textrm{ and }\\ \kn{\al}\cap D_{0}  =
\kn{\theta_{1}}\cap D_{0} \}$ on
$\bigoplus_{\be\in\G_2}M_{\be}$.\vspace{0.2cm}

Fix any $\al' \in \G_{1}$ satisfying $D_{0}\not \subseteq \kn{\al'}$ and $\kn{\al'}\cap D_{0} = \kn{\theta_{1}}\cap D_{0}$. Fix any $\varrho' \in \G_{2}$ such that
\begin{equation}\label{eq:63}
D_{0}\not \subseteq \kn (\varrho'+\mu)  \textrm{ and } D_{0}\not \subseteq \kn (\varrho'+\mu+\al').
\end{equation}
We want to see how $x^{\al'}\parti$ act on $v_{\varrho'}$ for any $\parti \in \kn{\al'}\cap D_{0}$.

Since $\kn{\al'}\cap D_{0} = \kn{\theta_{1}}\cap D_{0}$, we have $\kn{\al'} \cap \kn{\es_1}\cap \kn{\es_{2}} \cap D_{0}=\{0\}$. Let
\begin{equation}
G_2= \mathbb{Z}\al' + \mathbb{Z}\es_1 + \mathbb{Z}\es_2.
\end{equation}
Then $(\cap_{\nu \in G_2} \kn \nu)\cap D_{0} =\{0\}$ and $G_2$ can be viewed as a subgroup of $D_{0}^{\ast}$.  Hence we get a simple generalized divergence-free Lie algebra ${\cal S}(G_2, D_{0})$, which is also a subalgebra of ${\cal S}(\G,D)$.
Note that (\ref{equ:4.214}) gives the action of $x^{\pm\es_s}\parti$ on $\bigoplus_{\be\in\G_2}M_{\be}$ for $\parti\in\kn\es_s$ and $s\in\{1,2\}$.
 By the similar arguments as those from (\ref{equ:8.15}) to (\ref{equ:8.16}), we can prove that
\begin{equation}\label{equ:8.18}
x^{\al'}\parti.  v_{\varrho'} =(\varrho'+\mu)(\parti)v_{\al'+\varrho'}
\end{equation}
for $\parti \in \kn{\al'}\cap D_{0}$, $\al' \in \G_{1}$ satisfying $D_{0}\not \subseteq \kn{\al'}$ and $\kn{\al'}\cap D_{0} = \kn{\theta_{1}}\cap D_{0}$, and $\varrho' \in \G_{2}$ satisfying $D_{0}\not \subseteq \kn (\varrho'+\mu)$ and $D_{0}\not \subseteq \kn (\al'+\varrho'+\mu)$.\vspace{0.2cm}

 {\it Fourthly}, we determine the action of $\{x^{\al}\parti \mid
\parti \in \kn{\al}\cap D_{0}, \al\in \G_2\backslash\G_1 \textrm{
satisfying } D_{0}\not \subseteq \kn{\al}\}$ on
$\bigoplus_{\be\in\G_2}M_{\be}$.\vspace{0.2cm}

Fix any $\al' \in \G_{2}\backslash\G_1$ satisfying $D_{0}\not \subseteq \kn{\al'}$. Write $\al'=k_0 \theta_{1} + \al_0$
 for some $k_0 \in\mathbb{Z}\backslash\{0\}$ and $\al_0\in \G_{1}$. If $\kn{\al'}\cap D_{0} \not = \kn{\theta_{1}}\cap D_{0}$, then $\al'=k_0 \theta_{1} + \al_0$ implies $D_{0}\not \subseteq\kn{\al_0}$ and $\kn{\al_0}\cap D_{0} \not = \kn{\theta_{1}}\cap D_{0}$. Thus we have $\kn{\theta_{1}} \cap \kn{\al_0}\cap \kn{\es_1} \cap D_{0}=\{0\}$ or $\kn{\theta_{1}} \cap \kn{\al_0}\cap \kn{\es_2} \cap D_{0}=\{0\}$. Without loss of generality, we assume that $\kn{\theta_{1}} \cap \kn{\al_0}\cap \kn{\es_1} \cap D_{0}=\{0\}$. Let
\begin{equation}
G_{3}= \mathbb{Z}\theta_{1} + \mathbb{Z}\es_1 + \mathbb{Z}\al_0.
\end{equation}
Then $(\cap_{\nu \in G_{3}} \kn \nu)\cap D_{0} =\{0\}$ and $G_{3}$ can be viewed as a subgroup of $D_{0}^{\ast}$. Hence we get a simple generalized divergence-free Lie algebra ${\cal S}(G_3, D_{0})$, which is also a subalgebra of ${\cal S}(\G,D)$. Note that $\al'\in G_{3}$. Moreover, since ${\cal X}_3=\{x^{\pm \sigma}\parti \mid \sigma\in\{\theta_{1}, \es_1, \al_0\},\,\parti\in \kn{\sigma}\cap D_{0}\}$ generates ${\cal S}( G_3, D_{0})$, and (\ref{equ:4.214}), (\ref{equ:8.16}) and (\ref{equ:8.18}) give the action of ${\cal X}_3$ on ${\cal S}(G_3, D_{0})$-module $\bigoplus_{\be\in\G_2}M_{\be}$, we can deduce
\begin{equation}\label{equ:4.217}
x^{\al'}\parti.  v_{\be} =(\be+\mu)(\parti)v_{\al'+\be}
\end{equation}
for $\parti \in \kn{\al'}\cap D_{0}$, $\be \in \G_{2}$ such that $D_{0}\not \subseteq \kn (\be+\mu)$ and $D_{0}\not \subseteq \kn (\be+\mu+\al')$.

If $\kn{\al'}\cap D_{0} = \kn{\theta_{1}}\cap D_{0}$, then $\al'=k_0 \theta_{1} + \al_0$ implies $D_{0}\not \subseteq\kn(\al_0-k_0 \es_1)$ and $\kn(\theta_{1}+\es_1) \cap D_{0}\not= \kn(\al_0-k_0 \es_1)\cap D_{0}$. Choose $\es_3 \in \G_{1}$ such that
\begin{equation}
\kn(\theta_{1}+\es_1) \cap \kn(\al_0-k_0 \es_1)\cap \kn{\es_3}\cap D_{0}=\{0\}.
\end{equation}
 Let
\begin{equation}
G_4= \mathbb{Z}(\theta_{1}+\es_1) + \mathbb{Z}(\al_0-k_0 \es_1) + \mathbb{Z}\es_3.
\end{equation}
Then $(\cap_{\nu \in G_{4}} \kn \nu)\cap D_{0} =\{0\}$ and $G_{4}$ can be viewed as a subgroup of $D_{0}^{\ast}$. Hence we get a simple generalized divergence-free Lie algebra ${\cal S}(G_4, D_{0})$, which is also a subalgebra of ${\cal S}(\G,D)$. Note that $\al'\in G_4$. Moreover, since the set \begin{equation}
{\cal X}_4=\{x^{\pm \sigma}\parti \mid \sigma\in\{\theta_{1}+\es_1, \al_0-k_0 \es_1, \es_3\},\,\parti\in \kn{\sigma}\cap D_{0}\}
 \end{equation}
 generates ${\cal S}(G_{4}, D_{0})$, and (\ref{equ:4.214}), (\ref{equ:8.16}) and (\ref{equ:8.18}) give the action of ${\cal X}_4$ on ${\cal S}(G_4, D_{0})$-module $\bigoplus_{\be\in\G_2}M_{\be}$, we can deduce
\begin{equation}\label{equ:4.218}
x^{\al'}\parti.  v_{\be} =(\be+\mu)(\parti)v_{\al'+\be}
\end{equation}
for $\parti \in \kn{\al'}\cap D_{0}$, $\be \in \G_{2}$ such that $D_{0}\not \subseteq \kn (\be+\mu)$ and $D_{0}\not \subseteq \kn (\be+\mu+\al')$.\vspace{0.2cm}

{\it To sum up}, (\ref{equ:8.16}), (\ref{equ:8.18}), (\ref{equ:4.217}) and (\ref{equ:4.218}) show that $\G_2$ satisfies $(p2)$ in (\Rmnum{2}). This contradicts the maximality of $\G_{1}$. On the other hand, if $\G_1$ satisfies $(p1)$ in (\Rmnum{2}), we can similarly prove that $\G_2$ satisfies $(p1)$ in (\Rmnum{2}), which also contradicts the maximality of $\G_{1}$; we omit the details here. So we must have $\G_1=\G$. Therefore the lemma follows. $\qquad\qquad \Box$

\begin{lemma}\label{le:6.2}
 ${\cal S}(\G,D)$ act on $M$ in one of the following two ways:

$(P'1)$ $x^{\al}\parti. M =\{0\}$
for $\parti \in \kn{\al}$ and $\al\in \G\backslash\{0\}$;

$(P'2)$ there exist $\mu \in D^{\ast}$ and $\{0\not=v_{\be}\in
M_{\be}\mid \be\in \G\backslash\{-\mu\}\}$ such that
\begin{equation}
x^{\al}\parti.  v_{\be} =(\be+\mu)(\parti)v_{\al+\be}
\end{equation}
for $\parti \in \kn{\al}$, $\al\in \G\backslash\{0\}$ and $\be \in \G$ satisfying $\be+\mu\not=0\not=\be+\mu+\al$.
\end{lemma}

\noindent{\bf Proof.}
Lemma \ref{le:6.1} and Zorn's Lemma imply that, there
exists a maximal subspace $D'\subseteq D$ such that

$(\Rmnum{1}')$  $\quad\dim D' \geqslant 3$,

$(\Rmnum{2}')$ $\{x^{\al}\parti \mid \parti \in \kn \al \cap D', \al\in \G \textrm{ satisfying } D'\not \subseteq \kn{\al}\}$ act on $M$ in one of the following two ways:

$(p'1)$ $x^{\al}\parti. M =\{0\}$
for $\parti \in \kn{\al}\cap D'$ and $\al\in \G$ satisfying $D'\not \subseteq \kn{\al}$;

$(p'2)$ there exist
$\mu \in D^{\ast}$ and $\{0\not=v_{\be}\in M_{\be} \mid \be\in \G\textrm{ satisfying } D'\not \subseteq \kn (\be+\mu)\}$
such that
\begin{equation}
x^{\al}\parti.  v_{\be} =(\be+\mu)(\parti)v_{\al+\be}
\end{equation}
for $\parti \in \kn{\al}\cap D'$, $\al\in \G$ satisfying $D'\not \subseteq \kn{\al}$, and $\be \in \G$ satisfying $D'\not \subseteq \kn (\be+\mu)$ and $D'\not \subseteq \kn (\be+\mu+\al)$.

To prove the lemma, it suffices to show $D'=D$. Suppose $D'\not=D$.
We will see that this leads to a contradiction.

Pick $\bar{\parti}\in D \backslash D'$. Set
\begin{equation}\label{eq:67}
D''=D'+\mathbb{F}\bar{\parti}.
\end{equation}
Observe that $D''$ satisfies $(\Rmnum{1}')$. We will also show $D''$ satisfies $(\Rmnum{2}')$.

Assume $D'$ satisfies $(p'2)$ in $(\Rmnum{2}')$; the case in which
$D'$ satisfies $(p'1)$ can be proved similarly. Namely, we can
choose $\mu \in D^{\ast}$ and $\{0\not=v_{\be}\in M_{\be} \mid
\be\in \G \textrm{ satisfying }D'\not \subseteq \kn (\be+\mu)\}$
such that
\begin{equation}\label{equ:4.222}
x^{\al}\parti.  v_{\be} =(\be+\mu)(\parti)v_{\al+\be}
\end{equation}
for $\parti \in \kn{\al}\cap D'$, $\al\in \G$ satisfying $D'\not \subseteq \kn{\al}$, and $\be \in \G$ satisfying $D'\not \subseteq \kn (\be+\mu)$ and $D'\not \subseteq \kn (\be+\mu+\al)$. Then we proceed our proof in several steps.\vspace{0.2cm}

{\it Firstly}, we give a claim:\vspace{0.2cm}

Pick any two linearly independent elements $\bar{\parti}_{1}, \bar{\parti}_{2} \in D'$. Set
\begin{equation}\label{eq:78}
\tilde{D}_{0}=\mathbb{F}\bar{\parti}+\mathbb{F}\bar{\parti}_{1}+\mathbb{F}\bar{\parti}_{2}.
\end{equation}
Lemma \ref{le:6.1} and (\ref{equ:4.222}) imply that there exist
$\mu' \in D^{\ast}$ and $\{0\not=u_{\be}\in M_{\be}\mid \be\in
\G\textrm{ satisfying }\tilde{D}_{0}\not\subseteq \kn(\be+\mu')\}$
such that
\begin{equation}\label{equ:4.223}
x^{\al}\parti.  u_{\be} =(\be+\mu')(\parti)u_{\al+\be}
\end{equation}
for $\parti \in \kn{\al}\cap \tilde{D}_{0}$, $\al\in \G$ satisfying $\tilde{D}_{0}\not \subseteq \kn{\al}$, and $\be \in \G$ satisfying $\tilde{D}_{0}\not \subseteq \kn (\be+\mu')$ and $\tilde{D}_{0}\not \subseteq \kn (\be+\mu'+\al)$. Similar arguments as those from (\ref{equ:4.169}) to (\ref{equ:6.3}), combined with (\ref{equ:4.222}) and (\ref{equ:4.223}), indicate
\begin{equation}\label{eq:68}
\mu'(\parti)=\mu(\parti) \ \textrm{ for }\parti \in \mbb{F}\bar{\parti}_{1}+\mbb{F}\bar{\parti}_{2}.
\end{equation}
Assume
\begin{equation}\label{eq:66}
u_{\be}=\tilde{a}_{\be}v_{\be} \quad \textrm{ for } \be\in \G \textrm{ satisfying }\mathbb{F}\bar{\parti}_{1}+\mathbb{F}\bar{\parti}_{2}\not\subseteq \kn(\be+\mu),
\end{equation}
  where $\tilde{a}_\be\in\mathbb{F}$.\vspace{0.3cm}

 \emph{Claim 1}. {\it  $\tilde{a}_{\be}=\tilde{a}$ is independent of $\be$. Moreover,
\begin{equation}
x^{\al}\parti.  v_{\be} =(\be+\mu')(\parti)v_{\al+\be}
\end{equation}
for $\parti \in \kn\al\cap \tilde{D}_{0}$, $\al\in \G$ satisfying
$\tilde{D}_{0}\not \subseteq \kn{\al}$, and $\be \in \G$ satisfying
$\mathbb{F}\bar{\parti}_{1}+\mathbb{F}\bar{\parti}_{2}\not\subseteq
\kn(\be+\mu)$ and
$\mathbb{F}\bar{\parti}_{1}+\mathbb{F}\bar{\parti}_{2}\not \subseteq
\kn (\be+\mu+\al)$.}\vspace{0.2cm}

Choose $\rho_{0} \in \G$ such that $(\rho_{0}+\mu)(\bar{\parti}_{1})\not=0$. Write
\begin{equation}
\parti'=(\rho_{0}+\mu)(\bar{\parti}_{1})\bar{\parti}_{2}-(\rho_{0}+\mu)(\bar{\parti}_{2})\bar{\parti}_{1}.
 \end{equation}
Then $(\rho_{0}+\mu)(\parti')=0$. For any $\al\in \G$ such that $ \al(\parti')\not=0$, making use of (\ref{eq:68}),
 we  deduce from (\ref{equ:4.222}) and (\ref{equ:4.223}) that
\begin{equation}\label{equ:4.224}
x^{\al}(\al(\parti')\bar{\parti}_{1}-\al(\bar{\parti}_{1})\parti').  v_{\rho_{0}}  =  \al(\parti')(\rho_{0}+\mu)(\bar{\parti}_{1})v_{\al+\rho_{0}}
\end{equation}
and
\begin{equation}\label{equ:4.225}
x^{\al}(\al(\parti')\bar{\parti}_{1}-\al(\bar{\parti}_{1})\parti').  u_{\rho_{0}}  =  \al(\parti')(\rho_{0}+\mu)(\bar{\parti}_{1})u_{\al+\rho_{0}},
\end{equation}
which implies
\begin{equation}\label{eq:64}
\tilde{a}_{\al+\rho_{0}}=\tilde{a}_{\rho_{0}}.
\end{equation}
Such $\al $'s do exist and we denote one of them by $\al_{0}$. For
any $\al'\in \G$ such that $ \al'(\parti')=0$ but
$(\al'+\rho_{0}+\mu)(\bar{\parti}_{1})\not=0$, we have $
(\al'-\al_{0})(\parti')\not=0$. Moreover, making use of
(\ref{eq:68}),  we derive from (\ref{equ:4.222}) and
(\ref{equ:4.223}) that
\begin{equation}\label{equ:4.226}
x^{\al'-\al_{0}}(\al_{0}(\parti')\bar{\parti}_{1}+(\al'-\al_{0})(\bar{\parti}_{1})\parti').  v_{\rho_{0}+\al_{0}}  =  \al_{0}(\parti')(\al'+\rho_{0}+\mu)(\bar{\parti}_{1})v_{\al'+\rho_{0}}
\end{equation}
and
\begin{equation}\label{equ:4.227}
x^{\al'-\al_{0}}(\al_{0}(\parti')\bar{\parti}_{1}+(\al'-\al_{0})(\bar{\parti}_{1})\parti').  u_{\rho_{0}+\al_{0}}  =  \al_{0}(\parti')(\al'+\rho_{0}+\mu)(\bar{\parti}_{1})u_{\al'+\rho_{0}},
\end{equation}
which implies
\begin{equation}\label{eq:65}
\tilde{a}_{\al'+\rho_{0}}=\tilde{a}_{\al_{0}+\rho_{0}}=\tilde{a}_{\rho_{0}}.
\end{equation}
Since $\{\parti',\bar{\parti}_{1}\}$ forms a basis of $\mathbb{F}\bar{\parti}_{1}+\mathbb{F}\bar{\parti}_{2}$, (\ref{eq:64}) and (\ref{eq:65}) show
\begin{equation}
\tilde{a}_{\be}=\tilde{a}\quad\textrm{ for } \be\in \G \textrm{ such that }\mathbb{F}\bar{\parti}_{1}+\mathbb{F}\bar{\parti}_{2}\not\subseteq \kn(\be+\mu).
\end{equation}
So (\ref{equ:4.223}) and (\ref{eq:66}) give
\begin{equation}\label{equ:4.229}
x^{\al}\parti.  v_{\be} =(\be+\mu')(\parti)v_{\al+\be}
\end{equation}
for $\parti \in \kn{\al}\cap \tilde{D}_{0}$, $\al\in \G$ satisfying $\tilde{D}_{0}\not \subseteq \kn{\al}$, and $\be \in \G$ satisfying $\mathbb{F}\bar{\parti}_{1}+\mathbb{F}\bar{\parti}_{2}\not\subseteq \kn(\be+\mu)$ and $\mathbb{F}\bar{\parti}_{1}+\mathbb{F}\bar{\parti}_{2}\not \subseteq \kn (\be+\mu+\al)$. Thus the claim follows.\vspace{0.3cm}

{\it Secondly}, we choose some $\mu_0\in D^{\ast}$ such that
\begin{equation}\label{eq:77}
\mu_0|_{D'}=\mu|_{D'} \textrm{ and } \mu_0(\bar{\parti})=\mu'(\bar{\parti}),
\end{equation}
where $\mu'$ is the same as in (\ref{equ:4.223}).

\vspace{0.2cm}

{\it Thirdly}, we give another claim:\vspace{0.2cm}

Pick another two linearly independent elements $\bar{\parti}_{1}', \bar{\parti}_{2}' \in D'$. Set
\begin{equation}\label{eq:79}
\tilde{D}_{0}'=\mathbb{F}\bar{\parti}+\mathbb{F}\bar{\parti}_{1}'+\mathbb{F}\bar{\parti}_{2}'.
\end{equation}
Lemma \ref{le:6.1} and (\ref{equ:4.222}) imply that there exist
$\mu'' \in D^{\ast}$ and $\{0\not=u'_{\be}\in M_{\be} \mid \be\in
\G\textrm{ satisfying }\tilde{D}_{0}'\not\subseteq \kn(\be+\mu'')\}$
such that
\begin{equation}\label{eq:69}
x^{\al}\parti.  u'_{\be} =(\be+\mu'')(\parti)u'_{\al+\be}
\end{equation}
for $\parti \in \kn{\al}\cap \tilde{D}_{0}'$, $\al\in \G$ satisfying $\tilde{D}_{0}'\not \subseteq \kn{\al}$, and $\be \in \G$ satisfying $\tilde{D}_{0}'\not \subseteq \kn (\be+\mu'')$ and $\tilde{D}_{0}'\not \subseteq \kn (\be+\mu''+\al)$.\vspace{0.3cm}

 \emph{Claim 2}. {\it  $\mu''|_{\tilde{D}_{0}'}=\mu_0|_{\tilde{D}_{0}'}$.}\vspace{0.2cm}

Similar arguments as those from (\ref{equ:4.169}) to (\ref{equ:6.3}), combined with (\ref{equ:4.222}), (\ref{eq:77}) and (\ref{eq:69}), imply
\begin{equation}\label{eq:70}
\mu''(\parti)=\mu(\parti)=\mu_0(\parti) \ \textrm{ for }\parti \in \mbb{F}\bar{\parti}_{1}'+\mbb{F}\bar{\parti}_{2}'.
\end{equation}
Moreover, Claim 1 indicates that
\begin{equation}\label{eq:71}
x^{\al}\parti.  v_{\be} =(\be+\mu'')(\parti)v_{\al+\be}
\end{equation}
for $\parti \in \kn{\al}\cap \tilde{D}_{0}'$, $\al\in \G$ satisfying $\tilde{D}_{0}'\not \subseteq \kn{\al}$, and $\be \in \G$ satisfying $\mathbb{F}\bar{\parti}_{1}'+\mathbb{F}\bar{\parti}_{2}'\not\subseteq \kn(\be+\mu)$ and $\mathbb{F}\bar{\parti}_{1}'+\mathbb{F}\bar{\parti}_{2}'\not \subseteq \kn (\be+\mu+\al)$.

Pick nonzero $\parti \in \mbb{F}\bar{\parti}_{1}+\mbb{F}\bar{\parti}_{2}$ and $\parti' \in \mbb{F}\bar{\parti}_{1}'+\mbb{F}\bar{\parti}_{2}'$ (c.f. (\ref{eq:78}), (\ref{eq:79})). Choose $\al\in\G$ such that $\al(\parti)\not=0$ and $\al(\parti')\not=0$. Then
\begin{equation}
\bar{\parti}-\frac{\al(\bar{\parti})}{\al(\parti)}\parti \in \kn\al \cap\tilde{D}_{0}, \quad  \bar{\parti}-\frac{\al(\bar{\parti})}{\al(\parti')}\parti' \in \kn\al \cap\tilde{D}_{0}'.
\end{equation}
Moreover, we choose $\be\in\G$ such that $\mathbb{F}\bar{\parti}_{1}+\mathbb{F}\bar{\parti}_{2}\not\subseteq \kn(\be+\mu)$, $\mathbb{F}\bar{\parti}_{1}'+\mathbb{F}\bar{\parti}_{2}'\not\subseteq \kn(\be+\mu)$, $\mathbb{F}\bar{\parti}_{1}+\mathbb{F}\bar{\parti}_{2}\not \subseteq \kn (\be+\mu+\al)$ and $\mathbb{F}\bar{\parti}_{1}'+\mathbb{F}\bar{\parti}_{2}'\not \subseteq \kn (\be+\mu+\al)$. Then (\ref{equ:4.229}) and (\ref{eq:71}) show
\begin{equation}\label{eq:72}
x^{\al}(\bar{\parti}-\frac{\al(\bar{\parti})}{\al(\parti)}\parti ).  v_{\be} =(\be+\mu')(\bar{\parti}-\frac{\al(\bar{\parti})}{\al(\parti)}\parti )v_{\al+\be},
\end{equation}
\begin{equation}\label{eq:73}
x^{\al}(\bar{\parti}-\frac{\al(\bar{\parti})}{\al(\parti')}\parti').  v_{\be} =(\be+\mu'')(\bar{\parti}-\frac{\al(\bar{\parti})}{\al(\parti')}\parti')v_{\al+\be}.
\end{equation}
Making use of (\ref{eq:68}) and (\ref{eq:70}), we deduce that the difference between (\ref{eq:72}) and (\ref{eq:73}) is
\begin{equation}\label{eq:74}
\al(\bar{\parti})x^{\al}(\frac{\parti'}{\al(\parti')}-\frac{\parti}{\al(\parti)}).  v_{\be} =\big( (\mu'-\mu'')(\bar{\parti})+\al(\bar{\parti})\cdot(\be+\mu)(\frac{\parti'}{\al(\parti')}-\frac{\parti}{\al(\parti)})\big)v_{\al+\be}.
\end{equation}
On the other hand, (\ref{equ:4.222}) shows
\begin{equation}\label{eq:75}
x^{\al}(\frac{\parti'}{\al(\parti')}-\frac{\parti}{\al(\parti)}).  v_{\be} =(\be+\mu)(\frac{\parti'}{\al(\parti')}-\frac{\parti}{\al(\parti)})v_{\al+\be}.
\end{equation}
Inserting (\ref{eq:75}) into (\ref{eq:74}) and making use of (\ref{eq:77}), we get
\begin{equation}\label{eq:76}
\mu''(\bar{\parti})=\mu'(\bar{\parti})=\mu_0(\bar{\parti}).
\end{equation}
Since $\tilde{D}_{0}'=\mathbb{F}\bar{\parti}+\mathbb{F}\bar{\parti}_{1}'+\mathbb{F}\bar{\parti}_{2}'$, this claim follows from (\ref{eq:70}) and (\ref{eq:76}).\vspace{0.3cm}

{\it Fourthly, we extend the set $\{0\not=v_{\be}\in M_{\be}\mid \be\in \G \textrm{ satisfying }D'\not \subseteq \kn (\be+\mu)\}$ to $\{0\not=v_{\be}\in M_{\be}\mid \be\in \G \textrm{ satisfying }D''\not \subseteq \kn (\be+\mu_0)\}$.}\vspace{0.2cm}

Since (\ref{equ:4.222}) gives the set
\begin{equation}
\{0\not=v_{\be}\in M_{\be}\mid \be\in \G \textrm{ satisfying }D'\not \subseteq \kn (\be+\mu)\}
\end{equation}
 and (\ref{eq:77}) shows
\begin{equation}
\{\be\in \G \mid D'\not \subseteq \kn (\be+\mu)\}=\{\be\in \G \mid D'\not \subseteq \kn (\be+\mu_0)\},
\end{equation}
we take $v_{\be}$ for $\be\in \G $ satisfying $D'\not \subseteq \kn (\be+\mu_0)$ as they were in (\ref{equ:4.222}).

Recall that $D''=D'+\mathbb{F}\bar{\parti}$ (c.f. (\ref{eq:67})). Then we only need to determine $v_{\be}$ for $\be\in \G$ such that $D'\subseteq \kn(\be+\mu_0)$ and $(\be+\mu_0)(\bar{\parti})\not=0$.

Fix two linearly independent elements $\hat{\parti}_{1},\hat{\parti}_{2}\in D'$. Set
\begin{equation}\label{eq:83}
\hat{D}=\mathbb{F}\bar{\parti}+\mathbb{F}\hat{\parti}_{1}+\mathbb{F}\hat{\parti}_{2}.
\end{equation}
Then Lemma \ref{le:6.1}, Claim 2 and (\ref{equ:4.222}) imply that
there exist $\{0\not=\hat{w}_{\be}\in M_{\be}\mid \be\in \G\textrm{
satisfying } \hat{D}\not\subseteq \kn(\be+\mu_0)\}$ such that
\begin{equation}\label{eq:95}
x^{\al}\parti. \hat{w}_{\be} =(\be+\mu_0)(\parti)\hat{w}_{\al+\be}
\end{equation}
for $\parti \in \kn{\al}\cap \hat{D}$, $\al\in \G$ satisfying $\hat{D}\not \subseteq \kn{\al}$, and $\be \in \G$ satisfying $\hat{D}\not \subseteq \kn (\be+\mu_0)$ and $\hat{D}\not \subseteq \kn (\be+\mu_0+\al)$. Moreover, Claim 1 shows that
\begin{equation}\label{eq:84}
\hat{w}_{\be}=\hat{a}v_{\be}\quad \textrm{ for } \be\in \G \textrm{ satisfying }\mathbb{F}\hat{\parti}_{1}+\mathbb{F}\hat{\parti}_{2}\not\subseteq \kn(\be+\mu_0),
\end{equation}
where $\hat{a}\in\mbb{F}$ is a nonzero constant. We then define
\begin{equation}\label{eq:80}
v_{\be}=\frac{1}{\hat{a}}\hat{w}_{\be}\quad \textrm{for }\be\in \G \textrm{ such that }D'\subseteq \kn(\be+\mu_0)\textrm{ and }(\be+\mu_0)(\bar{\parti})\not=0.
\end{equation}
Thus (\ref{eq:80}) together with (\ref{equ:4.222}) give the set
\begin{equation}
\{0\not=v_{\be}\in M_{\be}\mid \be\in \G \textrm{ satisfying }D''\not \subseteq \kn (\be+\mu_0)\}.
\end{equation}

{\it Fifthly}, we derive the action of $x^{\al}\parti$ on $v_{\be}$
for all $\parti \in \kn{\al}\cap D''$, $\al\in \G$ satisfying $D''\not \subseteq \kn{\al}$, and $\be \in \G$ satisfying $D''\not \subseteq \kn (\be+\mu_0)$ and $D''\not \subseteq \kn (\be+\mu_0+\al)$.\vspace{0.2cm}

{\it Case 1}. $\al\in \G$ satisfies $D'\subseteq \kn{\al}$ and $\al(\bar{\parti})\not=0$.\vspace{0.2cm}

{\it Subcase 1.1}. $\be\in \G$ satisfies $D'\subseteq \kn(\be+\mu_0)$, $(\be+\mu_0)(\bar{\parti})\not=0$ and $(\be+\mu_0+\al)(\bar{\parti})\not=0$.\vspace{0.2cm}

Pick any two linearly independent elements $\parti_{1}',\parti_{2}'\in\kn\al\cap D''=D'$. Set
\begin{equation}
\tilde{D}_{1}=\mathbb{F}\bar{\parti}+\mathbb{F}\parti_{1}'+\mathbb{F}\parti_{2}'.
\end{equation}
Then Lemma \ref{le:6.1}, Claim 2 and (\ref{equ:4.222}) imply that
there exist $\{0\not=w'_{\rho}\in M_{\rho}\mid \rho\in \G \textrm{
satisfying }\tilde{D}_{1}\not\subseteq \kn(\rho+\mu_0)\}$ such that
\begin{equation}\label{eq:81}
x^{\tau}\parti.  w'_{\rho} =(\rho+\mu_0)(\parti)w'_{\tau+\rho}
\end{equation}
for $\parti \in \kn{\tau}\cap \tilde{D}_{1}$, $\tau\in \G$ satisfying $\tilde{D}_1\not \subseteq \kn{\tau}$,
and $\rho \in \G$ satisfying $\tilde{D}_1\not \subseteq \kn (\rho+\mu_0)$ and $\tilde{D}_1\not \subseteq \kn
 (\rho+\mu_0+\tau)$. Since $\al(\bar{\parti})\not=0$, $(\be+\mu_0)(\bar{\parti})\not=0$ and $(\be+\mu_0+\al)
 (\bar{\parti})\not=0$, from (\ref{eq:81}) and $D'\subseteq \kn(\be+\mu_0)$, we see that
\begin{equation}
x^{\al}\parti.  w'_{\be} = (\be+\mu_0)(\parti)w'_{\be+\al}=0\ \textrm{ for } \parti \in \mathbb{F}\parti_{1}'
+\mathbb{F}\parti_{2}',
\end{equation}
which implies
\begin{equation}
x^{\al}\parti.  v_{\be} =0= (\be+\mu_0)(\parti)v_{\be+\al} \ \textrm{ for } \parti \in \mathbb{F}\parti_{1}'+\mathbb{F}\parti_{2}'.
\end{equation}
Since $\parti_{1}',\parti_{2}'\in\kn\al\cap D''$ are any two linearly independent elements, we see that
\begin{equation}
x^{\al}\parti.  v_{\be} =0= (\be+\mu_0)(\parti)v_{\be+\al} \  \textrm{ for }\parti\in\kn\al\cap D''.
\end{equation}

{\it Subcase 1.2}. $\be\in \G$ satisfies $D'\not\subseteq \kn(\be+\mu_0)$.\vspace{0.2cm}

Pick ${\parti}_{1}'\in D'\backslash \kn(\be+\mu_0)$. Moreover, we pick any ${\parti}_{2}' \in D'\backslash\mbb{F}\parti_{1}'$. Set
\begin{equation}
\tilde{D}_{2}=\mathbb{F}\bar{\parti}+\mathbb{F}{\parti}_{1}'+\mathbb{F}{\parti}_{2}'.
\end{equation}
From $(\be+\mu_0)(\parti_{1}')\not=0$, $D'\subseteq \kn{\al}$ and
$\al(\bar{\parti})\not=0$, we see that
$\mathbb{F}{\parti}_{1}'+\mathbb{F}{\parti}_{2}'\not\subseteq
\kn(\be+\mu_0)$, $\mathbb{F}{\parti}_{1}'
+\mathbb{F}{\parti}_{2}'\not\subseteq \kn(\be+\al+\mu_0)$ and
$\tilde{D}_{2}\not\subseteq \kn{\al}$. So Claim 1 and
Claim 2 show
\begin{equation}
x^{\al}\parti.  v_{\be} = (\be+\mu_0)(\parti)v_{\be+\al}\ \textrm{ for }\parti\in\kn\al\cap \tilde{D}_{2}=\mathbb{F}{\parti}_{1}'+\mathbb{F}{\parti}_{2}'.
\end{equation}
Since ${\parti}_{2}'\in D'\backslash\mbb{F}\parti_{1}'$ is arbitrary, we deduce
\begin{equation}\label{equ:4.238}
x^{\al}\parti.v_{\be}  =  (\be+\mu_0)(\parti)v_{\be+\al}\quad \textrm{ for } \parti \in D'=\kn\al\cap D''.
\end{equation}

{\it Case 2.} $\al\in \G$ satisfies $D'\not\subseteq \kn{\al}$.\vspace{0.2cm}

{\it Subcase 2.1}. $\be\in \G$ satisfies $D'\subseteq \kn(\be+\mu_0)$ and $(\be+\mu_0)(\bar{\parti})\not=0$.\vspace{0.2cm}

Similarly, in analogy with Subcase 1.1, we can prove
\begin{equation}\label{eq:94}
x^{\al}\parti.  v_{\be} =0= (\be+\mu_0)(\parti)v_{\be+\al} \  \textrm{ for }\parti\in\kn\al\cap D'.
\end{equation}

Pick $\parti_{1}' \in D'\backslash \kn{\al}$. Moreover, we pick any nonzero $\parti_{2}' \in \kn{\al}\cap (\mathbb{F}\hat{\parti}_{1}+\mathbb{F}\hat{\parti}_{2})$ (c.f. (\ref{eq:83})). Set
\begin{equation}
\tilde{D}_{3}=\mathbb{F}\bar{\parti}+\mathbb{F}{\parti}_{1}'+\mathbb{F}{\parti}_{2}'.
\end{equation}
Then Lemma \ref{le:6.1}, Claim 2 and (\ref{equ:4.222}) imply that
there exist $\{0\not=w'_{\rho}\in M_{\rho} \mid \rho\in \G \textrm{
satisfying }\tilde{D}_{3}\not\subseteq \kn(\rho+\mu_0)\}$ such that
\begin{equation}\label{eq:85}
x^{\tau}\parti.  w'_{\rho} =(\rho+\mu_0)(\parti)w'_{\tau+\rho}
\end{equation}
for $\parti \in \kn{\tau}\cap \tilde{D}_{3}$, $\tau\in \G$ satisfying $\tilde{D}_3\not \subseteq \kn{\tau}$, and $\rho \in \G$ satisfying $\tilde{D}_3\not \subseteq \kn (\rho+\mu_0)$ and $\tilde{D}_3\not \subseteq \kn (\rho+\mu_0+\tau)$. In particular,
\begin{equation}\label{eq:92}
x^{\al}(\bar{\parti}-\frac{\al(\bar{\parti})}{\al(\parti_{1}')}\parti_{1}').  w'_{\be} =(\be+\mu_0)(\bar{\parti}-\frac{\al(\bar{\parti})}{\al(\parti_{1}')}\parti_{1}')w'_{\be+\al}.
\end{equation}

Choose $\ga\in\G$ such that $\ga(\parti_{2}')\not=0$. Then (\ref{eq:85}) shows
\begin{equation}\label{eq:86}
x^{\ga}(\bar{\parti}-\frac{\ga(\bar{\parti})}{\ga(\parti_{2}')}\parti_{2}').  w'_{\be-\ga} =(\be+\mu_0)(\bar{\parti})w'_{\be}.
\end{equation}
On the other hand, as
$\bar{\parti}-\frac{\ga(\bar{\parti})}{\ga(\parti_{2}')}\parti_{2}'\in\hat{D}$,
we deduce from (\ref{eq:95}) that
\begin{equation}\label{eq:87}
x^{\ga}(\bar{\parti}-\frac{\ga(\bar{\parti})}{\ga(\parti_{2}')}\parti_{2}'). \hat{w}_{\be-\ga} =(\be+\mu_0)(\bar{\parti})\hat{w}_{\be}.
\end{equation}
Since $(\be-\ga+\mu_0)(\parti_{2}')=-\ga(\parti_{2}')\not=0$ and
$\parti_{2}'\in\mathbb{F}\hat{\parti}_{1}
+\mathbb{F}\hat{\parti}_{2}$,  we know from (\ref{eq:84}) that
$\hat{w}_{\be-\ga}=\hat{a}v_{\be-\ga}$. While (\ref{eq:80}),
$D'\subseteq \kn(\be+\mu_0)$ and $(\be+\mu_0)(\bar{\parti})\not=0$
imply $\hat{w}_{\be}=\hat{a}v_{\be}$. Inserting these results into
(\ref{eq:87}), we get
\begin{equation}\label{eq:88}
x^{\ga}(\bar{\parti}-\frac{\ga(\bar{\parti})}{\ga(\parti_{2}')}\parti_{2}'). v_{\be-\ga} =(\be+\mu_0)(\bar{\parti})v_{\be}.
\end{equation}
Since $(\al+\ga)(\parti_{2}')=\ga(\parti_{2}')\not=0$, (\ref{eq:85}) shows
\begin{equation}\label{eq:89}
x^{\al+\ga}(\parti_{1}'-\frac{(\al+\ga)(\parti_{1}')}{(\al+\ga)(\parti_{2}')}\parti_{2}').  w'_{\be-\ga} =\al(\parti_{1}')w'_{\be+\al}.
\end{equation}
As
$\parti_{1}'-\frac{(\al+\ga)(\parti_{1}')}{(\al+\ga)(\parti_{2}')}\parti_{2}'\in
D'$,  we derive from (\ref{equ:4.222}) that
\begin{equation}\label{eq:90}
x^{\al+\ga}(\parti_{1}'-\frac{(\al+\ga)(\parti_{1}')}{(\al+\ga)(\parti_{2}')}\parti_{2}').  v_{\be-\ga} =\al(\parti_{1}')v_{\be+\al}.
\end{equation}
Since $\al(\parti_{1}')\not=0$ and $(\be+\mu_0)(\bar{\parti})\not=0$, comparing (\ref{eq:90}) with (\ref{eq:89}) and comparing (\ref{eq:88}) with (\ref{eq:86}), we see that
\begin{equation}\label{eq:91}
w'_{\be-\ga}=cv_{\be-\ga},\ w'_{\be+\al}=cv_{\be+\al} \textrm{ and } w'_{\be}=cv_{\be} \textrm{ for some nonzero constant } c.
\end{equation}
So inserting (\ref{eq:91}) into (\ref{eq:92}), we get
\begin{equation}\label{eq:93}
x^{\al}(\bar{\parti}-\frac{\al(\bar{\parti})}{\al(\parti_{1}')}\parti_{1}').  v_{\be} =(\be+\mu_0)(\bar{\parti}-\frac{\al(\bar{\parti})}{\al(\parti_{1}')}\parti_{1}')v_{\be+\al}.
\end{equation}
Since $\kn\al\cap D''=\mbb{F}(\bar{\parti}-\frac{\al(\bar{\parti})}{\al(\parti_{1}')}\parti_{1}')+\kn\al\cap D'$, combining (\ref{eq:93}) with (\ref{eq:94}), we get
\begin{equation}
x^{\al}\parti.  v_{\be}= (\be+\mu_0)(\parti)v_{\be+\al} \  \textrm{ for }\parti\in\kn\al\cap D''.
\end{equation}

{\it Subcase 2.2}. $\be\in \G$ satisfies $D'\not\subseteq \kn(\be+\mu_0)$ and $D'\not\subseteq \kn(\be+\mu_0+\al)$.\vspace{0.2cm}

By (\ref{equ:4.222}) we have
\begin{equation}\label{equ:4.239}
x^{\al}\parti.v_{\be}  =  (\be+\mu_0)(\parti)v_{\al+\be}\ \textrm{ for } \parti \in \kn{\al} \cap D'.
\end{equation}
Choose ${\parti}_{1}'\in D'\backslash \kn{\al}$ and ${\parti}_{2}'\in D'\backslash (\kn(\be+\mu_0)\cup \kn(\be+\al+\mu_0))$ such that ${\parti}_{1}',  {\parti}_{2}'$ are linearly independent. Set
\begin{equation}
\tilde{D}_{4}=\mathbb{F}\bar{\parti}+\mathbb{F}{\parti}_{1}'+\mathbb{F}{\parti}_{2}'.
\end{equation}
Since $\mathbb{F}{\parti}_{1}'+\mathbb{F}{\parti}_{2}'\not\subseteq
\kn(\be+\mu_0)$ and
$\mathbb{F}{\parti}_{1}'+\mathbb{F}{\parti}_{2}'\not\subseteq
\kn(\be+\al+\mu_0)$, Claim 1 and Claim 2 imply
\begin{equation}\label{equ:4.240}
x^{\al}(\bar{\parti}-\frac{\al(\bar{\parti})}{\al(\parti_{1}')}\parti_{1}').v_{\be}  =
 (\be+\mu_0)(\bar{\parti}-\frac{\al(\bar{\parti})}{\al(\parti_{1}')}\parti_{1}')v_{\al+\be}.
\end{equation}
Since $\kn\al\cap D''=\mbb{F}(\bar{\parti}-\frac{\al(\bar{\parti})}{\al(\parti_{1}')}\parti_{1}')+\kn\al\cap D'$,
combining (\ref{equ:4.239}) with (\ref{equ:4.240}), we get
\begin{equation}\label{equ:4.241}
x^{\al}\parti.v_{\be}  =   (\be+\mu_0)(\parti)v_{\al+\be} \ \textrm{ for } \parti \in \kn{\al} \cap D''.
\end{equation}

{\it Subcase 2.3}. $\be\in \G$ satisfies $D'\not\subseteq \kn(\be+\mu_0)$, $D'\subseteq \kn(\be+\al+\mu_0)$ and $(\be+\al+\mu_0)(\bar{\parti})\not=0$.\vspace{0.2cm}

By the similar arguments as those in Subcase 1.1, we can prove
\begin{equation}\label{eq:96}
x^{\al}\parti.  v_{\be} =0= (\be+\mu_0)(\parti)v_{\be+\al} \ \ \textrm{ for }\parti\in\kn\al\cap D'.
\end{equation}

Pick $\parti_{1}' \in D'\backslash \kn{\al}$. Moreover, we pick any nonzero $\parti_{2}' \in \kn{\al}\cap (\mathbb{F}\hat{\parti}_{1}+\mathbb{F}\hat{\parti}_{2})$ (c.f. (\ref{eq:83})). Set
\begin{equation}
\tilde{D}_{5}=\mathbb{F}\bar{\parti}+\mathbb{F}{\parti}_{1}'+\mathbb{F}{\parti}_{2}'.
\end{equation}
Then Lemma \ref{le:6.1}, Claim 2 and (\ref{equ:4.222}) imply that
there exist $\{0\not=w'_{\rho}\in M_{\rho} \mid \rho\in \G \textrm{
satisfying }\tilde{D}_{5}\not\subseteq \kn(\rho+\mu_0)\}$ such that
\begin{equation}\label{eq:97}
x^{\tau}\parti.  w'_{\rho} =(\rho+\mu_0)(\parti)w'_{\tau+\rho}
\end{equation}
for $\parti \in \kn{\tau}\cap \tilde{D}_{5}$, $\tau\in \G$ satisfying $\tilde{D}_5\not \subseteq \kn{\tau}$, and $\rho \in \G$ satisfying $\tilde{D}_5\not \subseteq \kn (\rho+\mu_0)$ and $\tilde{D}_5\not \subseteq \kn (\rho+\mu_0+\tau)$. In particular, since $(\be+\mu_0)(\parti_{1}')=-\al(\parti_{1}')\not=0$ and $(\be+\al+\mu_0)(\bar{\parti})\not=0$, we have
\begin{equation}\label{eq:98}
x^{\al}(\bar{\parti}-\frac{\al(\bar{\parti})}{\al(\parti_{1}')}\parti_{1}').  w'_{\be} =(\be+\mu_0)(\bar{\parti}-\frac{\al(\bar{\parti})}{\al(\parti_{1}')}\parti_{1}')w'_{\be+\al}.
\end{equation}

Choose $\ga\in\G$ such that $\ga(\parti_{2}')\not=0$. Then (\ref{eq:97}) shows
\begin{equation}\label{eq:99}
x^{\ga}(\bar{\parti}-\frac{\ga(\bar{\parti})}{\ga(\parti_{2}')}\parti_{2}').  w'_{\be+\al-\ga} =(\be+\al+\mu_0)(\bar{\parti})w'_{\be+\al}.
\end{equation}
On the other hand, as
$\bar{\parti}-\frac{\ga(\bar{\parti})}{\ga(\parti_{2}')}\parti_{2}'\in\hat{D}$,
we deduce from (\ref{eq:95}) that
\begin{equation}\label{eq:100}
x^{\ga}(\bar{\parti}-\frac{\ga(\bar{\parti})}{\ga(\parti_{2}')}\parti_{2}'). \hat{w}_{\be+\al-\ga} =(\be+\al+\mu_0)(\bar{\parti})\hat{w}_{\be+\al}.
\end{equation}
Since $(\be+\al-\ga+\mu_0)(\parti_{2}')=-\ga(\parti_{2}')\not=0$ and
$\parti_{2}'\in\mathbb{F}\hat{\parti}_{1}
+\mathbb{F}\hat{\parti}_{2}$, we have
$\hat{w}_{\be+\al-\ga}=\hat{a}v_{\be+\al-\ga}$ by (\ref{eq:84}).
While (\ref{eq:80}), $D'\subseteq \kn(\be+\al+\mu_0)$ and
$(\be+\al+\mu_0)(\bar{\parti})\not=0$ imply
$\hat{w}_{\be+\al}=\hat{a}v_{\be+\al}$. Inserting these results into
(\ref{eq:100}), we get
\begin{equation}\label{eq:101}
x^{\ga}(\bar{\parti}-\frac{\ga(\bar{\parti})}{\ga(\parti_{2}')}\parti_{2}'). v_{\be+\al-\ga} =(\be+\al+\mu_0)(\bar{\parti})v_{\be+\al}.
\end{equation}
Since $(\ga-\al)(\parti_{2}')=\ga(\parti_{2}')\not=0$, (\ref{eq:97}) shows
\begin{equation}\label{eq:102}
x^{\ga-\al}(\parti_{1}'-\frac{(\ga-\al)(\parti_{1}')}{(\ga-\al)(\parti_{2}')}\parti_{2}').  w'_{\be+\al-\ga} =-\al(\parti_{1}')w'_{\be}.
\end{equation}
As
$\parti_{1}'-\frac{(\ga-\al)(\parti_{1}')}{(\ga-\al)(\parti_{2}')}\parti_{2}'\in
D'$,  we derive from (\ref{equ:4.222}) that
\begin{equation}\label{eq:103}
x^{\ga-\al}(\parti_{1}'-\frac{(\ga-\al)(\parti_{1}')}{(\ga-\al)(\parti_{2}')}\parti_{2}').  v_{\be+\al-\ga} =-\al(\parti_{1}')v_{\be}.
\end{equation}
Since $\al(\parti_{1}')\not=0$ and $(\be+\al+\mu_0)(\bar{\parti})\not=0$, comparing (\ref{eq:103}) with (\ref{eq:102}) and comparing (\ref{eq:101}) with (\ref{eq:99}), we see that
\begin{equation}\label{eq:104}
w'_{\be+\al-\ga}=c'v_{\be+\al-\ga},\ w'_{\be}=c'v_{\be} \textrm{ and } w'_{\be+\al}=c'v_{\be+\al} \textrm{ for some nonzero constant } c'.
\end{equation}
So inserting (\ref{eq:104}) into (\ref{eq:98}), we get
\begin{equation}\label{eq:105}
x^{\al}(\bar{\parti}-\frac{\al(\bar{\parti})}{\al(\parti_{1}')}\parti_{1}').  v_{\be} =(\be+\mu_0)(\bar{\parti}-\frac{\al(\bar{\parti})}{\al(\parti_{1}')}\parti_{1}')v_{\be+\al}.
\end{equation}
Since $\kn\al\cap D''=\mbb{F}(\bar{\parti}-\frac{\al(\bar{\parti})}{\al(\parti_{1}')}\parti_{1}')+\kn\al\cap D'$,
 combining (\ref{eq:105}) with (\ref{eq:96}), we find
\begin{equation}
x^{\al}\parti.  v_{\be} = (\be+\mu_0)(\parti)v_{\be+\al} \  \textrm{ for }\parti\in\kn\al\cap D''.
\end{equation}\vspace{0.1cm}

{\it To sum up}, the above two cases show
\begin{equation}\label{equ:4.242}
x^{\al}\parti.  v_{\be} =(\be+\mu_0)(\parti)v_{\al+\be}
\end{equation}
for $\parti \in \kn{\al}\cap D''$, $\al\in \G$ satisfying $D''\not \subseteq \kn{\al}$,
and $\be \in \G$ satisfying $D''\not \subseteq \kn (\be+\mu_0)$ and $D''\not \subseteq \kn (\be+\mu_0+\al)$. Namely, $D''$ satisfies $(p'2)$ in $(\Rmnum{2}')$. This contradicts the maximality of $D'$.

The case in which
$D'$ satisfies $(p'1)$ in $(\Rmnum{2}')$ similarly leads to a contradiction. We omit the details here. So we must have $D'=D$, from which
the lemma follows. $\qquad\Box$ \vspace{0.3cm}

\begin{lemma}\label{le:6.3} The ${\cal S}( \G,D )$-module
$M$ is isomorphic to
\begin{eqnarray}
&(\rmnum{1})&  \mathscr{M}_{\mu} \textrm{ for some }  \mu\in D^{\ast}; \\
&(\rmnum{2})&  \mathscr{A}_{\mu,\eta} \textrm{ for some }  \mu\in \G \textrm{ and } \eta\in D^{\ast}\backslash\{0\}; \\
&(\rmnum{3})&  \mathscr{B}_{\mu,\eta} \textrm{ for some }  \mu\in \G \textrm{ and } \eta\in D^{\ast}\backslash\{0\}; \\
&(\rmnum{4})&  \bigoplus_{\nu\in \G}\mathbb{F}w_{\nu}, \textrm{
where each component is a trivial submodule of }
{\cal S}( \G,D ).
   \end{eqnarray}
\end{lemma}

\noindent{\bf Proof.} If ${\cal S}(\G,D)$ act on $M$ as in $(P'1)$ of Lemma \ref{le:6.2}, i.e.,
\begin{equation}
x^{\al}\parti. M =\{0\} \textrm{ for }\parti \in \kn{\al} \textrm{ and }\al\in \G\backslash\{0\},
\end{equation}
we see that
\begin{equation}
M\simeq \bigoplus _{\nu \in \G}
\mathbb{F}w_{\nu}, \textrm{ where each component is a trivial submodule of }{\cal S} (\G, D).
\end{equation}

Assume that ${\cal S}(\G,D)$ act on $M$ as in $(P'2)$ of Lemma \ref{le:6.2}; i.e., there exist
$\mu \in D^{\ast}$ and $\{0\not=v_{\be}\in M_{\be}\mid \be\in \G\backslash\{-\mu\}\}$
such that
\begin{equation}\label{eq:106}
x^{\al}\parti.  v_{\be} =(\be+\mu)(\parti)v_{\al+\be}
\end{equation}
for $\parti \in \kn{\al}$, $\al\in \G\backslash\{0\}$ and $\be \in \G$ satisfying $\be+\mu\not=0\not=\be+\mu+\al$.

If $\mu\in D^{\ast}\backslash\G$,  (\ref{eq:106}) implies $M \simeq
\mathscr{M}_{\mu}$. Assume $\mu \in \G$. Pick some
$0\not=v_{-\mu}\in M_{-\mu}$. Write
\begin{equation}\label{equ:4.245}
x^{\al}\parti.  v_{-\mu} =g(\al, \parti)v_{\al-\mu} \textrm{ and }
x^{\al}\parti.  v_{-\mu-\al} =h(\al, \parti)v_{-\mu}
\end{equation}
for $\al\in \G\backslash\{0\}$ and $\parti \in \kn{\al}$. To complete the proof, we need to compute $g(\al, \parti)$ and $h(\al, \parti)$ for all $\al\in \G\backslash\{0\}$ and $\parti \in \kn{\al}$.
We now proceed our proof in three cases.\vspace{0.2cm}

{\it Case 1.} $g(\al,
\parti)=0$ and $h(\al, \parti)=0$ for all $\al\in \G\backslash\{0\}$ and
$\parti \in \kn{\al}$.\vspace{0.2cm}

By (\ref{eq:106}) and (\ref{equ:4.245}), we have $M \simeq
\mathscr{M}_{\mu}$. \vspace{0.2cm}

{\it Case 2.} $g(\al_{0}, \parti_{0})\not=0$ for some $\al_{0}\in \G\backslash\{0\}$ and $\parti_{0} \in \kn{\al_{0}}$.\vspace{0.2cm}

Applying (\ref{eq:106}) and (\ref{equ:4.245}) to
\begin{equation}\label{eq:109}
[x^{\al_{0}}\parti_{0},x^{\be}\parti].v_{-\mu-\be}=x^{\al_{0}+\be}(\be(\parti_{0})\parti-\al_{0}(\parti)\parti_{0}).v_{-\mu-\be},
\end{equation}
we get
\begin{equation}
g(\al_{0}, \parti_{0})h(\be, \parti)=0 \quad \textrm{ for } \be\in
\G\backslash\{0,\al_{0}\}, \;
\parti \in \kn{\be},
\end{equation}
which implies
\begin{equation}\label{equ:4.247}
h(\be, \parti)=0 \quad \textrm{ for } \be\in
\G\backslash\{0,\al_{0}\}, \;
\parti \in \kn{\be}.
\end{equation}

Fix any $\ga_1,\ga_2\in \G\backslash\{0\}$ such that
$\kn{\ga_1}\not=\kn{\ga_2}$. Pick any nonzero $\tilde{\parti} \in
\kn{\ga_1} \cap \kn{\ga_2}$, $\parti_{1} \in \kn{\ga_1} \backslash
\kn{\ga_2}$ and $\parti_{2} \in \kn{\ga_2} \backslash \kn{\ga_1}$.
Inserting (\ref{equ:4.245}) into
\begin{equation}
 [x^{-\ga_1}\parti_{1},[x^{\ga_2}\parti_{2},
x^{\ga_1}\tilde{\parti}]].v_{-\mu}=\ga_2(\parti_{1})\ga_1(\parti_{2})x^{\ga_2}\tilde{\parti}.v_{-\mu},
\end{equation}
we obtain
\begin{equation}\label{equ:8.19}
\ga_2(\parti_{1})\ga_1(\parti_{2})g(\ga_1, \tilde{\parti})-g(-\ga_1, \parti_1)h(\ga_1, \tilde{\parti})g(\ga_2, \parti_2)=\ga_2(\parti_{1})\ga_1(\parti_{2})g(\ga_2, \tilde{\parti}).
\end{equation}
When $\ga_1\not=\al_0$, using (\ref{equ:4.247}) in (\ref{equ:8.19}),
we get
\begin{equation}\label{eq:107}
g(\ga_2, \tilde{\parti})=g(\ga_1, \tilde{\parti}).
\end{equation}
When $\ga_1=\al_0$, applying (\ref{equ:4.245}) and (\ref{equ:4.247})
to
$x^{\al_0}\tilde{\parti}.x^{-\al_0}\parti_{1}.v_{-\mu}=x^{-\al_0}\parti_{1}.x^{\al_0}\tilde{\parti}.v_{-\mu}$,
we first have
\begin{equation}\label{equ:8.20}
g(-\al_0,\parti_1)h(\al_0,\tilde{\parti})=g(\al_0, \tilde{\parti})h(-\al_0, \parti_1)=0.
\end{equation}
Then inserting (\ref{equ:8.20}) into (\ref{equ:8.19}), we get
\begin{equation}\label{eq:108}
g(\ga_2, \tilde{\parti})=g(\al_0,
\tilde{\parti}).
\end{equation}
Since $\tilde{\parti} \in \kn{\ga_1} \cap \kn{\ga_2}$ is any nonzero
element, from (\ref{eq:107}) and (\ref{eq:108}), we see that
\begin{equation}\label{equ:4.248}
g(\ga_2, \parti)=g(\ga_1, \parti)
\end{equation}
for  $\parti \in
\kn{\ga_1} \cap \kn{\ga_2}$ and any $\ga_1, \ga_2\in \G\backslash\{0\}$
satisfying $\kn{\ga_1}\not=\kn{\ga_2}$.

Choose $\al_1\in \G\backslash\{0\}$ such that
$\kn{\al_1}\not=\kn{\al_{0}}$. Pick $\parti' \in
\kn{\al_1}\backslash \kn{\al_{0}}$. Since $\G=\kn{\al_{0}}+\mbb{F}\parti'$, we define $\eta \in D^{\ast}$ by
\begin{equation}\label{equ:4.249}
\eta(\parti)=g(\al_{0}, \parti) \textrm{ for any } \parti \in \kn{\al_{0}};\ \ \eta(\parti')=g(\al_1, \parti').
\end{equation}
Then by (\ref{equ:4.248}) and (\ref{equ:4.249}), we get
\begin{equation}\label{equ:4.250}
g(\al_1, \parti)=\eta(\parti) \  \textrm{ for any } \parti \in \kn{\al_1}.
\end{equation}
Choose $\al_{2}\in \G\backslash\{0\}$ such that $\textrm{codim}_{D}(\kn{\al_{0}}\cap \kn{\al_{1}}\cap \kn{\al_{2}})=3$. Then $\kn{\al_{2}}\not=\kn{\al_{1}}$, $\kn{\al_{2}}\not=\kn{\al_{0}}$ and $\kn{\al_{2}}=\kn{\al_{2}}\cap \kn{\al_{0}} + \kn{\al_{2}}\cap \kn{\al_{1}}$. By (\ref{equ:4.248}), (\ref{equ:4.249}) and (\ref{equ:4.250}), we derive
\begin{equation}\label{equ:4.251}
g(\al_{2}, \parti)=\eta(\parti)  \textrm{ for any } \parti \in \kn{\al_{2}}.
\end{equation}
For any $\be\in \G\backslash\{0\}$, we have $\textrm{codim}_{D}(\kn{\be}\cap \kn{\al_{i}}\cap \kn{\al_{j}})=3$ for some $i,  j \in \{0, 1,  2\}$ with $i\not=j$. So (\ref{equ:4.248})--(\ref{equ:4.251}) indicate,
\begin{equation}\label{equ:4.252}
g(\be, \parti)=\eta(\parti)  \ \textrm{ for any }  \be\in \G\backslash\{0\} \textrm{ and } \parti \in \kn{\be}.
\end{equation}

By (\ref{equ:4.252}) we have $g(-\al_{0}, \parti_{0})=\eta(\parti_{0})=g(\al_{0}, \parti_{0}) \not=0$. So, we get
\begin{equation}\label{equ:4.254}
h(\be, \parti)=0 \ \textrm{ for }\be \in \G\backslash \{0, -\al_0\}\textrm{ and } \parti \in \kn{\be}
\end{equation}
with $\al_0$ replaced by $-\al_0$ in the discussion from (\ref{eq:109}) to (\ref{equ:4.247}).
Hence, (\ref{eq:106}), (\ref{equ:4.245}), (\ref{equ:4.247}), (\ref{equ:4.252}) and
(\ref{equ:4.254}) imply $M \simeq \mathscr{A}_{\mu,\eta}$.\vspace{0.2cm}

{\it Case 3.} $g(\al, \parti)=0$ for all
$\al\in \G\backslash\{0\}$ and $\parti \in \kn{\al}$, while $h(\al_{0},
\parti_{0})\not=0$ for some $\al_{0}\in \G\backslash\{0\}$ and
$\parti_{0} \in \kn{\al_{0}}$. \vspace{0.2cm}

Fix any $\ga_1, \ga_2\in \G\backslash\{0\}$ such that
$\kn{\ga_1}\not=\kn{\ga_2}$. Pick any nonzero $\tilde{\parti} \in
\kn{\ga_1} \cap \kn{\ga_2}$, $\parti_{1} \in \kn{\ga_1} \backslash
\kn{\ga_2}$ and $\parti_{2} \in \kn{\ga_2} \backslash \kn{\ga_1}$.
Applying (\ref{eq:106}) and (\ref{equ:4.245}) to
\begin{equation}
[x^{-\ga_1}\parti_{1},[x^{\ga_2}\parti_{2},
x^{\ga_1}\tilde{\parti}]].v_{-\mu-\ga_2}=\ga_2(\parti_{1})\ga_1(\parti_{2})x^{\ga_2}\tilde{\parti}.v_{-\mu-\ga_2},
\end{equation}
we have
\begin{equation}
h(\ga_2, \tilde{\parti})=h(\ga_1, \tilde{\parti}).
\end{equation}
Since $\tilde{\parti} \in \kn{\ga_1} \cap
\kn{\ga_2}$ is any nonzero element, we get
\begin{equation}\label{equ:4.255}
h(\ga_2, \parti)=h(\ga_1, \parti)
\end{equation}
for $\parti \in \kn{\ga_1} \cap \kn{\ga_2}$ and $\ga_1, \ga_2\in \G\backslash\{0\}$ satisfying $ \kn{\ga_1}\not=\kn{\ga_2}$.

Choose $\al_1\in \G\backslash\{0\}$ such that $\kn{\al_1}\not=\kn{\al_{0}}$. Pick $\parti' \in \kn{\al_1}\backslash \kn{\al_{0}}$. We define $\eta \in D^{\ast}$ by
\begin{equation}\label{equ:8.21}
 \eta(\parti)=h(\al_{0}, \parti) \textrm{ for any } \parti \in \kn{\al_{0}};\ \ \eta(\parti')=h(\al_1, \parti').
\end{equation}
By the similar arguments as those from (\ref{equ:4.249}) to
(\ref{equ:4.252}), we can prove
\begin{equation}\label{equ:8.22}
h(\be, \parti)=\eta(\parti)\  \textrm{ for any }  \be\in \G\backslash\{0\} \textrm{ and } \parti \in \kn{\be}.
\end{equation}
Hence, from (\ref{eq:106}), (\ref{equ:4.245}) and
(\ref{equ:8.22}) we see that $M \simeq \mathscr{B}_{\mu,\eta}({\cal S}(\G,
D))$.

Thus, we complete the proof of the lemma. $\qquad\qquad\qquad \Box$

In summary, Lemma \ref{le:6.3}, together with Lemma \ref{le:3.1} and Theorem \ref{th:3.2}, implies

\begin{lemma}
If $\dim D \geq 4$, then Theorem \ref{th:4.1} holds.
\end{lemma}

\section*{Acknowledgment}
I would like to express my deep gratitude to Professor Xiaoping Xu
for all his advice, instructions and encouragements.

\end{document}